\documentclass[reqno, a4paper]{amsart}
\usepackage{amsmath}
\usepackage{amssymb}
\usepackage{amsthm}

\usepackage{thmtools} 

\usepackage[scale=0.9]{geometry}

\usepackage{natbib}
\usepackage{bibentry} 

\usepackage[english]{babel}
\usepackage[utf8]{inputenc}

\usepackage{subfig}
\usepackage{graphicx}

\usepackage[unicode, hypertexnames=false]{hyperref}
\usepackage[active]{srcltx} 

\usepackage{mathbbol}
\usepackage{bm} 
\usepackage{MnSymbol} 
\usepackage{gensymb}
\usepackage{eurosym}

\usepackage{units}
\usepackage{tensor}
\usepackage{accents}

\usepackage{enumitem}

\usepackage{lineno}

\usepackage{todonotes}

\usepackage{placeins}

\usepackage{xcolor}
\usepackage{mdframed}
\usepackage{newfloat}

\usepackage{tikz}
\usetikzlibrary{arrows.meta}

\DeclareMathOperator{\divergence}{div}

\DeclareMathOperator{\Tr}{Tr}

\DeclareMathOperator{\bregmandivergence}{BD} 

\newcommand{\smallo}[1]{\ensuremath{o\left(#1 \right)}}

\newcommand{\exponential}[1]{\ensuremath{{\mathrm e}^{#1}}}

\newcommand{\reference}{\mathrm{ref}}
\newcommand{\boundary}{\mathrm{bdr}}
\newcommand{\initial}{\mathrm{init}}

\newcommand{\bydefinition}{\mathrm{def}}
\newcommand{\traceless}[1]{{#1}_{\delta}}

\newcommand{\diff}{\mathrm{d}}
\newcommand{\Diff}[1][]{\mathrm{D}_{#1}} 
\newcommand{\hessian}[1]{\tensorq{H}_{#1}} 

\renewcommand{\vec}[1]{\ensuremath{\mathbf{#1}}}

\makeatletter
\@ifpackageloaded{bm}%
{\renewcommand{\vec}[1]{\ensuremath{\bm{#1}}}%
}{%
\relax
}
\makeatother

\newcommand{\tensorq}[1]{\ensuremath{\mathbb{#1}}}      

\newcommand{\conjugate}[1]{#1^\star} 
\newcommand{\transpose}[1]{#1^\top} 

\newcommand{\identity}{\ensuremath{\tensorq{I}}} 
\newcommand{\tensorzero}{{\mathbb{O}}} 

\newcommand{\cstress}{\tensorq{T}}

\makeatletter
\@ifpackageloaded{bm}%
{%
}{%
}

\@ifpackageloaded{bm}%
{%
 
}{%

}

\@ifpackageloaded{bm}%
{%
}{%
}
\makeatother

\newcommand{\generictensor}{{\tensorq{A}}}

\newcommand{\gradsym}{\ensuremath{\tensorq{D}}}

\newcommand{\R}{\ensuremath{{\mathbb R}}}

\newcommand{\ienergy}{\ensuremath{e}} 
\newcommand{\volienergy}{\overline{\ienergy}} 
\newcommand{\voltenergy}{\volienergy_{\mathrm{tot}}} 

\newcommand{\fenergy}{\ensuremath{\psi}} 

\newcommand{\entropy}{\ensuremath{\eta}} 
\newcommand{\volentropy}{\overline{\entropy}} 
\newcommand{\volentropygen}{\overline{\entropy}_{\mathrm{gen}}} 

\newcommand{\gibbs}{\ensuremath{g}} 

\newcommand{\temp}{\ensuremath{\theta}} 
\newcommand{\thpressure}{\ensuremath{p_{\mathrm{th}}}} 

\newcommand{\nettenergy}{\ensuremath{E}_{\mathrm{tot}}} 
\newcommand{\netkenergy}{\ensuremath{E_{\mathrm{kin}}}} 
\newcommand{\netentropy}{\ensuremath{S}} 

\newcommand{\ctentropy}{S}
\newcommand{\ctenergy}{U}
\newcommand{\ctvolume}{V}
\newcommand{\cttemperature}{T}
\newcommand{\ctpressure}{P}
\newcommand{\cthelmholtz}{F}

\newcommand{\cheatvol}{\ensuremath{c_{\mathrm{V}}}}
\newcommand{\cheatvolref}{\ensuremath{c_{\mathrm{V}\!,\,\reference}}} 

\newcommand{\hfluxc}{\vec{j}_{q}}     

\newcommand{\entfluxc}{\vec{j}_{\entropy}} 

\newcommand{\entprodc}{\xi} 

\newcommand{\pd}[2]{\ensuremath{\frac{\partial {#1}}{\partial {#2}}}}
\newcommand{\ppd}[2]{\ensuremath{\frac{\partial^2 {#1}}{\partial {#2^2}}}}
\newcommand{\dd}[2]{\ensuremath{\frac{\diff {#1}}{\diff {#2}}}}

\newcommand{\ddd}[2]{\ensuremath{\frac{\diff^2 {#1}}{\diff {#2}^2}}}

\makeatletter
\@ifpackageloaded{tensor}
{

}{%

}
\makeatother

\makeatletter
\@ifpackageloaded{tensor}
{

}{%

}
\makeatother

\newcommand{\norm}[2][]{\ensuremath{\left\|#2\right\|_{#1}}}
\newcommand{\absnorm}[1]{\ensuremath{\left|#1\right|}}

\makeatletter
\@ifundefined{volume}{%
}%
{%
}
\makeatother

\newcommand{\cvolumee}{{\diff} \mathrm{v}}

\newcommand{\csurfacees}{{\diff} \mathrm{s}}

\makeatletter

\@ifpackageloaded{MnSymbol} 
{
\newcommand{\tensordot}[2]{\ensuremath{#1 \vdotdot #2}} 
}{%
\newcommand{\tensordot}[2]{\ensuremath{#1 : #2}} 
}

\@ifpackageloaded{MnSymbol} 
{
   
}{%
   
}
\makeatother

\newcommand{\vectordot}[2]{\ensuremath{#1 \bullet #2}}

\newcommand{\sleb}[2]{\ensuremath{L}^{#1} \left(#2 \right)}             

\newcommand{\bregmandivergenceop}[3]{\bregmandivergence_{#1} \left(\left.\! #2 \right\| \! #3 \right)}

\newcommand{\tempbdr}{\temp_{\mathrm{bdr}}}

\newcommand{\tempref}{\ensuremath{\temp_{\reference}}}
\newcommand{\rhoref}{\ensuremath{\rho_{\reference}}}

\newcommand{\tempeq}{\ensuremath{\temp_{\equilibrium}}}
\newcommand{\rhoeq}{\ensuremath{\rho_{\equilibrium}}}

\newcommand{\tempsteady}{\temp_{\steady}}
\newcommand{\rhosteady}{\rho_{\steady}}

\newcommand{\lyapunovfn}[3]{\mathcal{V}_{#1} \left(\left.\! #2 \right\| \! #3 \right)}
\newcommand{\lyapunovfnup}[4]{\mathcal{V}_{#1}^{#2} \left(\left.\! #3 \right\| \! #4 \right)}

\newcommand{\equilibrium}{\mathrm{eq}}
\newcommand{\mequilibrium}{\mathrm{meq}}
\newcommand{\nequilibrium}{\mathrm{non-eq}}
\newcommand{\steady}{\mathrm{steady}}
\newcommand{\Bregman}{\mathrm{Bregman}}

\newcommand{\linmomentum}{\vec{m}}

\DeclareFloatingEnvironment[fileext=sum,placement={!ht},name=Summary]{summaryflt}

\newenvironment{summary}[1][]
    {
    \begin{summaryflt}[tb]
        \begin{summarycont}[#1] 
    }
    {
        \end{summarycont}
        \end{summaryflt}
    }

\mdfdefinestyle{summarystyle}{%
linecolor=black,linewidth=3pt,%
frametitlerule=true,%
frametitlebackgroundcolor=gray!20,
innertopmargin=\topskip,
}

\mdtheorem[style=summarystyle]{summarycont}{Summary}

\newmdenv[style=scholionstyle]{scholion}
\mdfdefinestyle{scholionstyle}{
  skipabove=1em,
  skipbelow=1em,
  backgroundcolor=gray!10,
  roundcorner=10pt
}

\declaretheorem[shaded]{theorem}
\declaretheorem[shaded,sibling=theorem]{definition}

\declaretheorem[shaded,sibling=theorem]{lemma}

\declaretheorem[sibling=theorem]{remark}

\numberwithin{equation}{section}

\makeatletter

\providecommand{\solutionname}{Solution}
\makeatother

\declaretheorem[shaded]{question}

\title[Thermodynamics and stability of (non)-equilibrium steady states]{Thermodynamics and stability of\\ equilibrium/non-equilibrium steady states in thermodynamically isolated/open systems\\---\\ case study for compressible heat conducting fluid}

\date{\today}

\author{V\'{\i}t Pr\r{u}\v{s}a}
\address{
Faculty of Mathematics and Physics\\
Charles University\\
Sokolovsk\'a 83\\
Praha 8 -- Karl\'{\i}n\\
CZ 186\;75\\
Czech Republic
}
\email{prusv@karlin.mff.cuni.cz}

\keywords{continuum thermodynamics, stability}
\subjclass[2000]{%
  76-02, 
  80-02
}

\begin{document}

\begin{abstract}
  We review all the calculations necessary for the construction of a Lyapunov like functional for nonlinear stability analysis of steady states in thermodynamically isolated/open systems wherein the substance of interest is a compressible heat conducting fluid.
\end{abstract}

\maketitle

\tableofcontents



\section{Introduction}
\label{sec:introduction}

\subsection{Compressible Navier--Stokes--Fourier equations}
\label{sec:compr-navi-stok}
Let us consider a compressible heat conducting fluid whose motion is described by the compressible Navier--Stokes--Fourier equations, that is by the equations
\begin{subequations}
  \label{eq:navier-stokes-fourier-equations}
  \begin{align}
    \label{eq:1}
    \dd{\rho}{t} + \rho \divergence \vec{v} &= 0, 
    \\
    \label{eq:2}
    \rho \dd{\vec{v}}{t}
                                            &=
                                              \divergence
                                              \left(
                                              -
                                              \thpressure(\temp, \rho) \identity
                                              +
                                              \tilde{\lambda} \left(\divergence \vec{v}\right) \identity
                                              +
                                              2 \nu \traceless{\gradsym}
                                              \right), 
    \\
    \label{eq:3}
    \rho \cheatvol(\temp, \rho)
    \dd{\temp}{t}
                                            &=
                                              -
                                              \temp
                                              \pd{\thpressure}{\temp}(\temp, \rho)
                                              \divergence \vec{v}
                                              +
                                              2 \nu \tensordot{\traceless{\gradsym}}{\traceless{\gradsym}}
                                              +
                                              \tilde{\lambda} \left( \divergence \vec{v}\right)^2
                                              +
                                              \divergence \left(\kappa \nabla \temp \right),
  \end{align}
  where $\rho$ denotes the density, $\temp$ denotes the thermodynamic temperature, $\thpressure$ denotes the thermodynamic pressure, $\cheatvol$ denotes the specific heat capacity at constant volume and $\vec{v}$ denotes the (Eulerian) velocity field, and $\gradsym$ denotes the symmetric part of the velocity gradient. The viscosities $\tilde{\lambda}$ and $\nu$ are positive constants; the thermal conductivity $\kappa$ is a positive constant as well. The symbol $\dd{}{t} =_{\bydefinition} \pd{}{t} + \vectordot{\vec{v}}{\nabla}$ denotes the material time derivative, and the symbol $\traceless{\generictensor} =_{\bydefinition} \generictensor - \frac{1}{3} \left( \Tr \generictensor \right) \identity$ denotes the traceless/deviatoric part of the corresponding tensor. 

  The thermodynamic pressure $\thpressure$ and the specific heat capacity at constant volume $\cheatvol$, as functions of the temperature and density, are given in terms of the Helmholtz free energy function $\fenergy = \fenergy(\temp, \rho)$ as
  \begin{align}
    \label{eq:4}
    \cheatvol &= - \temp  \ppd{\fenergy}{\temp}, 
    \\
    \label{eq:5}
    \thpressure &= \rho^2 \pd{\fenergy}{\rho}.
\end{align}
\end{subequations}
For further reference we recall that other thermodynamic potentials are given by the formulae
\begin{equation}
  \label{eq:205}
  \entropy (\temp, \rho) = - \pd{\fenergy}{\temp}(\temp, \rho)
\end{equation}
and
\begin{equation}
  \label{eq:289}
  \ienergy (\temp, \rho) = \fenergy(\temp, \rho) + \temp \entropy(\temp, \rho).  
\end{equation}
Note however that these formulae do not give the potentials in terms of their natural variables but in terms of the temperature and the density, which are the variables used in practice. Finally, we recall that the quantity $\hfluxc = - \kappa \nabla \temp$ is referred to as the heat flux vector.

\subsection{Calorically perfect ideal gas}
\label{sec:calor-perf-ideal}
In particular, for the \emph{calorically perfect ideal gas} the Helmholtz free energy takes the form
\begin{equation}
  \label{eq:free-energy-ideal-gas-introduction}
  \fenergy(\temp, \rho)
  =
  -
  \cheatvolref
  \temp
  \left(
    \ln
    \left(
      \frac{\temp}{\tempref}
    \right)
    -
    1
  \right)
  +
  \cheatvolref
  \temp
  \left(
    \gamma -1
  \right)
  \ln
  \left(
    \frac{\rho}{\rhoref}
  \right)
  ,
\end{equation}
where $\cheatvolref$, $\tempref$, $\rhoref$ and $\gamma$ are constants. All these constants are positive and the adiabatic exponent $\gamma$ is greater than one,~$\gamma > 1$.

Furthermore, the constants $\temp_\reference$ and $\rhoref$ serve just for normalisation purposes---we can fix them at will as the choice of normalisation constants has no impact on the final formula for the specific heat capacity at constant volume or the formula for the thermodynamic pressure. In other words, the formulae for the specific heat capacity at constant volume and the thermodynamic pressure, which \emph{have} an actual impact on the dynamical behaviour of the corresponding substance, remain the same for different choices of normalisation constants. Indeed, if we explicitly evaluate the formulae for the specific heat capacity at constant volume~\eqref{eq:4} and the thermodynamic pressure~\eqref{eq:5}, then we get
\begin{subequations}
  \label{eq:7}
  \begin{align}
    \label{eq:8}
    \cheatvol &= \cheatvolref, 
    \\
    \label{eq:9}
    \thpressure &= \cheatvolref \left(\gamma -1 \right) \temp \rho,
  \end{align}
\end{subequations}
and the formulae are free from normalisation constants.

The entropy and the internal energy are---in this special case---given as
\begin{subequations}
  \label{eq:10}
  \begin{align}
    \label{eq:11}
    \entropy (\temp, \rho)
    &=
    \cheatvolref
    \ln
    \left[
    \frac{\temp}{\tempref}
    \left(
    \frac{\rho}{\rhoref}
    \right)^{1 - \gamma}
    \right]
      ,
    \\
    \label{eq:12}
    \ienergy (\temp, \rho) &= \cheatvolref \temp,
  \end{align}
\end{subequations}
while in the natural variables we get
\begin{subequations}
  \label{eq:13}
  \begin{align}
    \label{eq:14}
    \ienergy(\entropy, \rho)
    &=
      \cheatvolref
      \tempref
      \left(
      \frac{\rho}{\rhoref}
      \right)^{\gamma-1}
      \exponential{\frac{\entropy}{\cheatvolref}}
      , \\
    \label{eq:15}
    \entropy(\ienergy, \rho)
    &=
    \cheatvolref
    \ln
    \left(
      \frac{\ienergy}{\cheatvolref \tempref}
      \left(
        \frac{\rho}{\rhoref}
      \right)^{1 - \gamma}
    \right)
    .
  \end{align}
\end{subequations}

\subsection{Questions regarding stability of equilibrium/non-equilibrium steady states in thermodynamically isolated/open systems}
\label{sec:quest-regard-stab}
The first question we are interested in is the following.

\begin{question}[Stability of spatially homogeneous rest state in a thermodynamically isolated system]
  \label{q:isolated}
  Consider a heat conducting fluid described by the Navier--Stokes--Fourier equations~\eqref{eq:navier-stokes-fourier-equations}
  \begin{subequations}
    \label{eq:navier-stokes-fourier-question-one}
    \begin{align}
      \label{eq:16}
      \dd{\rho}{t} + \rho \divergence \vec{v} &= 0, 
      \\
      \label{eq:17}
      \rho \dd{\vec{v}}{t}
                                              &=
                                                \divergence
                                                \left(
                                                -
                                                \thpressure(\temp, \rho) \identity
                                                +
                                                \tilde{\lambda} \left( \divergence \vec{v} \right) \identity
                                                +
                                                2 \nu \traceless{\gradsym}
                                                \right), 
      \\
      \label{eq:18}
      \rho \cheatvol(\temp, \rho)
      \dd{\temp}{t}
                                              &=
                                                -
                                                \temp
                                                \pd{\thpressure}{\temp}(\temp, \rho)
                                                \divergence \vec{v}
                                                +
                                                2 \nu \tensordot{\traceless{\gradsym}}{\traceless{\gradsym}}
                                                +
                                                \tilde{\lambda} \left( \divergence \vec{v}\right)^2
                                                +
                                                \divergence \left(\kappa \nabla \temp \right),
    \end{align}
  \end{subequations}
  where the material functions $\thpressure$ and $\cheatvol$ are derived from a general Helmholtz free energy $\fenergy$, see~\eqref{eq:4} and \eqref{eq:5}. Assume that the fluid occupies a thermodynamically isolated container $\Omega$, which means that the boundary conditions read
  \begin{subequations}
    \label{eq:boundary-conditions-isolated}
    \begin{align}
      \label{eq:19}
      \left. \vec{v} \right|_{\partial \Omega} &= \vec{0}, \\
      \label{eq:20}
      \left. \kappa \vectordot{\nabla \temp}{\vec{n}} \right|_{\partial \Omega} &= 0.
    \end{align}
  \end{subequations}
  The first condition guarantees no mechanical energy exchange with the surroundings, while the second condition guarantees no heat exchange with the surroundings.

  Let us now consider a spatially homogeneous rest state
  \begin{equation}
    \label{eq:spatially-homogenoeus-rest-state}
    \begin{bmatrix}
      \widehat{\rho} \\
      \widehat{\vec{v}} \\
      \widehat{\temp}
    \end{bmatrix}
    =_{\bydefinition}
    \begin{bmatrix}
      \rhoeq \\
      \vec{0} \\
      \tempeq
    \end{bmatrix}
    ,
  \end{equation}
  where $\rhoeq$ and $\tempeq$ are constants both in space and time. The spatially homogeneous rest state is clearly a solution to~\eqref{eq:navier-stokes-fourier-question-one} with boundary conditions~\eqref{eq:boundary-conditions-isolated}.

  Let us further consider the solution to the initial--boundary value problem~\eqref{eq:navier-stokes-fourier-question-one} with the initial condition
  \begin{equation}
    \label{eq:inital-condition-isolated}
    \left.
      \begin{bmatrix}
        \rho \\
        \vec{v} \\
        \temp
      \end{bmatrix}
    \right|_{t=0}
    =
    \begin{bmatrix}
        \rho_{\initial} \\
        \vec{v}_{\initial} \\
        \temp_{\initial}
    \end{bmatrix},
  \end{equation}
  such that the net total energy/net mass of the initial state~\eqref{eq:inital-condition-isolated} is the same as the net total energy/net mass of the spatially homogeneous rest state~\eqref{eq:spatially-homogenoeus-rest-state} meaning that
  \begin{subequations}
    \label{eq:23}
    \begin{align}
      \label{eq:24}
      \int_{\Omega} \left( \frac{1}{2} \rho_{\initial} \absnorm{\vec{v}_{\initial}}^2 + \rho_{\initial} \ienergy \left( \temp_{\initial}, \rho_{\initial} \right) \right) \, \cvolumee
      &=
        \int_{\Omega}  \rhoeq \ienergy \left( \tempeq, \rhoeq \right)  \, \cvolumee, \\
      \label{eq:25}
      \int_{\Omega}  \rho_{\initial} \, \cvolumee
      &=
        \int_{\Omega}  \rhoeq \, \cvolumee.
    \end{align}
  \end{subequations}

  Based on everyday experience we expect that the solution to the initial-boundary value problem~\eqref{eq:navier-stokes-fourier-question-one}, \eqref{eq:boundary-conditions-isolated}, \eqref{eq:inital-condition-isolated} with \emph{arbitrary initial data} \eqref{eq:inital-condition-isolated} converges, in some sense, to the spatially homogeneous rest state~\eqref{eq:spatially-homogenoeus-rest-state}, that is
  \begin{equation}
    \label{eq:convergenece-isolated}
    \begin{bmatrix}
      \rho \\
      \vec{v} \\
      \temp
    \end{bmatrix}
    \stackrel{t \to +\infty}{\to}
    \begin{bmatrix}
      \rhoeq \\
      \vec{0} \\
      \tempeq
    \end{bmatrix}
    .
  \end{equation}
  Can we show this using the Navier--Stokes--Fourier equations? 
\end{question}
  
Clearly, if we are not able to answer this question, then we are doomed---our laboriously built mathematical apparatus fails to recover elementary qualitative behaviour of heat conducting compressible fluids. Moreover, we would like to make use of \emph{some fancy thermodynamics in answering this question.} (Otherwise what would thermodynamics be good for?) If we manage to answer the first question, we can proceed with a more difficult one. In particular, we can ask the same question for a primitive thermodynamically \emph{open} system. The question is the following.

\begin{question}[Stability of spatially inhomogeneous steady state in a primitive thermodynamically open system]
  \label{q:open}
  Consider a heat conducting fluid described by the Navier--Stokes--Fourier equations~\eqref{eq:navier-stokes-fourier-equations}
  \begin{subequations}
    \label{eq:navier-stokes-fourier-question-two}
    \begin{align}
      \label{eq:27}
      \dd{\rho}{t} + \rho \divergence \vec{v} &= 0, 
      \\
      \label{eq:28}
      \rho \dd{\vec{v}}{t}
                                              &=
                                                \divergence
                                                \left(
                                                -
                                                \thpressure(\temp, \rho) \identity
                                                +
                                                \tilde{\lambda} \left( \divergence \vec{v} \right) \identity
                                                +
                                                2 \nu \traceless{\gradsym}
                                                \right), 
      \\
      \label{eq:29}
      \rho \cheatvol(\temp, \rho)
      \dd{\temp}{t}
                                              &=
                                                -
                                                \temp
                                                \pd{\thpressure}{\temp}(\temp, \rho)
                                                \divergence \vec{v}
                                                +
                                                2 \nu \tensordot{\traceless{\gradsym}}{\traceless{\gradsym}}
                                                +
                                                \tilde{\lambda} \left( \divergence \vec{v}\right)^2
                                                +
                                                \divergence \left(\kappa \nabla \temp \right),
    \end{align}
  \end{subequations}
  where the material functions $\thpressure$ and $\cheatvol$ are derived from  a general Helmholtz free energy $\fenergy$, see~\eqref{eq:4} and \eqref{eq:5}. Assume that the fluid occupies a mechanically isolated container $\Omega$, which means that the boundary conditions read
  \begin{subequations}
    \label{eq:boundary-conditions-dirichlet-temperature}
    \begin{align}
      \label{eq:30}
      \left. \vec{v} \right|_{\partial \Omega} &= \vec{0}, \\
      \label{eq:31}
      \left. \temp \right|_{\partial \Omega} &= \tempbdr.
    \end{align}
  \end{subequations}
  The first condition guarantees no mechanical energy exchange with the surroundings, while the second condition allows heat exchange with the surroundings---the temperature gradient, that is the heat flux, is not necessarily zero on the boundary. (Compare with~\eqref{eq:boundary-conditions-isolated}.) The boundary temperature value $\temp_{\boundary}$ might be a \emph{function of position}.

  Let $\tempsteady$ now denote the spatially inhomogeneous steady state which solves the boundary value problem 
  \begin{subequations}
    \label{eq:32}
    \begin{align}
      \label{eq:33}
      \divergence \left( \kappa \nabla \tempsteady \right) &= 0, \\
      \label{eq:34}
      \left. \tempsteady \right|_{\partial \Omega} &= \tempbdr,
    \end{align}
  \end{subequations}
  for the steady heat equation. Let
  \begin{equation}
    \label{eq:35}
    m_{\Omega} =_{\bydefinition} \int_{\Omega} \rho \, \cvolumee
  \end{equation}
  denote the net mass of fluid inside the container. (Thanks to the boundary condition~\eqref{eq:30} the total mass is conserved.) Furthermore, let $\rhosteady$ be a density field that leads to spatially constant pressure field, that is the density field obtained by the solution of pointwise algebraic equations
  \begin{equation}
    \label{eq:36}
    \thpressure(\tempsteady, \rhosteady) = C,
  \end{equation}
  for the given constant $C$ and the spatially dependent temperature field $\tempsteady$. (The constant $C$ is chosen such that the so obtained density field~$\rhosteady$ satisfies $\int_{\Omega} \rhosteady \, \cvolumee = m_{\Omega}$.) The density/velocity/temperature field
  \begin{equation}
    \label{eq:spatially-inhomogenoeus-steady-state-dirichlet}
    \begin{bmatrix}
      \widehat{\rho} \\
      \widehat{\vec{v}} \\
      \widehat{\temp}
    \end{bmatrix}
    =_{\bydefinition}
    \begin{bmatrix}
      \rhosteady \\
      \vec{0} \\
      \tempsteady
    \end{bmatrix}
    ,
  \end{equation}
  is clearly a solution to~\eqref{eq:navier-stokes-fourier-question-two} with boundary conditions~\eqref{eq:boundary-conditions-dirichlet-temperature}.

  Let us further consider the solution to the initial--boundary value problem~\eqref{eq:navier-stokes-fourier-question-two} with the initial condition
  \begin{equation}
    \label{eq:initial-condition-dirichlet}
    \left.
      \begin{bmatrix}
        \rho \\
        \vec{v} \\
        \temp
      \end{bmatrix}
    \right|_{t=0}
    =
    \begin{bmatrix}
        \rho_{\initial} \\
        \vec{v}_{\initial} \\
        \temp_{\initial}
    \end{bmatrix},
  \end{equation}
  and let us assume that the net mass of the initial state~\eqref{eq:initial-condition-dirichlet} is the same as the net mass of the spatially inhomogeneous steady state~\eqref{eq:spatially-inhomogenoeus-steady-state-dirichlet} meaning that
  \begin{equation}
    \label{eq:39}
    \int_{\Omega}  \rho_{\initial} \, \cvolumee
      =
        \int_{\Omega}  \rhosteady \, \cvolumee
        .
  \end{equation}
  Based on everyday experience we expect that the solution to the initial-boundary value problem~\eqref{eq:navier-stokes-fourier-question-two}, \eqref{eq:boundary-conditions-dirichlet-temperature}, \eqref{eq:initial-condition-dirichlet} with \emph{arbitrary initial data} \eqref{eq:initial-condition-dirichlet} converges, in some sense, to the spatially inhomogeneous steady state, that is
  \begin{equation}
    \label{eq:convergence-dirichlet}
    \begin{bmatrix}
      \rho \\
      \vec{v} \\
      \temp
    \end{bmatrix}
    \stackrel{t \to +\infty}{\to}
    \begin{bmatrix}
      \rhosteady \\
      \vec{0} \\
      \tempsteady
    \end{bmatrix}
    .
  \end{equation}
  Can we show this using the Navier--Stokes--Fourier equations?
\end{question}

We can continue phrasing similar questions for more and more complicated systems. In particular we can add a body force (gravitational force) to the right-hand side of balance of linear momentum~\eqref{eq:17}, which would lead to the thermal convection problem (Rayleigh--B\'enard convection for compressible heat conducting fluid). However, for the time being, we stick to Question~\ref{q:isolated} and Question~\ref{q:open}.

Before we proceed with analysis, we make a remark on the level of rigour. All calculations in the forthcoming sections are \emph{formal}. We assume that the solution to the corresponding initial--boundary value problem exists for all times, that it is sufficiently smooth so that all differentiations and integrations by parts are justified, and that the density and the temperature remain positive, that is $\rho > 0$ and $\temp > 0$ in $\Omega$ for all times. Consequently, our analysis does not provide a rigorous answer to Question~\ref{q:isolated} and Question~\ref{q:open}. It rather identifies the functionals that are, from the thermodynamic point of view, the natural candidates for the rigorous analysis. In fact, a functional of the same type---the relative energy---is the key tool in the rigorous analysis of weak solutions, see the discussion in Section~\ref{sec:funct-used-feir}.



\section{Isolated systems---spatially homogeneous rest state}
\label{sec:isol-syst-spat}
We start with Question~\ref{q:isolated}, and we \emph{first investigate what can be said if we assume that the perturbations from the spatially homogeneous equilibrium steady state are small}, that is if we describe the evolution close to the spatially homogeneous steady state by the \emph{linearised governing equations}. Clearly, if we do not succeed in the linearised setting, we can hardly expect that the nonlinearity will save us. This is the reason why we investigate the linearised equations.

\subsection{Linearised setting}
\label{sec:linearised-setting}

We linearise the Navier--Stokes--Fourier equations~\eqref{eq:navier-stokes-fourier-equations} in the neighbourhood of the spatially homogeneous rest state~\eqref{eq:spatially-homogenoeus-rest-state}. We rewrite the density, temperature and velocity as
\begin{subequations}
  \label{eq:41}
\begin{align}
  \label{eq:42}
  \rho &= \widehat{\rho} + \widetilde{\rho}, \\
  \label{eq:43}
  \vec{v} &= \widehat{\vec{v}} + \widetilde{\vec{v}}, \\
  \label{eq:44}
  \temp &= \widehat{\temp} + \widetilde{\temp},
\end{align}
\end{subequations}
and we want to find the leading order equations for the perturbation $[\widetilde{\rho}, \widetilde{\vec{v}}, \widetilde{\temp}]$.

\subsubsection{Linearised governing equations}
\label{sec:line-govern-equat}

Let us start with the balance of mass~\eqref{eq:1}. We see that
\begin{multline}
  \label{eq:45}
  \dd{\rho}{t} + \rho \divergence{\vec{v}}
  =
  \pd{\rho}{t}
  +
  \vectordot{
    \vec{v}
  }
  {
    \nabla \rho
  }
  +
  \rho
  \divergence \vec{v}
  =
  \pd{}{t}
  \left(
    \widehat{\rho} + \widetilde{\rho}
  \right)
  +
  \vectordot{
    \widehat{\vec{v}} + \widetilde{\vec{v}}
  }
  {
    \nabla \left(\widehat{\rho} + \widetilde{\rho}\right)
  }
  +
  \left(\widehat{\rho} + \widetilde{\rho}\right)
  \divergence \left( \widehat{\vec{v}} + \widetilde{\vec{v}} \right)
  \\
  =
  \pd{\widetilde{\rho}}{t}
  +
  \vectordot{
    \widetilde{\vec{v}}
  }
  {
    \nabla \widetilde{\rho}
  }
  +
  \left(\widehat{\rho} + \widetilde{\rho}\right)
  \divergence \widetilde{\vec{v}}
  \approx
  \pd{\widetilde{\rho}}{t}
  +
  \widehat{\rho}
  \divergence \widetilde{\vec{v}}
  ,
\end{multline}
where in the last equality we have neglected the terms that are nonlinear in the perturbation. (Recall that $\widehat{\rho}$ is a constant in space and time, and that $\widehat{\vec{v}}=\vec{0}$.) Consequently, the linearisation of the balance of mass~\eqref{eq:1} in the neighbourhood of the spatially homogeneous rest state reads
\begin{equation}
  \label{eq:46}
  \pd{\widetilde{\rho}}{t}
  +
  \widehat{\rho}
  \divergence \widetilde{\vec{v}}
  =
  0
  .
\end{equation}
Now we linearise the balance of linear momentum~\eqref{eq:2}. We use the Taylor expansion of the formula for the pressure,
\begin{equation}
  \label{eq:47}
  \thpressure(\temp, \rho)
  \approx
  \thpressure(\widehat{\temp}, \widehat{\rho})
  +
  \left. \pd{\thpressure}{\temp}(\temp, \rho) \right|_{\temp = \widehat{\temp}, \, \rho = \widehat{\rho}} \widetilde{\temp}
  +
  \left. \pd{\thpressure}{\rho}(\temp, \rho) \right|_{\temp = \widehat{\temp}, \, \rho = \widehat{\rho}} \widetilde{\rho}
  ,
\end{equation}
and we see that the linearised version of the balance of linear momentum reads
\begin{equation}
  \label{eq:48}
  \widehat{\rho} \pd{\widetilde{\vec{v}}}{t}
  =
  \divergence
  \left(
    -
    \left(
      \left. \pd{\thpressure}{\temp}(\temp, \rho) \right|_{\temp = \widehat{\temp}, \, \rho = \widehat{\rho}} \widetilde{\temp}
      +
      \left. \pd{\thpressure}{\rho}(\temp, \rho) \right|_{\temp = \widehat{\temp}, \, \rho = \widehat{\rho}} \widetilde{\rho}
    \right)
    \identity
    +
    \tilde{\lambda} \left( \divergence \widetilde{\vec{v}} \right) \identity
    +
    2 \nu \traceless{\widetilde{\gradsym}}
  \right), 
\end{equation}
where $\widetilde{\gradsym} =_{\bydefinition} \frac{1}{2} \left( \nabla \widetilde{\vec{v}} + \transpose{\nabla \widetilde{\vec{v}}} \right)$, and where we have used the fact that $\thpressure(\widehat{\temp}, \widehat{\rho})$ is a constant. Finally, we linearise the temperature evolution equation~\eqref{eq:3}, and we arrive at
\begin{equation}
  \label{eq:49}
  \widehat{\rho} \cheatvol\left(\widehat{\temp}, \widehat{\rho}\right) 
  \pd{\widetilde{\temp}}{t}
  =
  -
  \widehat{\temp}
  \left. \pd{\thpressure}{\temp}(\temp, \rho) \right|_{\temp = \widehat{\temp}, \, \rho = \widehat{\rho}}
  \divergence \widetilde{\vec{v}}
  +
  \divergence \left(\kappa \nabla \widetilde{\temp} \right)
  .
\end{equation}

Consequently, the linearised system of governing equations for the perturbation reads
\begin{subequations}
  \label{eq:navier-stokes-fourier-perturbation-linearised}
  \begin{align}
    \label{eq:50}
    \pd{\widetilde{\rho}}{t}
    +
    \widehat{\rho}
    \divergence \widetilde{\vec{v}}
    &=
      0
      ,
    \\
    \label{eq:51}
    \widehat{\rho} \pd{\widetilde{\vec{v}}}{t}
    &=
      \divergence
      \left(
      -
      \left(
      \left. \pd{\thpressure}{\temp}(\temp, \rho) \right|_{\temp = \widehat{\temp}, \, \rho = \widehat{\rho}} \widetilde{\temp}
      +
      \left. \pd{\thpressure}{\rho}(\temp, \rho) \right|_{\temp = \widehat{\temp}, \, \rho = \widehat{\rho}} \widetilde{\rho}
      \right)
      \identity
      +
      \tilde{\lambda} \left( \divergence \widetilde{\vec{v}} \right) \identity
      +
      2 \nu \traceless{\widetilde{\gradsym}}
      \right), 
    \\
    \label{eq:52}
    \widehat{\rho} \cheatvol\left(\widehat{\temp}, \widehat{\rho}\right) 
    \pd{\widetilde{\temp}}{t}
    &=
      -
      \widehat{\temp}
      \left. \pd{\thpressure}{\temp}(\temp, \rho) \right|_{\temp = \widehat{\temp}, \, \rho = \widehat{\rho}}
      \divergence \widetilde{\vec{v}}
      +
      \divergence \left(\kappa \nabla \widetilde{\temp} \right)
      .
  \end{align}
\end{subequations}

\subsubsection{Decay of perturbations in linearised setting}
\label{sec:decay-pert-line}
Using the system of linearised equations we would like to conclude that perturbations to the spatially homogeneous rest state decay, that is we want to show that
\begin{equation}
  \label{eq:53}
  \begin{bmatrix}
    \widetilde{\rho} \\
    \widetilde{\vec{v}} \\
    \widetilde{\temp}
  \end{bmatrix}
  \stackrel{t \to +\infty}{\to}
  \begin{bmatrix}
    0 \\
    \vec{0} \\
    0
  \end{bmatrix}
\end{equation}
for any (small) initial condition/perturbation. The meaning of the ``convergence'' in~\eqref{eq:53} is a matter for us to consider/specify; naturally, we stick to the notion of convergence that emerges as natural in the forthcoming analysis.  

We manipulate the linearised governing equations~\eqref{eq:navier-stokes-fourier-perturbation-linearised} as follows. First, we take the scalar product of~\eqref{eq:51} with $\widetilde{\vec{v}}$, we integrate over the domain $\Omega$, and then we use the integration by parts. The boundary condition~\eqref{eq:boundary-conditions-isolated} implies that $\left. \widetilde{\vec{v}} \right|_{\partial \Omega} = \vec{0}$, hence the boundary term in the integration by parts vanishes, and we get
\begin{equation}
  \label{eq:54}
  \int_{\Omega}\widehat{\rho} \frac{1}{2}\pd{}{t}\absnorm{\widetilde{\vec{v}}} ^2 \, \cvolumee
  =
  \int_{\Omega}
  \left. \pd{\thpressure}{\temp}(\temp, \rho) \right|_{\temp = \widehat{\temp}, \, \rho = \widehat{\rho}} \widetilde{\temp}
  \divergence{\widetilde{\vec{v}}}
  \, \cvolumee
  +
  \int_{\Omega}
  \left. \pd{\thpressure}{\rho}(\temp, \rho) \right|_{\temp = \widehat{\temp}, \, \rho = \widehat{\rho}} \widetilde{\rho}
  \divergence{\widetilde{\vec{v}}}
  \, \cvolumee
  -
  \int_{\Omega}
  \tilde{\lambda} \left(\divergence \widetilde{\vec{v}}\right)^2
  \, \cvolumee
  -
  \int_{\Omega}
  2 \nu \tensordot{\traceless{\widetilde{\gradsym}}}{\traceless{\widetilde{\gradsym}}}
  \, \cvolumee
  .
\end{equation}
The last two terms are negative for nonzero perturbation, which is convenient---these are the terms that force the time derivative on the left-hand side to be negative. It would be nice if we managed to cancel the first two terms without a definite sign. This is possible.

We divide the linearised balance of mass~\eqref{eq:50} by $\widehat{\rho}$, we multiply the equation by $\left. \pd{\thpressure}{\rho}(\temp, \rho) \right|_{\temp = \widehat{\temp}, \, \rho = \widehat{\rho}} \widetilde{\rho}$, and we integrate over the domain $\Omega$. We get
\begin{equation}
  \label{eq:55}
  \frac{1}{\widehat{\rho}}\left. \pd{\thpressure}{\rho}(\temp, \rho) \right|_{\temp = \widehat{\temp}, \, \rho = \widehat{\rho}}
  \int_{\Omega} \frac{1}{2}\pd{}{t} \widetilde{\rho}^2 \, \cvolumee
  =
  -
  \int_{\Omega}
  \left. \pd{\thpressure}{\rho}(\temp, \rho) \right|_{\temp = \widehat{\temp}, \, \rho = \widehat{\rho}} \widetilde{\rho}
  \divergence{\widetilde{\vec{v}}}
  \, \cvolumee.
\end{equation}
Similarly, if we multiply the linearised temperature evolution equation~\eqref{eq:52} by $\frac{\widetilde{\temp}}{\widehat{\temp}}$, and we integrate over the domain $\Omega$, then we get
\begin{equation}
  \label{eq:56}
  \widehat{\rho} \frac{\cheatvol\left(\widehat{\temp}, \widehat{\rho}\right)}{\widehat{\temp}}
  \int_{\Omega} \frac{1}{2}\pd{}{t} \widetilde{\temp}^2 \, \cvolumee
  =
  -
  \int_{\Omega}
  \left. \pd{\thpressure}{\temp}(\temp, \rho) \right|_{\temp = \widehat{\temp}, \, \rho = \widehat{\rho}} \widetilde{\temp}
  \divergence{\widetilde{\vec{v}}}
  \, \cvolumee
  -
  \int_{\Omega}
  \kappa
  \frac{
    \vectordot{\nabla \widetilde{\temp}}{\nabla{\widetilde{\temp}}}
  }
  {
    \widehat{\temp}
  }
  \, \cvolumee
  ,
\end{equation}
where we have again used the integration by parts and the fact that $\left. \kappa \vectordot{\nabla \widetilde{\temp}}{\vec{n}} \right|_{\partial \Omega} = 0$, see the zero heat flux boundary condition~\eqref{eq:20}.

Now we take the sum of~\eqref{eq:54}, \eqref{eq:55} and \eqref{eq:56}, and we get the equality
\begin{multline}
  \label{eq:57}
  \frac{1}{2}
  \dd{}{t}
  \left(
    \widehat{\rho}
    \int_{\Omega} \absnorm{\widetilde{\vec{v}}} ^2 \, \cvolumee
    +
    \frac{1}{\widehat{\rho}}\left. \pd{\thpressure}{\rho}(\temp, \rho) \right|_{\temp = \widehat{\temp}, \, \rho = \widehat{\rho}}
    \int_{\Omega} \widetilde{\rho}^2 \, \cvolumee
    +
    \widehat{\rho} \frac{\cheatvol\left(\widehat{\temp}, \widehat{\rho}\right)}{\widehat{\temp}}
    \int_{\Omega} \widetilde{\temp}^2 \, \cvolumee
  \right)
  \\
  =
  -
  \int_{\Omega}
  \tilde{\lambda} \left(\divergence \widetilde{\vec{v}}\right)^2
  \, \cvolumee
  -
  \int_{\Omega}
  2 \nu \tensordot{\traceless{\widetilde{\gradsym}}}{\traceless{\widetilde{\gradsym}}}
  \, \cvolumee
  -
  \int_{\Omega}
  \kappa
  \frac{
    \vectordot{\nabla \widetilde{\temp}}{\nabla{\widetilde{\temp}}}
  }
  {
    \widehat{\temp}
  }
  \, \cvolumee
  .
\end{multline}

If the coefficients $\widehat{\rho}$, $\frac{1}{\widehat{\rho}}\left. \pd{\thpressure}{\rho}(\temp, \rho) \right|_{\temp = \widehat{\temp}, \, \rho = \widehat{\rho}}$ and $\frac{\cheatvol\left(\widehat{\temp}, \widehat{\rho}\right)}{\widehat{\temp}}$ on the left-hand side of~\eqref{eq:57} are positive, then we see that~\eqref{eq:57} could lead to the desired stability result~\eqref{eq:53}. Indeed, in this case~\eqref{eq:57} is a statement regarding the decay of the Lebesgue norm~$\norm[\sleb{2}{\Omega}]{\cdot}$ of the perturbations to the velocity, density and temperature field. (We recall that the right-hand side is negative everywhere except at zero velocity perturbation $\widetilde{\vec{v}}$ and the spatially homogeneous temperature perturbation $\widetilde{\temp}$. The negativity of the right-hand side is in fact a consequence of the second law of thermodynamics, which forces us to fix the sign of the constants $\nu$, $\tilde{\lambda}$ and $\kappa$ in the convenient way from the perspective of the stability.) Consequently, we can conjecture that the stability of the spatially homogeneous rest state with respect to infinitesimal perturbations is guaranteed provided that we enforce additional requirements 
\begin{subequations}
  \label{eq:58}
  \begin{align}
    \label{eq:59}
    \cheatvol \left(\temp, \rho\right) &> 0, \\
    \label{eq:60}
    \pd{\thpressure}{\rho}(\temp, \rho) &> 0.
  \end{align}
\end{subequations}
(Since the spatially homogeneous rest state temperature and density can be arbitrary, we want the inequalities to be satisfied for all $\widehat{\temp}$ and $\widehat{\rho}$.) We summarise our findings in a concise form, see Summary~\ref{summary:decay-linearised}.

\begin{summary}[Decay equation in the linearised setting and a conjecture regarding the Helmholtz free energy]
  \label{summary:decay-linearised}
  Let us consider the problem in Question~\ref{q:isolated}. The Navier--Stokes--Fourier equations \emph{linearised} in the neighbourhood of the spatially homogeneous rest state~\eqref{eq:spatially-homogenoeus-rest-state} imply that small perturbations $
  [
  \widetilde{\rho}, \widetilde{\vec{v}}, \widetilde{\temp}
  ]
  $, see~\eqref{eq:41}, to the spatially homogeneous rest state~\eqref{eq:spatially-homogenoeus-rest-state} satisfy:
  \begin{multline}
    \label{eq:262}
    \frac{1}{2}
    \dd{}{t}
    \left(
      \widehat{\rho}
      \int_{\Omega} \absnorm{\widetilde{\vec{v}}} ^2 \, \cvolumee
      +
      \frac{1}{\widehat{\rho}}\left. \pd{\thpressure}{\rho}(\temp, \rho) \right|_{\temp = \widehat{\temp}, \, \rho = \widehat{\rho}}
      \int_{\Omega} \widetilde{\rho}^2 \, \cvolumee
      +
      \widehat{\rho} \frac{\cheatvol\left(\widehat{\temp}, \widehat{\rho}\right)}{\widehat{\temp}}
      \int_{\Omega} \widetilde{\temp}^2 \, \cvolumee
    \right)
    \\
    =
    -
    \int_{\Omega}
    \tilde{\lambda} \left(\divergence \widetilde{\vec{v}}\right)^2
    \, \cvolumee
    -
    \int_{\Omega}
    2 \nu \tensordot{\traceless{\widetilde{\gradsym}}}{\traceless{\widetilde{\gradsym}}}
    \, \cvolumee
    -
    \int_{\Omega}
    \kappa
    \frac{
      \vectordot{\nabla \widetilde{\temp}}{\nabla{\widetilde{\temp}}}
    }
    {
      \widehat{\temp}
    }
    \, \cvolumee
    .
  \end{multline}

  We thus conjecture that the stability of an arbitrary spatially homogeneous rest state~\eqref{eq:spatially-homogenoeus-rest-state} is guaranteed provided that the Helmholtz free energy is chosen in such a way that
  \begin{subequations}
    \label{eq:thermodynamic-stability-conditions-conjecture}
    \begin{align}
      \label{eq:264}
      \cheatvol \left(\temp, \rho\right) &> 0, \\
      \label{eq:265}
      \pd{\thpressure}{\rho}(\temp, \rho) &> 0.
    \end{align}
  Note that in virtue of the formulae for the specific heat capacity at constant volume~\eqref{eq:4} and the thermodynamic pressure~\eqref{eq:5} the requirements~\eqref{eq:thermodynamic-stability-conditions-conjecture} are requirements on the second derivatives of the Helmholtz free energy. 
\end{subequations}
\end{summary}

\subsection{Beyond linearised setting---net entropy}
\label{sec:beyond-line-sett}
If we want to go beyond the linearised setting we need to construct a Lyapunov type functional, that is a functional that decays in time and that controls the size of perturbations. One might think that the \emph{net entropy}
\begin{equation}
  \label{eq:267}
  \netentropy =_{\bydefinition}  \int_{\Omega} \rho \entropy \, \cvolumee, 
\end{equation}
can be useful here. In particular, in our thermodynamically isolated system we have
\begin{equation}
  \label{eq:266}
  \dd{}{t}
  \int_{\Omega} \rho \entropy \, \cvolumee
  =
  \int_{\Omega}
  \left(
    \frac{
      \tilde{\lambda} \left(\divergence \vec{v}\right)^2
      +
      2 \nu \tensordot{\traceless{\gradsym}}{\traceless{\gradsym}}
    }
    {
      \temp
    }
    +
    \frac{
      \kappa
      \vectordot{\nabla \temp}{\nabla{\temp}}
    }
    {
      \temp^2
    }
  \right)
  \, \cvolumee
  ,
\end{equation}
hence the negative net entropy, $- \netentropy$, decays in time, which is what we want for the Lyapunov functional. So far so good.

Unfortunately, if we write down the explicit formula for $\netentropy$ for example for the calorically perfect ideal gas, see~\eqref{eq:11}, we see that
\begin{equation}
  \label{eq:268}
  \netentropy
  =
  \int_{\Omega}
  \rho
  \entropy (\temp, \rho)
  \, \cvolumee
  =
  \int_{\Omega}
  \rho
  \cheatvolref
  \ln
  \left[
    \frac{\temp}{\tempeq}
    \left(
      \frac{\rho}{\rhoeq}
    \right)^{1 - \gamma}
  \right]
  \, \cvolumee
  =
  \int_{\Omega}
  \rho
  \cheatvolref
  \ln
  \left[
    \left(1 + \frac{\widetilde{\temp}}{\tempeq}\right)
    \left(
      1 + \frac{\widetilde{\rho}}{\rhoeq}
    \right)^{1 - \gamma}
  \right]
  \, \cvolumee
  ,
\end{equation}
where we have conveniently chosen the normalisation constants $\rhoref$ and $\tempref$ in~\eqref{eq:11} as $\rhoref =_{\bydefinition} \rhoeq$ and $\tempref =_{\bydefinition} \tempeq$. Formula~\eqref{eq:268} shows that~\eqref{eq:266} provides us, unlike~\eqref{eq:262}, no way to track the evolution of the size of the perturbation~$\widetilde{\temp}$ and~$\widetilde{\rho}$ to the spatially homogeneous rest state~\eqref{eq:spatially-homogenoeus-rest-state}. In particular, if the net entropy~\eqref{eq:268} is zero, it is not necessarily true that the perturbations $\widetilde{\temp}$ and $\widetilde{\rho}$ vanish in the whole domain $\Omega$. The same happens for more complex substances.

Consequently, tracking the net entropy evolution is a step in the right direction, but the net entropy alone is useless. The moral of this observation is summarised in Summary~\ref{summary:entropy-is-useless}. Note, however, that the net entropy can be still useful in the sense that its \emph{second derivative} matrix might induce a positive definite quadratic form that can be exploited to measure the size of the perturbations, see Section~\ref{sec:lead-order-expans} for details.

\begin{summary}[Entropy alone is useless in stability analysis]
  \label{summary:entropy-is-useless}
  The negative net entropy decays in time, but it fails to provide control over the ``size'' of perturbation.
\end{summary}

\subsection{Beyond linearised setting---Lyapunov type functional}
\label{sec:lyap-type-funct-1}
We have to rethink and refine the unsuccessful attempt with the net entropy. The key idea is that instead of the net entropy evolution we must track \emph{net entropy evolution subject to all possible constraints}. This is motivated by Clausius's famous statement, \cite[page 400]{clausius.r:ueber}: ``The energy of the world is constant. The entropy of the world strives to a maximum.'' In our case the constraints are the constant net total energy and constant net mass. 

In particular, we can guess that a convenient functional for the nonlinear stability analysis of spatially homogeneous rest state in a compressible heat conducting fluid is the functional
\begin{equation}
  \label{eq:61}
  \mathcal{V}_{\mathrm{eq}}
  =_{\bydefinition}
  -
  \left(
    \netentropy
    -
    \netentropy_{\max}
  \right)
  +
  \lambda_1
  \left( \nettenergy - \widehat{\nettenergy} \right)
  +
  \lambda_2
  \int_{\Omega}
  \left(
    \rho
    -
    \widehat{\rho}
  \right)
  \,
  \cvolumee
  ,
\end{equation}
where
\begin{subequations}
  \label{eq:62}
  \begin{align}
    \label{eq:63}
    \netentropy &=_{\bydefinition}  \int_{\Omega} \rho \entropy(\temp, \rho) \, \cvolumee, \\
    \label{eq:64}
    \nettenergy &=_{\bydefinition}
    \int_{\Omega}
                  \left(
                  \frac{1}{2} \rho \absnorm{\vec{v}}^2
                  +
                  \rho
                  \ienergy(\temp, \rho)
                  \right)
    \,
                  \cvolumee
                  ,
  \end{align}
\end{subequations}
denote the net entropy and the net total energy; $\entropy$ and $\ienergy$ denote the entropy density and the internal energy density respectively. The symbol $\widehat{\nettenergy}$ denotes the net total energy $\nettenergy$ evaluated at the target spatially homogeneous rest state~\eqref{eq:spatially-homogenoeus-rest-state}. The constant shift by $\netentropy_{\max}$ is just a convenient normalisation---we want the functional $\mathcal{V}_{\mathrm{eq}}$ to vanish at the state with the maximum entropy $\netentropy_{\max}$, which is the desired property for a Lyapunov functional.  

If we want to use the functional~\eqref{eq:61}, we must first identify the multipliers $\lambda_1$ and $\lambda_2$. We want to choose the multipliers in such a way that the net entropy subject to all constraints is maximal at the spatially homogeneous rest state~\eqref{eq:spatially-homogenoeus-rest-state}, that is we have $\netentropy_{\max} = \netentropy(\widehat{\temp}, \widehat{\rho}) =  \int_{\Omega} \widehat{\rho} \entropy( \widehat{\temp}, \widehat{\rho}) \, \cvolumee$. The identification of multipliers requires some long algebraic manipulations, and it is done in detail in Section~\ref{sec:ident-lagr-mult}, see also~\cite{bulcek.m.malek.j.ea:thermodynamics} for the incompressible case. The multipliers are identified as
\begin{subequations}
  \label{eq:65}
  \begin{align}
    \label{eq:66}
    \lambda_1 &= \frac{1}{\widehat{\temp}}, \\
    \label{eq:67}
    \lambda_2 &=
                -
                \frac{\widehat{\thpressure}}{\widehat{\rho} \widehat{\temp}}
                -
                \frac{\widehat{\fenergy}}{\widehat{\temp}}
  \end{align}
\end{subequations}
and if we further multiply~\eqref{eq:61} by the rest state temperature\footnote{Concerning the rationale behind the multiplication by $\widehat{\temp}$, see~\cite{bulcek.m.malek.j.ea:thermodynamics} for details.} $\widehat{\temp}$, then we finally arrive at the functional
\begin{multline}
  \label{eq:68}
  \mathcal{V}_{\mathrm{meq}, \, \widehat{\temp}, \, \widehat{\rho}}(\rho, \vec{v}, \temp)
    =
    -
    \int_{\Omega} \widehat{\temp} \left( \rho \entropy(\temp, \rho) - \widehat{\rho} \entropy(\widehat{\temp}, \widehat{\rho}) \right) \, \cvolumee
    +
    \int_{\Omega} \left(\frac{1}{2} \rho \absnorm{\vec{v}}^2 + \rho \ienergy(\temp, \rho) - \widehat{\rho} \ienergy(\widehat{\temp}, \widehat{\rho}) \right) \, \cvolumee
    -
    \int_{\Omega} \left( \frac{\thpressure(\widehat{\temp}, \widehat{\rho})}{\widehat{\rho}} + \fenergy(\widehat{\temp}, \widehat{\rho}) \right) \left(\rho  - \widehat{\rho} \right)  \, \cvolumee
    .
  \end{multline}
We note that the rest state density $\widehat{\rho}$ and temperature $\widehat{\temp}$ are constants both in space and time, and that $\widehat{\rho}$ and $\widehat{\temp}$ are parameters in the functional~$\mathcal{V}_{\mathrm{meq}, \, \widehat{\temp}, \, \widehat{\rho}}$; the functional itself takes as arguments only the triple $\rho$, $\vec{v}$, $\temp$ (the current state). The outlined construction leads to the functional $\mathcal{V}_{\mathrm{meq}, \, \widehat{\temp}, \, \widehat{\rho}}(\rho, \vec{v}, \temp)$ that is
\begin{enumerate}
\item decreasing in time,
\item nonnegative,
\item vanishing if and only if the perturbation vanishes.
\end{enumerate}
(These properties are proven in Section~\ref{sec:time-deriv-lyap} and Section~\ref{sec:nonnegativity}.) Consequently, the functional~\eqref{eq:68} might indeed be useful in nonlinear stability analysis of \emph{spatially homogeneous rest state} $[\widehat{\rho}, \vec{0}, \widehat{\temp}]$---it can help to answer the question whether any solution $[\rho, \vec{v}, \temp]$ to the governing equations that starts from an arbitrary initial condition eventually approaches the spatially homogeneous rest state $[\widehat{\rho}, \vec{0}, \widehat{\temp}]$, that is if we have, in some sense, the convergence property~\eqref{eq:convergenece-isolated}. 

In continuum thermodynamics setting---spatially distributed systems---the idea to use functionals of type~\eqref{eq:68} in nonlinear stability analysis was introduced by \cite{coleman.bd.greenberg.jm:thermodynamics}, \cite[Equation 2.6]{coleman.bd:on} and \cite{gurtin.me:thermodynamics*1,gurtin.me:thermodynamics}, though its origins can be traced back to~\cite{duhem.p:traite}. (See also~\cite{silhavy.m:mechanics} and~\cite{dafermos.cm:second}.) The rationale behind the choice of the functional, in particular the choice of the multipliers, is however not always clear in these works.

The utility of functional~\eqref{eq:68} in nonlinear stability analysis is evident. The guaranteed nonpositivity of the time derivative $\dd{\mathcal{V}_{\mathrm{meq}, \, \widehat{\temp}, \, \widehat{\rho}}}{t}$ and the nonnegativity of the functional $\mathcal{V}_{\mathrm{meq}, \, \widehat{\temp}, \, \widehat{\rho}}$ make the functional an ideal \emph{candidate} for the Lyapunov functional for the nonlinear stability analysis. (See~\cite{la-salle.j.lefschetz.s:stability}, \cite{yoshizawa.t:stability} and \cite{henry.d:geometric} for the concept of Lyapunov functional.) Note, however, that in the infinite dimensional setting---spatially distributed systems---the requirements on the Lyapunov functionals are in general stricter than in the finite dimensional setting. In particular, the genuine Lyapunov functional must be related to a \emph{suitable norm/metric} on the corresponding state space, which might be difficult to obtain. We do not discuss the issues related to the appropriate choice of norm/metric, hence we prefer to denote functionals of this type only as \emph{Lyapunov type} functionals.

\subsection{Time derivative of Lyapunov type functional}
\label{sec:time-deriv-lyap}
The time derivative of the functional~\eqref{eq:68} reads
\begin{equation}
  \label{eq:70}
  \dd{\mathcal{V}_{\mathrm{meq}, \, \widehat{\temp}, \, \widehat{\rho}}}{t}
  =
  -
  \widehat{\temp}
  \dd{\netentropy}{t}
  =
  -
  \widehat{\temp}
  \int_{\Omega}
  \entprodc
  \, \cvolumee
  \leq 0
  .
\end{equation}
Equation~\eqref{eq:70} follows from several facts.

\emph{First}, the net total energy $\nettenergy$ and the net mass are conserved in a thermodynamically isolated system. \emph{Second}, the generic entropy evolution equation reads
\begin{equation}
  \label{eq:generic-entropy-evolution-equation}
  \rho \dd{\entropy}{t}
  +
  \divergence \entfluxc
  =
  \entprodc
  ,
\end{equation}
where $\entfluxc$ denotes the entropy flux and $\entprodc$ denotes the nonnegative entropy production. In the thermodynamically isolated system there is no \emph{entropy flux $\entfluxc$ through the boundary}, that is we have
\begin{equation}
  \label{eq:71}
  \left. \vectordot{\entfluxc}{\vec{n}} \right|_{\partial \Omega} = 0.
\end{equation}
(In our case we have the zero heat flux boundary condition~\eqref{eq:boundary-conditions-isolated}, and we know that the entropy flux $\entfluxc$ is just proportional to the heat flux $\hfluxc$, namely $\entfluxc = \frac{\hfluxc}{\temp}$.) Consequently, the definition of net total entropy~\eqref{eq:63} and the entropy evolution equation~\eqref{eq:generic-entropy-evolution-equation} imply that
\begin{equation}
  \label{eq:72}
  \dd{\netentropy}{t}
  =
  \int_{\Omega}
  \entprodc
  \, \cvolumee
  -
  \int_{\partial \Omega} \vectordot{\entfluxc}{\vec{n}} \, \csurfacees
  =
  \int_{\Omega}
  \entprodc
  \, \cvolumee
  ,
\end{equation}
where the entropy production $\entprodc$ is a nonnegative quantity. \emph{Third}, the temperature at the equilibrium steady state $\widehat{\temp}$ is, in the case of a thermodynamically isolated system, constant, hence it can be taken out of the first integral in~\eqref{eq:68}. We note that the sign of the time derivative does not depend on the Lagrange multipliers $\lambda_1$ and $\lambda_2$; they do not enter the formula~\eqref{eq:70}.

For our compressible heat conducting fluid we thus have
\begin{equation}
  \label{eq:21}
  \dd{\mathcal{V}_{\mathrm{meq}, \, \widehat{\temp}, \, \widehat{\rho}}}{t}
  =
  -
  \widehat{\temp}
  \int_{\Omega}
  \left(
    \frac{
      \tilde{\lambda} \left(\divergence \vec{v}\right)^2
      +
      2 \nu \tensordot{\traceless{\gradsym}}{\traceless{\gradsym}}
    }
    {
      \temp
    }
    +
    \frac{
      \kappa
      \vectordot{\nabla \temp}{\nabla{\temp}}
    }
    {
      \temp^2
    }
  \right)
  \, \cvolumee
  ,
\end{equation}
where the density, velocity and temperature on the right-hand side are the density, velocity and temperature at the current state that is $\temp = \widehat{\temp} + \widetilde{\temp}$ and so forth. We also note that the right-hand side vanishes at the spatially homogeneous rest state---all gradients are gone at the spatially homogeneous rest state.

\subsection{Identification of Lagrange multipliers in Lyapunov type functional}
\label{sec:ident-lagr-mult}
The Lagrange multipliers in~\eqref{eq:61} can be identified via the solution of a constrained maximisation problem, see~\cite{bulcek.m.malek.j.ea:thermodynamics} for details. However, \cite{bulcek.m.malek.j.ea:thermodynamics} dealt with an incompressible material only; the treatment of a compressible material is discussed below.

We want the net entropy $\netentropy$ at the spatially homogeneous rest state~\eqref{eq:spatially-homogenoeus-rest-state} to be maximal subject to the corresponding constraints.  The auxiliary functional for the constrained maximisation problem is, up to the sign, the functional~\eqref{eq:61}. If we use the definitions of the net entropy and the net total energy, we get the auxiliary functional in the form
\begin{equation}
  \label{eq:auxiliary-functional}
  \mathcal{L}_{\lambda_1, \lambda_2} = _{\bydefinition}
  \int_{\Omega} \left( \rho \entropy - \widehat{\rho} \widehat{\entropy} \right) \, \cvolumee
  -
  \lambda_1 \int_{\Omega} \left(\frac{1}{2} \rho \absnorm{\vec{v}}^2 + \rho \ienergy - \widehat{\rho} \widehat{\ienergy} \right) \, \cvolumee
  -
  \lambda_2 \int_{\Omega} \left(\rho  - \widehat{\rho} \right)  \, \cvolumee
  .
\end{equation}
(In~\eqref{eq:auxiliary-functional} we use notation $\widehat{\entropy} \equiv \entropy(\widehat{\temp}, \widehat{\rho})$ and $\widehat{\ienergy} \equiv \ienergy(\widehat{\temp}, \widehat{\rho})$.) The spatially homogeneous rest state $\widehat{\rho}$, $\widehat{\vec{v}} = \vec{0}$ and $\widehat{\temp}$, see~\eqref{eq:spatially-homogenoeus-rest-state}, is a solution to the maximisation problem provided that the G\^ateaux derivative\footnote{We recall that the G\^ateaux derivative $\Diff \mathcal{M} (\vec{x})[\vec{y}]$ of a functional $\mathcal{M}$ at point $\vec{x}$ in the direction $\vec{y}$ is defined as
$
    \Diff \mathcal{M} (\vec{x})[\vec{y}]
    =_{\bydefinition}
    \lim_{s \to 0}
    \frac{
      \mathcal{M} (\vec{x} + s \vec{y}) - \mathcal{M} (\vec{x})
    }
    {
      s
    }
$ 
which is tantamount to 
$
\Diff \mathcal{M} (\vec{x})[\vec{y}]
=_{\bydefinition}
\left.
  \dd{}{s}
  \mathcal{M} (\vec{x} + s \vec{y})
\right|_{s=0}
$. The second G\^ateaux derivative $\Diff^2 \mathcal{M} (\vec{x})[\vec{y}, \vec{z}]$ of the functional $\mathcal{M}$ at point $\vec{x}$ is the bilinear form
$
\Diff^2 \mathcal{M} (\vec{x})[\vec{y}, \vec{z}]
=_{\bydefinition}
\left.
  \pd{^2}{s \partial t}
  \mathcal{M} (\vec{x} + s \vec{y} + t \vec{z})
\right|_{s = t = 0}
$,
and for $\vec{z} = \vec{y}$ we get the corresponding quadratic form
$
\Diff^2 \mathcal{M} (\vec{x})[\vec{y}, \vec{y}]
=
\left.
  \ddd{}{s}
  \mathcal{M} (\vec{x} + s \vec{y})
\right|_{s=0}
$. If it is necessary to emphasise the variable with respect to which we differentiate, we also write $\Diff[\vec{x}] \mathcal{M} (\vec{x})[\vec{y}]$ instead of $\Diff \mathcal{M} (\vec{x})[\vec{y}]$.
} 
of auxiliary functional~\eqref{eq:auxiliary-functional} at point $\widehat{\rho}$, $\widehat{\vec{v}}$ and~$\widehat{\temp}$  vanishes in every admissible direction $\widetilde{\rho}$, $\widetilde{\vec{v}}$ and $\widetilde{\temp}$.

Now we evaluate the G\^ateaux derivative in \emph{two different descriptions}. As we shall see, the identification of the second multiplier~$\lambda_2$ hinges on whether we impose the \emph{integral} mass conservation constraint on the perturbations or not, and the two descriptions provide two independent ways of identifying the multipliers~$\lambda_1$ and~$\lambda_2$. The idea is the following. The fact that the net entropy attains its maximum value must be true no matter whether our primitive variables are the temperature and the density, or the temperature and the pressure and so forth; consequently \emph{we can conveniently switch between various descriptions in order to get the desired piece of information}.

\subsubsection{G\^ateaux derivative in temperature--density representation}
\label{sec:gate-deriv-temp}

First, we interpret the entropy $\entropy$ and the internal energy $\ienergy$ in~\eqref{eq:auxiliary-functional} as functions of the density $\rho$ and the temperature $\temp$. The formula for the G\^ateaux derivative at point $\widehat{\rho}$, $\widehat{\vec{v}}$ and $\widehat{\temp}$ in the direction $\widetilde{\rho}$, $\widetilde{\vec{v}}$ and $\widetilde{\temp}$ reads
\begin{multline}
  \label{eq:74}
  \Diff \mathcal{L}_{\lambda_1, \lambda_2}(\widehat{\rho}, \vec{0}, \widehat{\temp})[\widetilde{\rho}, \widetilde{\vec{v}}, \widetilde{\temp}]
  =
  \\
  \left.
    \dd{}{s}
    \left\{
      \int_{\Omega}
      \left[
        \left( \widehat{\rho} + s\widetilde{\rho} \right) \entropy(\widehat{\temp} + s\widetilde{\temp}, \widehat{\rho} + s\widetilde{\rho})
        -
        \widehat{\rho}\entropy(\widehat{\temp}, \widehat{\rho})
      \right]\, \cvolumee
      -
      \lambda_1
      \int_{\Omega}
      \left(\frac{1}{2} \left( \widehat{\rho} + s\widetilde{\rho} \right) \absnorm{s \widetilde{\vec{v}}}^2
        +
        \left( \widehat{\rho} + s\widetilde{\rho} \right)
        \ienergy(\widehat{\temp} + s\widetilde{\temp}, \widehat{\rho} + s\widetilde{\rho}) - \widehat{\rho} \ienergy(\widehat{\temp}, \widehat{\rho}) \right)
      \, \cvolumee
    \right.
    \right.
    \\
    \left.
    \left.
      -
      \lambda_2 \int_{\Omega} \left(\widehat{\rho} + s\widetilde{\rho}  - \widehat{\rho} \right)  \, \cvolumee
    \right\}
  \right|_{s=0}
  ,
\end{multline}
which with a slight abuse of notation\footnote{The abuse of notation is about using $\pd{\entropy(\widehat{\temp}, \widehat{\rho})}{\widehat{\temp}}$ as an abbreviation for $\left. \pd{\entropy(\temp, \rho)}{\temp} \right|_{\left[\temp, \rho\right] = \left[\widehat{\temp}, \widehat{\rho}\right]}$.}  yields
\begin{equation}
  \label{eq:75}
  \Diff \mathcal{L}_{\lambda_1, \lambda_2}(\widehat{\rho}, \vec{0}, \widehat{\temp})[\widetilde{\rho}, \widetilde{\vec{v}}, \widetilde{\temp}]
  =
  \int_{\Omega}
  \widehat{\rho}
  \left(
    \pd{\entropy(\widehat{\temp}, \widehat{\rho})}{\widehat{\temp}}
    -
    \lambda_1
    \pd{\ienergy(\widehat{\temp}, \widehat{\rho})}{\widehat{\temp}}
  \right)
  \widetilde{\temp}
  \, \cvolumee
  ,
\end{equation}
where we have used the fact that $\int_{\Omega} \widetilde{\rho} \, \cvolumee = 0$, which is a consequence of the mass conservation constraint. This consequence of mass conservation in fact eliminates all terms that are linear in $\widetilde{\rho}$, hence also all terms that contain the second Lagrange multiplier $\lambda_2$.  (Recall also that the stationary state is spatially homogeneous, hence $\widehat{\temp}$ and $\widehat{\rho}$ are constants.) Using standard thermodynamic identities
\begin{subequations}
  \label{eq:76}
  \begin{align}
    \label{eq:77}
\pd{\entropy(\widehat{\temp}, \widehat{\rho})}{\widehat{\temp}}
&=
\frac{\cheatvol (\widehat{\temp}, \widehat{\rho})}{\widehat{\temp}},
    \\
    \label{eq:78}
    \pd{\ienergy(\widehat{\temp}, \widehat{\rho})}{\widehat{\temp}}
  &=
  \cheatvol (\widehat{\temp}, \widehat{\rho}),
  \end{align}
\end{subequations}
where $\cheatvol$ denotes the specific heat capacity at constant volume, we see that~\eqref{eq:75} reduces to
\begin{equation}
  \label{eq:79}
  \Diff \mathcal{L}_{\lambda_1, \lambda_2}(\widehat{\rho}, \vec{0}, \widehat{\temp})[\widetilde{\rho}, \widetilde{\vec{v}}, \widetilde{\temp}]
  =
  \int_{\Omega}
  \widehat{\rho}
  \left(
    \frac{1}{\widehat{\temp}}
    -
    \lambda_1
  \right)
  \cheatvol (\widehat{\temp}, \widehat{\rho})
  \widetilde{\temp}
  \, \cvolumee
  .
\end{equation}
The G\^ateaux derivative therefore vanishes for arbitrary $\widetilde{\temp}$ provided that we fix the Lagrange multiplier as
\begin{equation}
  \label{eq:80}
  \lambda_1 = \frac{1}{\widehat{\temp}}.
\end{equation}

The second Lagrange multiplier $\lambda_2$ is however still unidentified, because we are exploiting the \emph{integral} mass conservation constraint $\int_{\Omega} \widetilde{\rho} \, \cvolumee = 0$, which eliminates all terms linear in $\widetilde{\rho}$.

If we had instead required the stronger \emph{local stationarity} condition in the sense that the \emph{integrand} in~\eqref{eq:auxiliary-functional} is stationary at the spatially homogeneous rest state \emph{pointwise}, that is for arbitrary perturbations $\widetilde{\rho}$ and $\widetilde{\temp}$ that are not necessarily subject to the integral constraint, then we would get stationarity conditions
\begin{subequations}
  \label{eq:247}
  \begin{align}
    \label{eq:248}
    \widehat{\rho}
    \left(
    \pd{\entropy(\widehat{\temp}, \widehat{\rho})}{\widehat{\temp}}
    -
    \lambda_1
    \pd{\ienergy(\widehat{\temp}, \widehat{\rho})}{\widehat{\temp}}
    \right)
    &=
      0,
    \\
    \label{eq:336}
    \entropy(\widehat{\temp}, \widehat{\rho})
    -
    \lambda_1  \ienergy(\widehat{\temp}, \widehat{\rho})
    +
    \widehat{\rho}
    \left(
    \pd{\entropy}{\rho}(\widehat{\temp}, \widehat{\rho})
    -
    \lambda_1 \pd{\ienergy}{\rho}(\widehat{\temp}, \widehat{\rho})
    \right)
    -
    \lambda_2
    &=
      0
      .
  \end{align}
\end{subequations}
The first condition would again give us $\lambda_1 = \frac{1}{\widehat{\temp}}$. Concerning the second condition, we can, upon using $\lambda_1 = \frac{1}{\widehat{\temp}}$, rewrite it as
\begin{equation}
  \label{eq:337}
  -
  \frac{\fenergy(\widehat{\temp}, \widehat{\rho})}{\widehat{\temp}}
  +
  \widehat{\rho}
  \left(
    \pd{\entropy}{\rho}(\widehat{\temp}, \widehat{\rho})
    -
    \frac{1}{\widehat{\temp}}
    \pd{\ienergy}{\rho}(\widehat{\temp}, \widehat{\rho})
  \right)
  -
  \lambda_2
  =
  0
  .
\end{equation}
The formulae~\eqref{eq:5}, \eqref{eq:205} and \eqref{eq:289} imply that $\pd{\ienergy}{\rho}(\temp, \rho) = \temp \pd{\entropy}{\rho}(\temp, \rho) + \frac{\thpressure(\temp, \rho)}{\rho^2}$, hence we further get
\begin{equation}
  \lambda_2
  =
  - \frac{\widehat{\fenergy}}{\widehat{\temp}} - \frac{\widehat{\thpressure}}{\widehat{\rho} \widehat{\temp}}
  ,
\end{equation}
where we have denoted $\widehat{\fenergy} = \fenergy(\widehat{\temp}, \widehat{\rho})$ and $\widehat{\thpressure} = \thpressure(\widehat{\temp}, \widehat{\rho})$. This is the same $\lambda_2$ value as the one obtained below in the temperature--pressure representation, see~\eqref{eq:89}. In this sense the temperature--pressure representation serves as an independent cross-check, \emph{but} we note that in the temperature--pressure representation we require the stationarity of the \emph{integrated} functional, and we avoid the \emph{integral} constraint on the density perturbation. Since the perturbations $\widetilde{\temp}$ and $\widetilde{\thpressure}$ are then arbitrary, this is equivalent to the \emph{local} stationarity condition.

\subsubsection{G\^ateaux derivative in temperature--pressure representation}
\label{sec:gate-deriv-temp-1}
 We interpret the entropy $\entropy$ and the internal energy $\ienergy$ in~\eqref{eq:auxiliary-functional} as functions of the temperature $\temp$ and the thermodynamic pressure $\thpressure$. The formula for the G\^ateaux derivative at point $\widehat{\thpressure}$, $\widehat{\vec{v}}$ and $\widehat{\temp}$ in the direction $\widetilde{\thpressure}$, $\widetilde{\vec{v}}$ and $\widetilde{\temp}$ reads
\begin{multline}
  \label{eq:81}
  \Diff \mathcal{L}_{\lambda_1, \lambda_2}(\widehat{\thpressure}, \vec{0}, \widehat{\temp})[\widetilde{\thpressure}, \widetilde{\vec{v}}, \widetilde{\temp}]
  =
  \left.
    \dd{}{s}
    \left\{
      \int_{\Omega}
      \left(
        \vphantom{\frac{1}{2}}
        \rho(\widehat{\temp} + s\widetilde{\temp}, \widehat{\thpressure} + s\widetilde{\thpressure})
        \entropy(\widehat{\temp} + s\widetilde{\temp}, \widehat{\thpressure} + s\widetilde{\thpressure})
        -
        \rho(\widehat{\temp}, \widehat{\thpressure}) \entropy(\widehat{\temp}, \widehat{\thpressure})
      \right)
      \, \cvolumee
    \right.
    \right.
    \\
    \left.
    \left.
      -
      \lambda_1
      \int_{\Omega}
      \left(
        \frac{1}{2}  \rho(\widehat{\temp} + s\widetilde{\temp}, \widehat{\thpressure} + s\widetilde{\thpressure}) \absnorm{s \widetilde{\vec{v}}}^2
        +
        \rho(\widehat{\temp} + s\widetilde{\temp}, \widehat{\thpressure} + s\widetilde{\thpressure})
        \ienergy(\widehat{\temp} + s\widetilde{\temp}, \widehat{\thpressure} + s\widetilde{\thpressure})
      \right) \, \cvolumee
    \right.
    \right.
    \\
    \left.
    \left.
      -
      \lambda_2 \int_{\Omega} \rho(\widehat{\temp} + s\widetilde{\temp}, \widehat{\thpressure} + s\widetilde{\thpressure}) \, \cvolumee
    \right\}
  \right|_{s=0}
  .
\end{multline}
(The density is now interpreted as a function of the primitive variables---the temperature and the thermodynamic pressure.) Straightforward calculation reveals that
\begin{multline}
  \label{eq:82}
  \Diff \mathcal{L}_{\lambda_1, \lambda_2}(\widehat{\thpressure}, \vec{0}, \widehat{\temp})[\widetilde{\thpressure}, \widetilde{\vec{v}}, \widetilde{\temp}]
  \\
  =
  \int_{\Omega}
  \left\{
    \pd{\rho(\widehat{\temp}, \widehat{\thpressure})}{\widehat{\temp}}
    \left(
      \widehat{\entropy}
      -
      \lambda_1
      \widehat{\ienergy}
    \right)
    +
    \widehat{\rho}
    \left(
      \pd{\entropy(\widehat{\temp}, \widehat{\thpressure})}{\widehat{\temp}}
      -
      \lambda_1
      \pd{\ienergy(\widehat{\temp}, \widehat{\thpressure})}{\widehat{\temp}}
    \right)
    -
    \lambda_2
    \pd{\rho(\widehat{\temp}, \widehat{\thpressure})}{\widehat{\temp}}
  \right\}
  \widetilde{\temp}
  \,
  \cvolumee
  \\
  +
    \int_{\Omega}
    \left\{
         \pd{\rho(\widehat{\temp}, \widehat{\thpressure})}{\widehat{\thpressure}}
    \left(
      \widehat{\entropy}
      -
      \lambda_1
      \widehat{\ienergy}
    \right)
    +
    \widehat{\rho}
    \left(
      \pd{\entropy(\widehat{\temp}, \widehat{\thpressure})}{\widehat{\thpressure}}
      -
      \lambda_1
      \pd{\ienergy(\widehat{\temp}, \widehat{\thpressure})}{\widehat{\thpressure}}
    \right)
    -
    \lambda_2
    \pd{\rho(\widehat{\temp}, \widehat{\thpressure})}{\widehat{\thpressure}}
  \right\}
  \widetilde{\thpressure}
  \,
  \cvolumee
  ,
\end{multline}
where we have denoted $\widehat{\rho} = \rho \left(\widehat{\temp}, \widehat{\thpressure} \right)$,  $\widehat{\entropy} = \entropy \left(\widehat{\temp}, \widehat{\thpressure} \right)$ and $\widehat{\ienergy} = \ienergy\left(\widehat{\temp}, \widehat{\thpressure} \right)$, and where we have again slightly abused the notation.

Now we have to find the multipliers, and we preferably want to find the multipliers \emph{without} invoking the \emph{integral} mass conservation constraint $\int_{\Omega} \left( \rho - \widehat{\rho} \right) \, \cvolumee = 0$. This might, later on, give us a chance to work with \emph{pointwise} formulae. The G\^ateaux derivative~\eqref{eq:82} vanishes in an arbitrary direction $\widetilde{\temp}$ and $\widetilde{\thpressure}$ if the multipliers $\lambda_1$ and $\lambda_2$ are chosen in such a way that the following holds
\begin{subequations}
  \label{eq:111}
  \begin{align}
    \label{eq:308}
    \pd{\rho(\widehat{\temp}, \widehat{\thpressure})}{\widehat{\temp}}
    \left(
      \widehat{\entropy}
      -
      \lambda_1
      \widehat{\ienergy}
    \right)
    +
    \widehat{\rho}
    \left(
      \pd{\entropy(\widehat{\temp}, \widehat{\thpressure})}{\widehat{\temp}}
      -
      \lambda_1
      \pd{\ienergy(\widehat{\temp}, \widehat{\thpressure})}{\widehat{\temp}}
    \right)
    -
    \lambda_2
    \pd{\rho(\widehat{\temp}, \widehat{\thpressure})}{\widehat{\temp}}
    &=
      0,
    \\
    \label{eq:309}
    \pd{\rho(\widehat{\temp}, \widehat{\thpressure})}{\widehat{\thpressure}}
    \left(
      \widehat{\entropy}
      -
      \lambda_1
      \widehat{\ienergy}
    \right)
    +
    \widehat{\rho}
    \left(
      \pd{\entropy(\widehat{\temp}, \widehat{\thpressure})}{\widehat{\thpressure}}
      -
      \lambda_1
      \pd{\ienergy(\widehat{\temp}, \widehat{\thpressure})}{\widehat{\thpressure}}
    \right)
    -
    \lambda_2
    \pd{\rho(\widehat{\temp}, \widehat{\thpressure})}{\widehat{\thpressure}}
    &=
    0
    .
  \end{align}
\end{subequations}
Differentiating the equality $\entropy(\widehat{\temp}, \widehat{\thpressure}) = \entropy(\widehat{\temp}, \rho(\widehat{\temp}, \widehat{\thpressure}))$ we see that
\begin{subequations}
  \label{eq:319}
  \begin{align}
    \label{eq:320}
    \pd{\entropy(\widehat{\temp}, \widehat{\thpressure})}{\widehat{\temp}}
    -
    \lambda_1
    \pd{\ienergy(\widehat{\temp}, \widehat{\thpressure})}{\widehat{\temp}}
    &=
      \left(
      \pd{\entropy(\widehat{\temp}, \widehat{\rho})}{\widehat{\temp}}
      -
      \lambda_1
      \pd{\ienergy(\widehat{\temp}, \widehat{\rho})}{\widehat{\temp}}
      \right)
      +
      \left(
      \pd{\entropy(\widehat{\temp}, \widehat{\rho})}{\widehat{\rho}}
      -
      \lambda_1
      \pd{\ienergy(\widehat{\temp}, \widehat{\rho})}{\widehat{\rho}}
      \right)
      \pd{\rho(\widehat{\temp}, \widehat{\thpressure})}{\widehat{\temp}}
      ,
    \\
    \label{eq:316}
    \pd{\entropy(\widehat{\temp}, \widehat{\thpressure})}{\widehat{\thpressure}}
    -
    \lambda_1
    \pd{\ienergy(\widehat{\temp}, \widehat{\thpressure})}{\widehat{\thpressure}}
    &=
      \left(
      \pd{\entropy(\widehat{\temp}, \widehat{\rho})}{\widehat{\rho}}
      -
      \lambda_1
      \pd{\ienergy(\widehat{\temp}, \widehat{\rho})}{\widehat{\rho}}
      \right)
      \pd{\rho(\widehat{\temp}, \widehat{\thpressure})}{\widehat{\thpressure}}
      ,
  \end{align}
\end{subequations}
which upon substitution into the system of equations~\eqref{eq:111} yields
\begin{subequations}
  \label{eq:136}
  \begin{align}
    \label{eq:314}
    \pd{\rho(\widehat{\temp}, \widehat{\thpressure})}{\widehat{\temp}}
    \left[
    \left(
    \widehat{\entropy}
      -
    \lambda_1
    \widehat{\ienergy}
    \right)
    +
    \widehat{\rho}
    \left(
    \pd{\entropy(\widehat{\temp}, \widehat{\rho})}{\widehat{\rho}}
    -
    \lambda_1
    \pd{\ienergy(\widehat{\temp}, \widehat{\rho})}{\widehat{\rho}}
    \right)
    -
    \lambda_2
    \right]
    +
    \widehat{\rho}
    \left(
    \pd{\entropy(\widehat{\temp}, \widehat{\rho})}{\widehat{\temp}}
    -
    \lambda_1
    \pd{\ienergy(\widehat{\temp}, \widehat{\rho})}{\widehat{\temp}}
    \right)
    &=
      0,
    \\
    \label{eq:313}
    \pd{\rho(\widehat{\temp}, \widehat{\thpressure})}{\widehat{\thpressure}}
    \left[
    \left(
      \widehat{\entropy}
      -
      \lambda_1
      \widehat{\ienergy}
    \right)
    +
    \widehat{\rho}
    \left(
    \pd{\entropy(\widehat{\temp}, \widehat{\rho})}{\widehat{\rho}}
      -
      \lambda_1
      \pd{\ienergy(\widehat{\temp}, \widehat{\rho})}{\widehat{\rho}}
    \right)
    -
    \lambda_2
    \right]
    \:
    \hphantom{
    +
    \widehat{\rho}
    \left(
    \pd{\entropy(\widehat{\temp}, \widehat{\rho})}{\widehat{\temp}}
    -
    \lambda_1
    \pd{\ienergy(\widehat{\temp}, \widehat{\rho})}{\widehat{\temp}}
    \right)
    }
    &=
      0
      .
  \end{align}
\end{subequations}
Now we must solve this system for $\lambda_1$ and $\lambda_2$, which is straightforward. First, we need to fix the Lagrange multiplier $\lambda_1$ such that
\begin{equation}
  \label{eq:315}
  \pd{\entropy(\widehat{\temp}, \widehat{\rho})}{\widehat{\temp}}
  -
  \lambda_1
  \pd{\ienergy(\widehat{\temp}, \widehat{\rho})}{\widehat{\temp}}
  =
  0
  ,
\end{equation}
which upon using thermodynamic identities~\eqref{eq:76} yields
\begin{equation}
  \label{eq:317}
  \lambda_1 = \frac{1}{\widehat{\temp}}.
\end{equation}
(This recovers the result we have already obtained in calculating the G\^ateaux derivative in the temperature--density representation.) Next we need to solve the equation
\begin{equation}
  \label{eq:318}
      \left(
      \widehat{\entropy}
      -
      \lambda_1
      \widehat{\ienergy}
    \right)
    +
    \widehat{\rho}
    \left(
    \pd{\entropy(\widehat{\temp}, \widehat{\rho})}{\widehat{\rho}}
      -
      \lambda_1
      \pd{\ienergy(\widehat{\temp}, \widehat{\rho})}{\widehat{\rho}}
    \right)
    -
    \lambda_2
    =
    0
    .
\end{equation}
Since we have already identified the Lagrange multiplier $\lambda_1$ as $\lambda_1 = \frac{1}{\widehat{\temp}}$, we see that
\begin{subequations}
  \label{eq:321}
  \begin{align}
    \label{eq:322}
    \widehat{\entropy}
    -
    \lambda_1
    \widehat{\ienergy}
    &=
      \widehat{\entropy}
      -
      \frac{1}{\widehat{\temp}}
      \widehat{\ienergy}
      =
      -
      \frac{1}{\widehat{\temp}}
      \widehat{\fenergy},
    \\
    \label{eq:323}
    \pd{\entropy(\widehat{\temp}, \widehat{\rho})}{\widehat{\rho}}
    -
    \lambda_1
    \pd{\ienergy(\widehat{\temp}, \widehat{\rho})}{\widehat{\rho}}
    &=
      \pd{\entropy(\widehat{\temp}, \widehat{\rho})}{\widehat{\rho}}
      -
      \frac{1}{\widehat{\temp}}
      \pd{\ienergy(\widehat{\temp}, \widehat{\rho})}{\widehat{\rho}}
      =
      -
      \frac{1}{\widehat{\temp}}
      \left(
      \pd{\fenergy(\widehat{\temp}, \widehat{\rho})}{\widehat{\rho}}
      \right)
      =
      -
      \frac{\widehat{\thpressure}}{\widehat{\rho}^2\widehat{\temp}}
      ,
  \end{align}
\end{subequations}
where we have used the relation between the Helmholtz free energy and the thermodynamic pressure~\eqref{eq:5}. Observations~\eqref{eq:321} allow us to rewrite the equation~\eqref{eq:318} as
\begin{equation}
  \label{eq:324}
  -
  \frac{\widehat{\fenergy}}{\widehat{\temp}}
  -
  \frac{\widehat{\thpressure}}{\widehat{\rho} \widehat{\temp}}
  -
  \lambda_2
  =
  0,
\end{equation}
which shows that we need to fix the second Lagrange multiplier~$\lambda_2$ as
\begin{equation}
  \label{eq:89}
  \lambda_2
  =
  -
  \frac{\widehat{\thpressure}}{\widehat{\rho} \widehat{\temp}}
  -
  \frac{\widehat{\fenergy}}{\widehat{\temp}}
  .
\end{equation}
Let us summarise our findings concerning the maximisation procedure, see Summary~\ref{summary:lagrange-multipliers}. We emphasise that Lagrange multipliers have been found without using the mass conservation constraint in the integral form.

\begin{summary}[Entropy maximisation with respect to all available constraints]
  \label{summary:lagrange-multipliers}
  If we want the functional
  \begin{equation}
    \label{eq:40}
    \mathcal{V}_{\mathrm{eq}}
    =_{\bydefinition}
    -
    \left(
      \netentropy
      -
      \netentropy_{\max}
    \right)
    +
    \lambda_1
    \left( \nettenergy - \widehat{\nettenergy} \right)
    +
    \lambda_2
    \int_{\Omega}
    \left(
      \rho
      -
      \widehat{\rho}
    \right)
    \,
    \cvolumee
  \end{equation}
  to have zero G\^ateaux derivative at the spatially homogeneous rest state~\eqref{eq:spatially-homogenoeus-rest-state}, then we must set the Lagrange multipliers $\lambda_1$ and $\lambda_2$ as
  \begin{subequations}
    \label{eq:69}
    \begin{align}
      \label{eq:73}
      \lambda_1 &= \frac{1}{\widehat{\temp}}, \\
      \label{eq:263}
      \lambda_2 &= - \frac{\widehat{\thpressure}}{\widehat{\temp} \widehat{\rho}} - \frac{\widehat{\fenergy}}{\widehat{\temp}}.
    \end{align}
   
  \end{subequations}
  Out of the functional $\mathcal{V}_{\mathrm{eq}}$ with the properly adjusted Lagrange multipliers we can construct (by simple multiplication) another functional
  \begin{equation}
    \label{eq:270}
    \mathcal{V}_{\mathrm{meq}, \, \widehat{\temp}, \, \widehat{\rho}}(\rho, \vec{v}, \temp)
    =_{\bydefinition}
    \widehat{\temp}
    \mathcal{V}_{\mathrm{eq}},
  \end{equation}
  which if written explicitly reads
  \begin{multline}
    \label{eq:271}
    \mathcal{V}_{\mathrm{meq}, \, \widehat{\temp}, \, \widehat{\rho}}(\rho, \vec{v}, \temp)
    =
    -
    \int_{\Omega} \widehat{\temp} \left( \rho \entropy(\temp, \rho) - \widehat{\rho} \entropy(\widehat{\temp}, \widehat{\rho}) \right) \, \cvolumee
    +
    \int_{\Omega} \left(\frac{1}{2} \rho \absnorm{\vec{v}}^2 + \rho \ienergy(\temp, \rho) - \widehat{\rho} \ienergy(\widehat{\temp}, \widehat{\rho}) \right) \, \cvolumee
    \\
    -
    \int_{\Omega} \left( \frac{\thpressure(\widehat{\temp}, \widehat{\rho})}{\widehat{\rho}} + \fenergy(\widehat{\temp}, \widehat{\rho}) \right) \left(\rho  - \widehat{\rho} \right)  \, \cvolumee
    .
  \end{multline}
  The time derivative of the functional $\mathcal{V}_{\mathrm{meq}, \, \widehat{\temp}, \, \widehat{\rho}}(\rho, \vec{v}, \temp)$ reads
  \begin{equation}
    \label{eq:269}
    \dd{\mathcal{V}_{\mathrm{meq}, \, \widehat{\temp}, \, \widehat{\rho}}}{t}
    =
    -
    \widehat{\temp}
    \int_{\Omega}
    \left(
      \frac{
        \tilde{\lambda} \left(\divergence \vec{v}\right)^2
        +
        2 \nu \tensordot{\traceless{\gradsym}}{\traceless{\gradsym}}
      }
      {
        \temp
      }
      +
      \frac{
        \kappa
        \vectordot{\nabla \temp}{\nabla{\temp}}
      }
      {
        \temp^2
      }
    \right)
    \, \cvolumee
    ,
  \end{equation}
  while on the right-hand side we see the entropy production multiplied by the spatially homogeneous rest state temperature $\widehat{\temp}$.
  
\end{summary}

\subsection{Nonnegativity of Lyapunov type functional}
\label{sec:nonnegativity}
Having obtained the functional $\mathcal{V}_{\mathrm{meq}, \, \widehat{\temp}, \, \widehat{\rho}}(\rho, \vec{v}, \temp)$, see~\eqref{eq:68},
\begin{multline}
  \label{eq:91}
  \mathcal{V}_{\mathrm{meq}, \, \widehat{\temp}, \, \widehat{\rho}}(\rho, \vec{v}, \temp)
    =
    -
    \int_{\Omega} \widehat{\temp} \left( \rho \entropy(\temp, \rho) - \widehat{\rho} \entropy(\widehat{\temp}, \widehat{\rho}) \right) \, \cvolumee
    +
    \int_{\Omega} \left(\frac{1}{2} \rho \absnorm{\vec{v}}^2 + \rho \ienergy(\temp, \rho) - \widehat{\rho} \ienergy(\widehat{\temp}, \widehat{\rho}) \right) \, \cvolumee
    \\
    -
    \int_{\Omega} \left( \frac{\thpressure(\widehat{\temp}, \widehat{\rho})}{\widehat{\rho}} + \fenergy(\widehat{\temp}, \widehat{\rho}) \right) \left(\rho  - \widehat{\rho} \right)  \, \cvolumee,
  \end{multline}
  we would like to show that the functional is nonnegative and that it vanishes if and only if the density, temperature and velocity fields correspond to the spatially homogeneous rest state~\eqref{eq:spatially-homogenoeus-rest-state}. We thus want to show that
  \begin{equation}
    \label{eq:38}
    \mathcal{V}_{\mathrm{meq}, \, \widehat{\temp}, \, \widehat{\rho}}(\rho, \vec{v}, \temp) \geq 0,
  \end{equation}
  while $\mathcal{V}_{\mathrm{meq}, \, \widehat{\temp}, \, \widehat{\rho}}(\rho, \vec{v}, \temp) = 0$ if and only if $\left[\rho, \vec{v}, \temp \right] = \left[ \widehat{\rho}, \widehat{\vec{v}}, \widehat{\temp} \right]$.


  \subsubsection{Proof of nonnegativity}
  \label{sec:proof-nonnegativity}
  The integrand in~\eqref{eq:91} is easy to analyse provided that we rewrite the entropy $\entropy$ and the internal energy $\ienergy$ in terms of Helmholtz free energy. (The Helmholtz free energy has the temperature $\temp$ and the density~$\rho$ as its natural variables, which indicates that the Helmholtz free energy might be the most convenient thermodynamic potential for our task.) We have
\begin{subequations}
  \label{eq:92}
  \begin{align}
    \label{eq:93}
    \entropy (\temp, \rho) &= - \pd{\fenergy}{\temp}(\temp, \rho), \\
    \label{eq:94}
    \ienergy (\temp, \rho) &=  \fenergy(\temp, \rho) + \temp \entropy(\temp, \rho),
  \end{align}
\end{subequations}
and consequently also
\begin{equation}
  \label{eq:95}
  \ienergy (\temp, \rho) =  \fenergy(\temp, \rho) - \temp \pd{\fenergy}{\temp}(\temp, \rho).
\end{equation}
Furthermore we recall that we want the specific heat capacity at constant volume $\cheatvol$ and the thermodynamic pressure $\thpressure$ to have the properties
\begin{subequations}
  \label{eq:96}
  \begin{align}
    \label{eq:97}
    \cheatvol \left(\temp, \rho\right) &> 0, \\
    \label{eq:98}
    \pd{\thpressure}{\rho}(\temp, \rho) &> 0.
  \end{align}
\end{subequations}
These are the so-called thermodynamic stability conditions, which in fact place restrictions on the second derivative of the Helmholtz free energy. As shown in Section~\ref{sec:linearised-setting} these conditions are in fact conditions on the stability of the linearised system of governing equations; see also~\cite{dostalk.m.prusa.v:non-linear} for a detailed discussion of the origin of these classical stability conditions.

Substituting \eqref{eq:92} and \eqref{eq:95} into~\eqref{eq:91} yields
\begin{multline}
  \label{eq:99}
  \mathcal{V}_{\mathrm{meq}, \, \underline{\widehat{\temp}},  \, \underline{\widehat{\rho}}}(\rho, \vec{v}, \temp)
  =
  \int_{\Omega}
  \left[
    \underline{\widehat{\temp}} \rho \pd{\fenergy}{\temp}(\temp, \rho)
    +
    \rho
    \left(
      \fenergy(\temp, \rho) - \temp \pd{\fenergy}{\temp}(\temp, \rho)
    \right)
    -
    \left(
      \frac{\thpressure(\underline{\widehat{\temp}}, \underline{\widehat{\rho}})}{\underline{\widehat{\rho}}}
      +
      \fenergy(\underline{\widehat{\temp}}, \underline{\widehat{\rho}})
    \right)
    \left(\rho  - \underline{\widehat{\rho}} \right)
  \right]
  \, \cvolumee
  \\
  -
  \int_{\Omega}
  \underline{\widehat{\rho}}
  \left[
    \ienergy(\underline{\widehat{\temp}}, \underline{\widehat{\rho}}) - \underline{\widehat{\temp}} \entropy(\underline{\widehat{\temp}}, \underline{\widehat{\rho}})
  \right]
  \, \cvolumee
  +
  \int_{\Omega} \frac{1}{2} \rho \absnorm{\vec{v}}^2 \, \cvolumee,
\end{multline}
which can be rewritten as
\begin{equation}
  \label{eq:100}
  \mathcal{V}_{\mathrm{meq}, \, \underline{\widehat{\temp}}, \, \underline{\widehat{\rho}}}(\rho, \vec{v}, \temp)
  =
  \int_{\Omega}
  \left[
    \rho
    \left(
      \fenergy\left(\temp, \rho\right)
      +
      \pd{\fenergy}{\temp}(\temp, \rho)
      \left(
        \underline{\widehat{\temp}}
        -
        \temp
      \right)
    \right)
    -
    \left(
      \frac{\thpressure(\underline{\widehat{\temp}}, \underline{\widehat{\rho}})}{\underline{\widehat{\rho}}}
      +
      \fenergy(\underline{\widehat{\temp}}, \underline{\widehat{\rho}})
    \right)
    \left(\rho  - \underline{\widehat{\rho}} \right)
  \right]
  \, \cvolumee
  -
  \int_{\Omega}
  \underline{\widehat{\rho}} \fenergy(\underline{\widehat{\temp}}, \underline{\widehat{\rho}})
  \, \cvolumee
  +
  \int_{\Omega} \frac{1}{2} \rho \absnorm{\vec{v}}^2 \, \cvolumee,
\end{equation}
where we have temporarily used the notation $\underline{\widehat{\temp}}$ and $\underline{\widehat{\rho}}$ instead of $\widehat{\temp}$ and $\widehat{\rho}$ to emphasise that the fields $\underline{\widehat{\temp}}$ and $\underline{\widehat{\rho}}$ in $\mathcal{V}_{\mathrm{meq}, \, \underline{\widehat{\temp}}, \, \underline{\widehat{\rho}}}$ enter the functional as parameters, the functional itself acts only on the fields $\rho$, $\vec{v}$ and $\temp$. Note also that the formula for the thermodynamic pressure~\eqref{eq:5} implies that~\eqref{eq:100} can be rewritten using the derivatives of the Helmholtz free energy~$\fenergy$.

We now investigate the \emph{integrand} in the first integral in~\eqref{eq:100}. Note that since we now work with the integrand only, we can not simply say that the last term vanishes because of the mass conservation. The term
\begin{equation}
  \label{eq:101}
  \left(
    \frac{\thpressure(\underline{\widehat{\temp}}, \underline{\widehat{\rho}})}{\underline{\widehat{\rho}}}
    +
    \fenergy(\underline{\widehat{\temp}}, \underline{\widehat{\rho}})
  \right)
  \left(\rho  - \underline{\widehat{\rho}} \right) 
\end{equation}
vanishes only if it is integrated over the whole domain, it \emph{does not} necessarily vanish in a pointwise sense. This is the point where we capitalise on our effort regarding the identification of the Lagrange multiplier $\lambda_2$, see~\eqref{eq:89}.

We first deal with the term
\begin{equation}
  \label{eq:102}
  \rho
  \left[
    \fenergy\left(\temp, \rho\right)
    +
    \pd{\fenergy}{\temp}(\temp, \rho)
    \left(
      \underline{\widehat{\temp}}
      -
      \temp
    \right)
  \right]
  .
\end{equation}
Taking the partial derivative of the function in the square bracket with respect to $\temp$ gives us
\begin{equation}
  \label{eq:103}
  \pd{}{\temp}
  \left[
    \fenergy\left(\temp, \rho\right)
    +
    \pd{\fenergy}{\temp}(\temp, \rho)
    \left(
      \underline{\widehat{\temp}}
      -
      \temp
    \right)
  \right]
  =
  \ppd{\fenergy}{\temp}(\temp, \rho)
  \left(
    \underline{\widehat{\temp}}
    -
    \temp
  \right)
  =
  -
  \frac{\cheatvol (\temp, \rho)}{\temp}
   \left(
     \underline{\widehat{\temp}}
     -
     \temp
   \right)
   ,
\end{equation}
where we have used the definition of the specific heat capacity at constant volume,
\begin{equation}
  \label{eq:104}
  \cheatvol(\temp, \rho)
  =
  -
  \temp
  \ppd{\fenergy}{\temp}(\temp, \rho)
  .
\end{equation}
Since the specific heat capacity at constant volume $\cheatvol$ is always positive, see~\eqref{eq:97}, we see that~\eqref{eq:103} implies that the function in the square bracket in~\eqref{eq:102}, interpreted as a function of $\temp$ only, has a strict global minimum at the point $\temp =  \underline{\widehat{\temp}}$ regardless of the value of $\rho$. We thus have
\begin{equation}
  \label{eq:105}
  \rho
  \left[
    \fenergy\left(\temp, \rho\right)
    +
    \pd{\fenergy}{\temp}(\temp, \rho)
    \left(
      \underline{\widehat{\temp}}
      -
      \temp
    \right)
  \right]
  \geq
  \rho
  \left[
    \fenergy\left(\temp, \rho\right)
    +
    \pd{\fenergy}{\temp}(\temp, \rho)
    \left(
      \underline{\widehat{\temp}}
      -
      \temp
    \right)
  \right]_{\temp = \underline{\widehat{\temp}}}
  =
  \rho
  \fenergy\left(\underline{\widehat{\temp}}, \rho\right)
  .
\end{equation}
Concerning the complete integrand in~\eqref{eq:100} we can thus write
\begin{equation}
  \label{eq:106}
  \rho
  \left[
    \fenergy\left(\temp, \rho\right)
    +
    \pd{\fenergy}{\temp}(\temp, \rho)
    \left(
      \underline{\widehat{\temp}}
      -
      \temp
    \right)
  \right]
  -
  \left(
    \frac{\thpressure(\underline{\widehat{\temp}}, \underline{\widehat{\rho}})}{\underline{\widehat{\rho}}}
    +
    \fenergy(\underline{\widehat{\temp}}, \underline{\widehat{\rho}})
  \right)
  \left(\rho  - \underline{\widehat{\rho}} \right)
  \geq
  \rho
  \fenergy\left(\underline{\widehat{\temp}}, \rho\right)
   -
  \left(
    \frac{\thpressure(\underline{\widehat{\temp}}, \underline{\widehat{\rho}})}{\underline{\widehat{\rho}}}
    +
    \fenergy(\underline{\widehat{\temp}}, \underline{\widehat{\rho}})
  \right)
  \left(\rho  - \underline{\widehat{\rho}} \right)
  .
\end{equation}

Now we investigate the function $\rho \fenergy\left(\underline{\widehat{\temp}}, \rho\right)$ as a function of $\rho$. Its \emph{second} derivative with respect to $\rho$ reads
\begin{equation}
  \label{eq:107}
  \ppd{}{\rho}
  \left[
    \rho \fenergy\left(\underline{\widehat{\temp}}, \rho\right)
  \right]
  =
  \pd{}{\rho}
  \left[
    \fenergy\left(\underline{\widehat{\temp}}, \rho\right)
    +
    \rho
    \pd{\fenergy\left(\underline{\widehat{\temp}}, \rho\right)}{\rho}
  \right]
  =
  2
  \pd{\fenergy\left(\underline{\widehat{\temp}}, \rho\right)}{\rho}
  +
  \rho
  \ppd{\fenergy\left(\underline{\widehat{\temp}}, \rho\right)}{\rho}
  =
  2
  \frac{\thpressure\left(\underline{\widehat{\temp}}, \rho\right)}{\rho^2}
  +
  \rho
  \pd{}{\rho}
  \left(
    \frac{\thpressure\left(\underline{\widehat{\temp}}, \rho\right)}{\rho^2}
  \right)
  =
  \frac{1}{\rho}
  \pd{\thpressure}{\rho}(\underline{\widehat{\temp}}, \rho)
  ,
\end{equation}
where we have used the standard formula for the thermodynamic pressure
\begin{equation}
  \label{eq:108}
  \thpressure \left( \temp, \rho \right) = \rho^2 \pd{\fenergy}{\rho}\left( \temp, \rho \right).
\end{equation}
Using the stability condition~\eqref{eq:98} we thus see that
\begin{equation}
  \label{eq:109}
  \ppd{}{\rho}
  \left[
    \rho \fenergy\left(\underline{\widehat{\temp}}, \rho\right)
  \right]
  =
  \frac{1}{\rho}
  \pd{\thpressure}{\rho}(\underline{\widehat{\temp}}, \rho)
  >
  0,
\end{equation}
hence the function $\rho \fenergy\left(\underline{\widehat{\temp}}, \rho\right)$, as a function of $\rho$, is strictly convex. A strictly convex function $f(x)$ satisfies for all $z$ and $y$ the inequality $f(z) \geq f(y) + \left. \dd{f}{x} \right|_{x = y}(z -y)$, and if we apply this characterisation of the convex function to $\rho \fenergy\left(\underline{\widehat{\temp}}, \rho\right)$, we get
\begin{equation}
  \label{eq:110}
  \rho \fenergy\left(\underline{\widehat{\temp}}, \rho\right)
  \geq
  \underline{\widehat{\rho}} \fenergy\left(\underline{\widehat{\temp}}, \underline{\widehat{\rho}} \right)
  +
  \left[
    \fenergy\left(\underline{\widehat{\temp}}, \rho\right)
    +
    \rho
    \pd{\fenergy\left(\underline{\widehat{\temp}}, \rho\right)}{\rho}
  \right]_{\rho = \underline{\widehat{\rho}}}
  \left(
    \rho
    -
    \underline{\widehat{\rho}}
  \right).
\end{equation}
Using the definition of the thermodynamic pressure we can further rewrite the inequality~\eqref{eq:110} as 
\begin{equation}
  \label{eq:113}
  \rho \fenergy\left(\underline{\widehat{\temp}}, \rho\right)
  \geq
  \underline{\widehat{\rho}} \fenergy\left(\underline{\widehat{\temp}}, \underline{\widehat{\rho}} \right)
  +
  \left(
    \fenergy\left(\underline{\widehat{\temp}},  \underline{\widehat{\rho}} \right)
    +
    \frac{\thpressure(\underline{\widehat{\temp}}, \underline{\widehat{\rho}})}{\underline{\widehat{\rho}}}
  \right)
  \left(
    \rho
    -
    \underline{\widehat{\rho}}
  \right)
  ,
\end{equation}
and exploiting this inequality in~\eqref{eq:106} finally gives us
\begin{equation}
  \label{eq:114}
  \rho
  \left[
    \fenergy\left(\temp, \rho\right)
    +
    \pd{\fenergy}{\temp}(\temp, \rho)
    \left(
      \underline{\widehat{\temp}}
      -
      \temp
    \right)
  \right]
  -
  \left(
    \frac{\thpressure(\underline{\widehat{\temp}}, \underline{\widehat{\rho}})}{\underline{\widehat{\rho}}}
    +
    \fenergy(\underline{\widehat{\temp}}, \underline{\widehat{\rho}})
  \right)
  \left(\rho  - \underline{\widehat{\rho}} \right)
  \geq
  \underline{\widehat{\rho}}  \fenergy \left(\underline{\widehat{\temp}}, \underline{\widehat{\rho}} \right)
  ,
\end{equation}
which leads to the \emph{pointwise} inequality
\begin{equation}
  \label{eq:26}
  \rho
  \left[
    \fenergy\left(\temp, \rho\right)
    +
    \pd{\fenergy}{\temp}(\temp, \rho)
    \left(
      \underline{\widehat{\temp}}
      -
      \temp
    \right)
  \right]
  -
  \left(
    \frac{\thpressure(\underline{\widehat{\temp}}, \underline{\widehat{\rho}})}{\underline{\widehat{\rho}}}
    +
    \fenergy(\underline{\widehat{\temp}}, \underline{\widehat{\rho}})
  \right)
  \left(\rho  - \underline{\widehat{\rho}} \right)
  -
  \underline{\widehat{\rho}}  \fenergy \left(\underline{\widehat{\temp}}, \underline{\widehat{\rho}} \right)
  \geq
  0
  .
\end{equation}
However, the left-hand side is precisely the combined integrand of the first two integrals in the functional $\mathcal{V}_{\mathrm{meq}, \, \underline{\widehat{\temp}}, \, \underline{\widehat{\rho}}}(\rho, \vec{v}, \temp)$, see~\eqref{eq:100}. Consequently, concerning the functional $\mathcal{V}_{\mathrm{meq}, \, \underline{\widehat{\temp}}, \, \underline{\widehat{\rho}}}(\rho, \vec{v}, \temp)$ we get the following inequality
\begin{multline}
  \label{eq:115}
  \mathcal{V}_{\mathrm{meq}, \, \underline{\widehat{\temp}}, \, \underline{\widehat{\rho}}}(\rho, \vec{v}, \temp)
  =
  \int_{\Omega}
  \left(
    \rho
    \left[
      \fenergy\left(\temp, \rho\right)
      +
      \pd{\fenergy}{\temp}(\temp, \rho)
      \left(
        \underline{\widehat{\temp}}
        -
        \temp
      \right)
    \right]
    -
    \left(
      \frac{\thpressure(\underline{\widehat{\temp}}, \underline{\widehat{\rho}})}{\underline{\widehat{\rho}}}
      +
      \fenergy(\underline{\widehat{\temp}}, \underline{\widehat{\rho}})
    \right)
    \left(\rho  - \underline{\widehat{\rho}} \right)
    -
    \underline{\widehat{\rho}}  \fenergy \left(\underline{\widehat{\temp}}, \underline{\widehat{\rho}} \right)
  \right)
  \, \cvolumee
  +
  \int_{\Omega} \frac{1}{2} \rho \absnorm{\vec{v}}^2 \, \cvolumee
  \geq
  0
  .
\end{multline}
Note that the inequality is independent of the choice of $\underline{\widehat{\temp}}$ and $\underline{\widehat{\rho}}$ as long as these quantities are constants. In particular, we can choose $\underline{\widehat{\temp}} \equiv \widehat{\temp}$ and $\underline{\widehat{\rho}} \equiv \widehat{\rho}$.

We also note that if we were dealing with the calorically perfect ideal gas, we would need no extra normalisation property regarding the choice of $\rhoref$ and $\tempref$, the pointwise inequality~\eqref{eq:26} is insensitive to the choice of normalisation constants $\rhoref$ and $\tempref$. Finally, we recall that the terms in~\eqref{eq:115} that are constant multiples of $\rho - \widehat{\rho}$ vanish when integrated over the domain $\Omega$. (This is a consequence of balance of mass.) Therefore these terms are not necessary for the validity of \emph{integral} inequality~\eqref{eq:115}. However, they are essential for the validity of \emph{pointwise} inequality~\eqref{eq:26}.

\subsubsection{G\^ateaux derivatives}
\label{sec:gateaux-derivatives}
Furthermore, we calculate the G\^ateaux derivatives of the functional $ \mathcal{V}_{\mathrm{meq}, \, \widehat{\temp}, \, \widehat{\rho}}$, see~\eqref{eq:99}, at the spatially homogeneous rest state, and we do the calculation in the temperature--density representation. We already know that the first G\^ateaux derivative of this functional at point $\widehat{\rho}$, $\widehat{\vec{v}}$ and $\widehat{\temp}$ (spatially homogeneous rest state, $\widehat{\vec{v}} = \vec{0}$) vanishes in an arbitrary direction $\widetilde{\rho}$, $\widetilde{\vec{v}}$ and $\widetilde{\temp}$. This is guaranteed by the identification of Lagrange multipliers in Section~\ref{sec:ident-lagr-mult}. The leading order non-trivial terms in the functional~\eqref{eq:91} are thus quadratic in the perturbation $\widetilde{\rho}$, $\widetilde{\vec{v}}$ and $\widetilde{\temp}$.

Since we have the balance of mass, we can in~\eqref{eq:91} ignore all terms linear in $\rho - \widehat{\rho}$. (Note that now we are operating at the level of the functionals, not at the level of pointwise inequalities for the integrand. This contrasts with the previous discussion on the nonnegativity of the integrand.) Furthermore, the G\^ateaux derivatives are easy to find for the kinetic energy term, hence we focus on calculation of G\^ateaux derivatives for the core part of the functional,
\begin{equation}
  \label{eq:116}
  \mathcal{V}_{\mathrm{meq}, \, \mathrm{core}, \, \underline{\widehat{\temp}}, \, \underline{\widehat{\rho}} }(\temp, \rho)
  =_{\bydefinition}
  \int_{\Omega}
    \rho
    \left[
      \fenergy\left(\temp, \rho\right)
      +
      \pd{\fenergy}{\temp}(\temp, \rho)
      \left(
        \underline{\widehat{\temp}}
        -
        \temp
      \right)
    \right]
  \, \cvolumee
  .
\end{equation}
Now we are ready to calculate the first and the second G\^ateaux derivative of this functional. The definition of the first and the second G\^ateaux derivative reads
\begin{subequations}
  \label{eq:117}
  \begin{align}
    \label{eq:118}
      \Diff \mathcal{V}_{\mathrm{meq}, \, \mathrm{core}, \, \underline{\widehat{\temp}}, \, \underline{\widehat{\rho}}}(\widehat{\temp}, \widehat{\rho})
  \left[
    \begin{bmatrix}
      \widetilde{\temp} \\
      \widetilde{\rho}
    \end{bmatrix}
  \right]
  &=
  \left.
    \dd{}{s}
    \left(
       \mathcal{V}_{\mathrm{meq}, \, \mathrm{core}, \, \underline{\widehat{\temp}}, \, \underline{\widehat{\rho}}} \left(\widehat{\temp} + s \widetilde{\temp}, \widehat{\rho} + s \widetilde{\rho} \right)
    \right)
  \right|_{s=0}
    ,
    \\
    \label{eq:119}
    \Diff^2 \mathcal{V}_{\mathrm{meq}, \, \mathrm{core}, \, \underline{\widehat{\temp}}, \, \underline{\widehat{\rho}}}(\widehat{\temp}, \widehat{\rho})
  \left[
    \begin{bmatrix}
      \widetilde{\temp} \\
      \widetilde{\rho}
    \end{bmatrix}
    ,
    \begin{bmatrix}
      \widetilde{\temp} \\
      \widetilde{\rho}
    \end{bmatrix}
  \right]
  &=
  \left.
    \ddd{}{s}
    \left(
       \mathcal{V}_{\mathrm{meq}, \, \mathrm{core}, \, \underline{\widehat{\temp}}, \, \underline{\widehat{\rho}}} \left(\widehat{\temp} + s \widetilde{\temp}, \widehat{\rho} + s \widetilde{\rho} \right)
    \right)
  \right|_{s=0}
    .
  \end{align}
\end{subequations}
As we have already noted, the first G\^ateaux derivative vanishes by construction. Concerning the second derivative we see that
\begin{equation}
  \label{eq:120}
  \underline{\widehat{\temp}} \rho \pd{\fenergy}{\temp}(\widehat{\temp} + s \widetilde{\temp}, \widehat{\rho} + s \widetilde{\rho})
  =
  \cdots
  +
  \underline{\widehat{\temp}} \widehat{\rho}
  \frac{1}{2}
  \left(
    \pd{^3\fenergy}{\temp^3}(\widehat{\temp}, \widehat{\rho})
    \widetilde{\temp}^2
    +
    2
    \pd{^3\fenergy}{\temp^2 \partial \rho}(\widehat{\temp}, \widehat{\rho})
    \widetilde{\temp} \widetilde{\rho}
    +
    \pd{^3\fenergy}{\temp \partial \rho^2}(\widehat{\temp}, \widehat{\rho})
    \widetilde{\rho}^2
   \right)
   s^2
   +
   \underline{\widehat{\temp}} \widetilde{\rho}
   \left(
     \pd{^2\fenergy}{\temp^2}(\widehat{\temp}, \widehat{\rho})
     \widetilde{\temp}
     +
     \pd{^2\fenergy}{\temp \partial \rho}(\widehat{\temp}, \widehat{\rho})
     \widetilde{\rho}
   \right)
   s^2
   +
   \cdots
   ,
\end{equation}
where we include only the second order terms. Furthermore,
\begin{multline}
  \label{eq:121}
  \rho
  \left(
    \fenergy(\temp, \rho) - \temp \pd{\fenergy}{\temp}(\temp, \rho)
  \right)
  =
  \cdots
  \\
  +
  \widehat{\rho}
  \left[
    \frac{1}{2}
    \left(
      \pd{^2\fenergy}{\temp^2}(\widehat{\temp}, \widehat{\rho})
      \widetilde{\temp}^2
      +
      2
      \pd{^2\fenergy}{\temp \partial \rho}(\widehat{\temp}, \widehat{\rho})
      \widetilde{\temp} \widetilde{\rho}
      +
      \pd{^2\fenergy}{\rho^2}(\widehat{\temp}, \widehat{\rho})
      \widetilde{\rho}^2
    \right)
    -
    \frac{1}{2}
    \widehat{\temp}
      \left(
    \pd{^3\fenergy}{\temp^3}(\widehat{\temp}, \widehat{\rho})
    \widetilde{\temp}^2
    +
    2
    \pd{^3\fenergy}{\temp^2 \partial \rho}(\widehat{\temp}, \widehat{\rho})
    \widetilde{\temp} \widetilde{\rho}
    +
    \pd{^3\fenergy}{\temp \partial \rho^2}(\widehat{\temp}, \widehat{\rho})
    \widetilde{\rho}^2
   \right)
 \right]
 s^2
 \\
 -
 \widehat{\rho}
 \widetilde{\temp}
 \left(
   \pd{^2\fenergy}{\temp^2}(\widehat{\temp}, \widehat{\rho})
   \widetilde{\temp}
   +
   \pd{^2\fenergy}{\temp \partial \rho}(\widehat{\temp}, \widehat{\rho})
   \widetilde{\rho}
 \right)
 s^2
 \\
 +
 \widetilde{\rho}
 \left[
   \left(
     \pd{\fenergy}{\temp}(\widehat{\temp}, \widehat{\rho})
     \widetilde{\temp}
     +
     \pd{\fenergy}{\rho}(\widehat{\temp}, \widehat{\rho})
     \widetilde{\rho}
   \right)
   -
   \widehat{\temp}
   \left(
     \pd{^2\fenergy}{\temp^2}(\widehat{\temp}, \widehat{\rho})
     \widetilde{\temp}
     +
     \pd{^2\fenergy}{\temp \partial \rho }(\widehat{\temp}, \widehat{\rho})
     \widetilde{\rho}
   \right)
 \right]
 s^2
 -
 \widetilde{\rho}
 \widetilde{\temp}
 \pd{\fenergy}{\temp}(\widehat{\temp}, \widehat{\rho})
 s^2
 .
\end{multline}
Using the just derived formulae, we see that
\begin{multline}
  \label{eq:122}
  \mathcal{V}_{\mathrm{meq}, \, \mathrm{core}, \, \underline{\widehat{\temp}}, \, \underline{\widehat{\rho}}}  \left(\widehat{\temp} + s \widetilde{\temp}, \widehat{\rho} + s \widetilde{\rho} \right)
  =
  \cdots
  +
  \int_{\Omega}
  \frac{1}{2}
  \widehat{\rho}
  \left(
    \underline{\widehat{\temp}} - \widehat{\temp}
  \right)
  \left(
    \pd{^3\fenergy}{\temp^3}(\widehat{\temp}, \widehat{\rho})
    \widetilde{\temp}^2
    +
    2
    \pd{^3\fenergy}{\temp^2 \partial \rho}(\widehat{\temp}, \widehat{\rho})
    \widetilde{\temp} \widetilde{\rho}
    +
    \pd{^3\fenergy}{\temp \partial \rho^2}(\widehat{\temp}, \widehat{\rho})
    \widetilde{\rho}^2
  \right)
  s^2
  \, \cvolumee
  \\
  +
  \int_{\Omega}
  \left(
    \underline{\widehat{\temp}} - \widehat{\temp}
  \right)
  \widetilde{\rho}
   \left(
     \pd{^2\fenergy}{\temp^2}(\widehat{\temp}, \widehat{\rho})
     \widetilde{\temp}
     +
     \pd{^2\fenergy}{\temp \partial \rho}(\widehat{\temp}, \widehat{\rho})
     \widetilde{\rho}
   \right)
   s^2  
  \, \cvolumee
  \\
  -
  \int_{\Omega}
  \frac{1}{2}
  \widehat{\rho}
  \pd{^2\fenergy}{\temp^2}(\widehat{\temp}, \widehat{\rho})
  \widetilde{\temp}^2
  s^2
  \, \cvolumee
  +
  \int_{\Omega}
  \left(
    \frac{1}{2}
    \widehat{\rho}
    \pd{^2\fenergy}{\rho^2}(\widehat{\temp}, \widehat{\rho})
    +
    \pd{\fenergy}{\rho}(\widehat{\temp}, \widehat{\rho})
  \right)
  \widetilde{\rho}^2
  s^2
  \, \cvolumee
  +
  \cdots
  .
\end{multline}
Now we make use of thermodynamic identities. We know that
\begin{align}
  \label{eq:123}
  \cheatvol \left(\widehat{\temp}, \widehat{\rho} \right)
  &=
  - \widehat{\temp}  \ppd{\fenergy}{\temp}(\widehat{\temp}, \widehat{\rho})
  ,
  \\
  \label{eq:124}
  \thpressure \left(\widehat{\temp}, \widehat{\rho} \right)
  &=
  \widehat{\rho}^2 \pd{\fenergy}{\rho}(\widehat{\temp}, \widehat{\rho}),
\end{align}
which implies that
\begin{equation}
  \label{eq:125}
  \frac{1}{2}
  \widehat{\rho}
  \pd{^2\fenergy}{\rho^2}(\widehat{\temp}, \widehat{\rho})
  +
  \pd{\fenergy}{\rho}(\widehat{\temp}, \widehat{\rho})
  =
  \frac{1}{2} \frac{1}{\widehat{\rho}} \pd{\thpressure}{\rho}(\widehat{\temp}, \widehat{\rho}).
\end{equation}
Using these identities and \emph{evaluating the functional at the point} $\widehat{\temp} \equiv \underline{\widehat{\temp}}$ and  $\widehat{\rho} \equiv \underline{\widehat{\rho}}$ we see that the second order terms in the functional expansion are
\begin{equation}
  \label{eq:126}
  \mathcal{V}_{\mathrm{meq}, \, \mathrm{core}, \, \underline{\widehat{\temp}}, \, \underline{\widehat{\rho}}}  \left(\widehat{\temp} + s \widetilde{\temp}, \widehat{\rho} + s \widetilde{\rho} \right)
  =
  \cdots
  +
  \left[
    \int_{\Omega}
    \left(
      \frac{1}{2}
      \widehat{\rho}
      \frac{\cheatvol \left(\widehat{\temp}, \widehat{\rho}\right)}{\widehat{\temp}}
      \widetilde{\temp}^2
      +
      \frac{1}{2} \frac{1}{\widehat{\rho}} \pd{\thpressure}{\rho}(\widehat{\temp}, \widehat{\rho})
      \widetilde{\rho}^2
    \right)
    \, \cvolumee
  \right]
  s^2
  +
  \cdots
  ,
\end{equation}
hence
\begin{equation}
  \label{eq:127}
  \frac{1}{2}
  \Diff^2 \mathcal{V}_{\mathrm{meq}, \, \mathrm{core}, \, \underline{\widehat{\temp}}, \, \underline{\widehat{\rho}}}(\widehat{\temp}, \widehat{\rho})
  \left[
    \begin{bmatrix}
      \widetilde{\temp} \\
      \widetilde{\rho}
    \end{bmatrix}
    ,
    \begin{bmatrix}
      \widetilde{\temp} \\
      \widetilde{\rho}
    \end{bmatrix}
  \right]
  =
  \frac{1}{2}
  \left.
    \ddd{}{s}
    \left(
       \mathcal{V}_{\mathrm{meq}, \, \mathrm{core}, \, \underline{\widehat{\temp}}, \, \underline{\widehat{\rho}}} \left(\widehat{\temp} + s \widetilde{\temp}, \widehat{\rho} + s \widetilde{\rho} \right)
    \right)
  \right|_{s=0}
  =
  \int_{\Omega}
  \left(
    \frac{1}{2}
    \widehat{\rho}
    \frac{\cheatvol \left(\widehat{\temp}, \widehat{\rho}\right)}{\widehat{\temp}}
    \widetilde{\temp}^2
    +
    \frac{1}{2} \frac{1}{\widehat{\rho}} \pd{\thpressure}{\rho}(\widehat{\temp}, \widehat{\rho})
    \widetilde{\rho}^2
  \right)
  \, \cvolumee
  .  
\end{equation}
This is by no means surprising. The functional~\eqref{eq:127} is in fact the quadratic functional that can be obtained by the stability analysis of the \emph{linearised governing equations} in the neighbourhood of the steady state $\widehat{\temp}$ and $\widehat{\rho}$, see Section~\ref{sec:linearised-setting}, formula~\eqref{eq:57}, and also~\cite{dostalk.m.prusa.v:non-linear}.

\subsubsection{Conclusion}
\label{sec:conclusion}
The analysis in this section shows that the thermodynamic conditions
\begin{subequations}
  \label{eq:128}
  \begin{align}
    \label{eq:129}
    \cheatvol \left(\temp, \rho\right) &> 0, \\
    \label{eq:130}
    \pd{\thpressure}{\rho}(\temp, \rho) &> 0
  \end{align}
\end{subequations}
guessed by the stability analysis of the \emph{linearised} equations in fact allow us to construct the Lyapunov type functional~\eqref{eq:91} that is suitable for \emph{nonlinear} stability analysis. Our findings are summarised in Summary~\ref{summary:decay-equation-nonlinear-rest-state}.

\begin{summary}[Thermodynamic stability conditions imply nonnegativity of the functional~$\mathcal{V}_{\mathrm{meq}, \, \widehat{\temp}, \, \widehat{\rho}}(\rho, \vec{v}, \temp)$]
\label{summary:decay-equation-nonlinear-rest-state}
  If thermodynamic stability inequalities
  \begin{subequations}
    \label{eq:277}
    \begin{align}
      \label{eq:278}
    \cheatvol \left(\temp, \rho\right) &> 0, \\
      \label{eq:279}
    \pd{\thpressure}{\rho}(\temp, \rho) &> 0,
    \end{align}
  \end{subequations}
  hold for all $\temp$, $\rho$ pairs, then the functional~$\mathcal{V}_{\mathrm{meq}, \, \widehat{\temp}, \, \widehat{\rho}}(\rho, \vec{v}, \temp)$ introduced in Summary~\ref{summary:lagrange-multipliers} is nonnegative and it vanishes if and only if the system is at the spatially homogeneous rest state, that is
    \begin{equation}
      \label{eq:281}
      \forall \text{ admissible } \left[ \rho, \vec{v}, \temp \right]: \quad \mathcal{V}_{\mathrm{meq}, \, \widehat{\temp}, \, \widehat{\rho}}(\rho, \vec{v}, \temp) \geq 0,
    \end{equation}
    where $\mathcal{V}_{\mathrm{meq}, \, \widehat{\temp}, \, \widehat{\rho}}(\rho, \vec{v}, \temp) = 0$ if and only if $\left[\rho, \vec{v}, \temp \right] = \left[ \widehat{\rho}, \widehat{\vec{v}}, \widehat{\temp} \right]$. The admissibility means that the state $\left[ \rho, \vec{v}, \temp \right]$ has the same net mass and the net total energy as the spatially homogeneous rest state.

    \bigskip
    Moreover, we already know that the time derivative of the functional is nonpositive, and that it vanishes if and only if the velocity field vanishes and the temperature field is spatially homogeneous,
    \begin{equation}
      \label{eq:283}
      \dd{\mathcal{V}_{\mathrm{meq}, \, \widehat{\temp}, \, \widehat{\rho}}}{t}
      =
      -
      \widehat{\temp}
      \int_{\Omega}
      \left(
        \frac{
          \tilde{\lambda} \left(\divergence \vec{v}\right)^2
          +
          2 \nu \tensordot{\traceless{\gradsym}}{\traceless{\gradsym}}
        }
        {
          \temp
        }
        +
        \frac{
          \kappa
          \vectordot{\nabla \temp}{\nabla{\temp}}
        }
        {
          \temp^2
        }
      \right)
      \, \cvolumee
      .
    \end{equation}
    Equation~\eqref{eq:283} is a universal equation that suggests that any compatible initial state eventually decays to the spatially homogeneous rest state.
\end{summary}

\subsection{Example---Lyapunov type functional and its time derivative for calorically perfect ideal gas}
\label{sec:example-calor-perf}
The calorically perfect ideal gas is a substance with the Helmholtz free energy given by the formula~\eqref{eq:free-energy-ideal-gas-introduction}.
 Now we substitute into the formula for the functional~$\mathcal{V}_{\mathrm{meq}, \, \widehat{\temp}, \, \widehat{\rho}}(\rho, \vec{v}, \temp)$, that is to the formula
\begin{multline}
  \label{eq:138}
  \mathcal{V}_{\mathrm{meq}, \, \widehat{\temp}, \, \widehat{\rho}}(\rho, \vec{v}, \temp)
  =
  -
    \int_{\Omega} \widehat{\temp} \left( \rho \entropy(\temp, \rho) - \widehat{\rho} \entropy(\widehat{\temp}, \widehat{\rho}) \right) \, \cvolumee
    +
    \int_{\Omega} \left(\frac{1}{2} \rho \absnorm{\vec{v}}^2 + \rho \ienergy(\temp, \rho) - \widehat{\rho} \ienergy(\widehat{\temp}, \widehat{\rho}) \right) \, \cvolumee
    \\
    -
    \int_{\Omega} \left( \frac{\thpressure(\widehat{\temp}, \widehat{\rho})}{\widehat{\rho}} + \fenergy(\widehat{\temp}, \widehat{\rho}) \right) \left(\rho  - \widehat{\rho} \right)  \, \cvolumee,
  \end{multline}
and after some algebra we get
\begin{equation}
  \label{eq:139}
  \mathcal{V}_{\mathrm{meq}, \, \widehat{\temp}, \, \widehat{\rho}}(\rho, \vec{v}, \temp)
  =
  \int_{\Omega}
  \frac{1}{2}
  \rho \absnorm{\vec{v}}^2
  \, \cvolumee
  +
  \int_{\Omega}
  \rho \widehat{\temp} \cheatvolref
  \left[
    \frac{\temp}{\widehat{\temp}}
    -
    1
    -
    \ln \frac{\temp}{\widehat{\temp}}
  \right]
  \, \cvolumee
  +
  \int_{\Omega}
  \cheatvolref \left(\gamma -1\right) \widehat{\temp} \widehat{\rho}
  \left[
    \frac{\rho}{\widehat{\rho}}
    \ln
    \frac{\rho}{\widehat{\rho}}
    -
    \frac{\rho}{\widehat{\rho}}
    +
    1
  \right]
  \, \cvolumee.
\end{equation}
(Note that the normalisation constants $\rhoref$ and $\tempref$ in formulae for the calorically perfect ideal gas, see Section~\ref{sec:calor-perf-ideal}, do not enter the final formula~\eqref{eq:139}.) It is straightforward to verify that the terms in the square brackets in~\eqref{eq:139} are nonnegative, and that they vanish if and only if $\temp = \widehat{\temp}$ and $\rho = \widehat{\rho}$. Furthermore, we see that
\begin{subequations}
  \label{eq:140}
  \begin{align}
  \label{eq:141}
    \frac{\temp}{\widehat{\temp}}
    -
    1
    -
    \ln \frac{\temp}{\widehat{\temp}}
    &=
    \frac{\widetilde{\temp}}{\widehat{\temp}}
    -
    \ln
    \left(
      1
      +
      \frac{\widetilde{\temp}}{\widehat{\temp}}
    \right)
    \approx
    \frac{1}{2}
    \left(
      \frac{\widetilde{\temp}}{\widehat{\temp}}
      \right)^2,
    \\
    \label{eq:142}
    \frac{\rho}{\widehat{\rho}}
    \ln
    \frac{\rho}{\widehat{\rho}}
    -
    \frac{\rho}{\widehat{\rho}}
    +
    1
    &=
      \left(1 + \frac{\widetilde{\rho}}{\widehat{\rho}} \right)
      \ln
      \left(
      1
      +
      \frac{\widetilde{\rho}}{\widehat{\rho}}
      \right)
      -
      \frac{\widetilde{\rho}}{\widehat{\rho}}
      \approx
      \frac{1}{2}
      \left(
            \frac{\widetilde{\rho}}{\widehat{\rho}}
      \right)^2
      ,
  \end{align}
\end{subequations}
which confirms our earlier findings regarding the second G\^ateaux derivative of~$\mathcal{V}_{\mathrm{meq}, \, \widehat{\temp}, \, \widehat{\rho}}$, see, for example, formula~\eqref{eq:127}. Complete formulae for the calorically perfect ideal gas are shown in Summary~\ref{summary:decay-equation-calorically-perfect-gas}.

\begin{summary}[Decay equation in the nonlinear setting---calorically perfect ideal gas]
  \label{summary:decay-equation-calorically-perfect-gas}
   Let us consider the problem in Question~\ref{q:isolated} and let the fluid of interest be the calorically perfect ideal gas. If we work out the formulae in Summary~\ref{summary:decay-equation-nonlinear-rest-state} for the calorically perfect ideal gas, then we see that any solution to the Navier--Stokes--Fourier equations starting from a compatible initial condition~\eqref{eq:inital-condition-isolated} satisfies  
  \begin{multline}
    \label{eq:282}
    \dd{}{t}
    \int_{\Omega}
    \left(
      \frac{1}{2}
      \rho \absnorm{\vec{v}}^2
      +
      \rho \widehat{\temp} \cheatvolref
      \left[
        \frac{\temp}{\widehat{\temp}}
        -
        1
        -
        \ln \frac{\temp}{\widehat{\temp}}
      \right]
      +
      \cheatvolref \left(\gamma -1\right) \widehat{\temp} \widehat{\rho}
      \left[
        \frac{\rho}{\widehat{\rho}}
        \ln
        \frac{\rho}{\widehat{\rho}}
        -
        \frac{\rho}{\widehat{\rho}}
        +
        1
      \right]
    \right)
    \, \cvolumee
    \\
    =
    -
    \widehat{\temp}
    \int_{\Omega}
    \left(
      \frac{
        \tilde{\lambda} \left(\divergence \vec{v}\right)^2
        +
        2 \nu \tensordot{\traceless{\gradsym}}{\traceless{\gradsym}}
      }
      {
        \temp
      }
      +
      \frac{
        \kappa
        \vectordot{\nabla \temp}{\nabla{\temp}}
      }
      {
        \temp^2
      }
    \right)
    \, \cvolumee
    ,
  \end{multline}
  where
  $
  [
  \widehat{\rho}, \widehat{\vec{v}}, \widehat{\temp}
  ]
  $
  denotes the spatially homogeneous rest state~\eqref{eq:spatially-homogenoeus-rest-state}.
\end{summary}

\subsection{Revisiting the construction from the perspective of classical thermodynamics of spatially homogeneous systems}
\label{sec:revis-constr-from}
The identification of the multipliers described in Section~\ref{sec:ident-lagr-mult} is, to some extent, a tedious continuum mechanics version (spatially distributed systems) of the following formal construction known in classical thermodynamics, see, for example, \cite{callen.hb:thermodynamics} and~\cite{muller.i:thermodynamics}. (For detailed computational machinery see also~\cite{prestipino.s.giaquinta.pv:concavity}.) We consider the entropy of a system with the net total energy $\widehat{\nettenergy}$ and the volume $\widehat{\ctvolume}$. In the classical setting there is no macroscopic motion, that is $\vec{v} = \vec{0}$, and the net total energy $\nettenergy$ coincides with the net internal energy $\ctenergy$, $\nettenergy \equiv \ctenergy$. We can therefore write the entropic equation of state as 
\begin{equation}
  \label{eq:143}
  \ctentropy = \ctentropy (\ctenergy, \ctvolume).
\end{equation}
Now we consider a substance with the entropic equation of state~\eqref{eq:143}, while the substance occupies a box of volume $\ctvolume_{\text{box}}$ and energy $\ctenergy_{\text{box}}$. The box is isolated from the outside environment and it is divided into two compartments by a movable and thermally conductive piston. (This means that the two compartments can exchange energy in any form, but no energy can be exchanged with the outside environment.) The compartments have energies $\ctenergy_1$ and $\ctenergy_2$ and volumes $\ctvolume_1$ and $\ctvolume_2$, meaning that the box energy/volume $\ctenergy_{\text{box}}$/$\ctvolume_{\text{box}}$ can be arbitrarily redistributed into the compartments---both compartments are still spatially homogeneous but the box is not. We have obvious constraints
\begin{subequations}
  \label{eq:144}
  \begin{align}
    \label{eq:145}
    \ctenergy_{\text{box}} &= \ctenergy_1 + \ctenergy_2, \\
    \label{eq:146}
    \ctvolume_{\text{box}} &= \ctvolume_1 + \ctvolume_2.
  \end{align}
\end{subequations}

The task is to maximise the net entropy of the whole box $\ctentropy_{\text{box}}$ subject to all constraints. This should give us the most
favourable redistribution of the energy/volume between the compartments. The maximisation is done using the Lagrange multipliers,
\begin{equation}
  \label{eq:147}
  \ctentropy_{\text{box}} - \lambda_1 \left (\ctenergy_1 + \ctenergy_2 - \ctenergy_{\text{box}}\right)  - \lambda_2 \left(\ctvolume_1 + \ctvolume_2 - \ctvolume_{\text{box}}\right).
\end{equation}

The entropy of the box $\ctentropy_{\text{box}}$ is however the sum of the entropies of the compartments. (Entropy is an extensive quantity.) We thus have
\begin{equation}
  \label{eq:148}
  \ctentropy_{\text{box}} =  \ctentropy_1 + \ctentropy_2,
\end{equation}
which allows us to write~\eqref{eq:147} as
\begin{equation}
  \label{eq:149}
  \ctentropy(\ctenergy_1, \ctvolume_1) + \ctentropy(\ctenergy_2, \ctvolume_2) - \lambda_1 \left (\ctenergy_1 + \ctenergy_2 - \ctenergy_{\text{box}}\right)  - \lambda_2 \left(\ctvolume_1 + \ctvolume_2 - \ctvolume_{\text{box}}\right).
\end{equation}
Our constrained maximisation problem thus reads
\begin{equation}
  \label{eq:150}
  \max_{\ctenergy_1, \ctenergy_2, \ctvolume_1, \ctvolume_2} \left\{  \ctentropy(\ctenergy_1, \ctvolume_1) + \ctentropy(\ctenergy_2, \ctvolume_2) - \lambda_1 \left (\ctenergy_1 + \ctenergy_2 - \ctenergy_{\text{box}}\right)  - \lambda_2 \left(\ctvolume_1 + \ctvolume_2 - \ctvolume_{\text{box}}\right) \right\}.
\end{equation}
We first investigate the conditions for the extremum. Taking the partial derivatives with respect to $\ctenergy_1, \ctenergy_2, \ctvolume_1$ and $\ctvolume_2$ we get
\begin{subequations}
    \label{eq:151}
    \begin{align}
      \label{eq:152}
      \pd{\ctentropy}{\ctenergy}(\ctenergy_1, \ctvolume_1)  - \lambda_1 & =0, \\
      \label{eq:153}
      \pd{\ctentropy}{\ctenergy}(\ctenergy_2, \ctvolume_2)  - \lambda_1 & =0, \\
      \label{eq:154}
      \pd{\ctentropy}{\ctvolume}(\ctenergy_1, \ctvolume_1) - \lambda_2 & =0, \\
      \label{eq:155}
      \pd{\ctentropy}{\ctvolume}(\ctenergy_2, \ctvolume_2) - \lambda_2 & =0.
    \end{align}
  \end{subequations}
  Here the notation $\pd{\ctentropy}{\ctenergy}(\ctenergy_1, \ctvolume_1)$ means that we take the partial derivative of the entropic equation of state~\eqref{eq:143} and evaluate the result at $\ctvolume = \ctvolume_1$ and $\ctenergy = \ctenergy_1$, that is $\left. \pd{\ctentropy}{\ctenergy}(\ctenergy, \ctvolume) \right|_{\left[\ctenergy, \ctvolume\right] = \left[\ctenergy_1, \ctvolume_1\right]}$. From~\eqref{eq:151} we can thus easily identify the Lagrange multipliers
  \begin{subequations}
    \label{eq:156}
     \begin{align}
      \label{eq:157}
       \lambda_1 &= \pd{\ctentropy}{\ctenergy}(\ctenergy_1, \ctvolume_1), \\
       \label{eq:158}
       \lambda_1 &= \pd{\ctentropy}{\ctenergy}(\ctenergy_2, \ctvolume_2), \\
       \label{eq:159}
       \lambda_2 &= \pd{\ctentropy}{\ctvolume}(\ctenergy_1, \ctvolume_1), \\
       \label{eq:160}
       \lambda_2 &= \pd{\ctentropy}{\ctvolume}(\ctenergy_2, \ctvolume_2).
    \end{align}
  \end{subequations}
  We however know that
  \begin{equation}
    \label{eq:161}
    \pd{\ctentropy}{\ctenergy}(\ctenergy, \ctvolume) = \frac{1}{\cttemperature},
  \end{equation}
  where $\cttemperature$ denotes the thermodynamic temperature. Furthermore, we know that
  \begin{equation}
    \label{eq:162}
    \pd{\ctenergy}{\ctvolume}(\ctentropy, \ctvolume) = - \ctpressure,    
  \end{equation}
  where $\ctpressure$ denotes the thermodynamic pressure. (With the sign convention used in the classical thermodynamics.) With the standard abuse of notation we also have 
  \begin{equation}
    \label{eq:163}
    \ctenergy (\ctentropy ( \ctenergy, \ctvolume), \ctvolume ) = \ctenergy,
  \end{equation}
  which upon differentiation with respect to $\ctvolume$ at constant $\ctenergy$ gives
  \begin{equation}
    \pd{\ctenergy}{\ctentropy}(\ctentropy, \ctvolume) \pd{\ctentropy}{\ctvolume}(\ctenergy, \ctvolume) + \pd{\ctenergy}{\ctvolume}(\ctentropy, \ctvolume) = 0,
  \end{equation}
  hence
  \begin{equation}
    \label{eq:164}
    \cttemperature \pd{\ctentropy}{\ctvolume}(\ctenergy, \ctvolume) - \ctpressure = 0.
  \end{equation}
  Concerning the Lagrange multipliers we thus have
 \begin{subequations}
   \label{eq:165}
     \begin{align}
       \label{eq:166}
       \lambda_1 &= \frac{1}{\cttemperature_1}, \\
       \label{eq:167}
       \lambda_1 &= \frac{1}{\cttemperature_2}, \\
       \label{eq:168}
       \lambda_2 &= \frac{\ctpressure_1}{\cttemperature_1}, \\
       \label{eq:169}
       \lambda_2 &= \frac{\ctpressure_2}{\cttemperature_2},
    \end{align}
  \end{subequations}
  where $\cttemperature_1$, $\cttemperature_2$, $\ctpressure_1$ and $\ctpressure_2$ denote the temperature/pressure in the given compartment. Consequently, we see that the entropy attains the extremum at the state where the pressure and the temperature are the same in both compartments. This means that at the extremum entropy the box is in a spatially homogeneous state, the energy and the volume are split equally in between the compartments.

  Now we have to figure out\footnote{This is a two way procedure---we can either design our entropic equation of state in such a way that it leads to the maximum, or we want the spatially homogeneous state to be the state of maximum entropy and investigate whether this places some restrictions on the possible structure of entropic equations of state.} whether the extremum is a (local) maximum. We characterise the extremum by the second-derivative test. If we perturb the extremal state $\ctenergy_1 = \ctenergy_2 = \frac{\ctenergy_{\text{box}}}{2}$, $\ctvolume_1 = \ctvolume_2 = \frac{\ctvolume_{\text{box}}}{2}$ in a way that is compatible with the constraints~\eqref{eq:144}, that is if we set $\ctenergy_1 = \frac{\ctenergy_{\text{box}}}{2} + \delta \ctenergy$, $\ctenergy_2 = \frac{\ctenergy_{\text{box}}}{2} - \delta \ctenergy$, $\ctvolume_1 = \frac{\ctvolume_{\text{box}}}{2} + \delta \ctvolume$ and $\ctvolume_2 = \frac{\ctvolume_{\text{box}}}{2} - \delta \ctvolume$, then the first order terms cancel, and the box entropy~\eqref{eq:148} reads
  \begin{equation}
    \label{eq:170}
    \ctentropy_{\text{box}}
    =
    2 \ctentropy \left( \frac{\ctenergy_{\text{box}}}{2}, \frac{\ctvolume_{\text{box}}}{2} \right)
    +
    \begin{bmatrix}
      \delta \ctenergy & \delta \ctvolume
    \end{bmatrix}
    \left.
    \begin{bmatrix}
      \ppd{\ctentropy}{\ctenergy}(\ctenergy, \ctvolume) & \pd{^2\ctentropy}{\ctenergy \partial \ctvolume}(\ctenergy, \ctvolume) \\
      \pd{^2\ctentropy}{\ctenergy \partial \ctvolume}(\ctenergy, \ctvolume) &  \ppd{\ctentropy}{\ctvolume}(\ctenergy, \ctvolume)
    \end{bmatrix}
    \right|_{\left[\ctenergy, \ctvolume\right] = \left[\frac{\ctenergy_{\text{box}}}{2}, \frac{\ctvolume_{\text{box}}}{2}\right]}
    \begin{bmatrix}
      \delta \ctenergy \\ \delta \ctvolume
    \end{bmatrix}
    +
    \cdots
    .
  \end{equation}
  The extremum is thus a local maximum provided that the entropy is a concave function of its natural variables $\ctenergy$ and $\ctvolume$ at the compartment state $\left[\frac{\ctenergy_{\text{box}}}{2}, \frac{\ctvolume_{\text{box}}}{2}\right]$, that is provided that the second derivatives matrix
  \begin{equation}
    \label{eq:171}
    \left.
    \begin{bmatrix}
      \ppd{\ctentropy}{\ctenergy}(\ctenergy, \ctvolume) & \pd{^2\ctentropy}{\ctenergy \partial \ctvolume}(\ctenergy, \ctvolume) \\
      \pd{^2\ctentropy}{\ctenergy \partial \ctvolume}(\ctenergy, \ctvolume) &  \ppd{\ctentropy}{\ctvolume}(\ctenergy, \ctvolume)
    \end{bmatrix}
    \right|_{\left[\ctenergy, \ctvolume\right] = \left[\frac{\ctenergy_{\text{box}}}{2}, \frac{\ctvolume_{\text{box}}}{2}\right]}
  \end{equation}
  is negative definite. (The entropy is one-homogeneous in the energy, the volume and the amount of substance, while the amount of substance in each compartment is fixed. If the entropy $\ctentropy(\ctenergy, \ctvolume)$ were one-homogeneous in $\ctenergy$ and $\ctvolume$ alone, then its matrix of second derivatives would be singular.)
  The principal minors test guarantees that the matrix of second derivatives is negative definite provided that
  \begin{subequations}
    \label{eq:172}
    \begin{align}
      \label{eq:173}
      \ppd{\ctentropy}{\ctenergy}(\ctenergy, \ctvolume) &< 0, \\
      \label{eq:174}
      \ppd{\ctentropy}{\ctenergy}(\ctenergy, \ctvolume)\ppd{\ctentropy}{\ctvolume}(\ctenergy, \ctvolume) - \left( \pd{^2\ctentropy}{\ctenergy \partial \ctvolume}(\ctenergy, \ctvolume) \right)^2 &>0,
    \end{align}
  \end{subequations}
  which implies
  \begin{subequations}
    \label{eq:175}
    \begin{align}
      \label{eq:176}
      \ppd{\ctentropy}{\ctenergy}(\ctenergy, \ctvolume) &< 0, \\
      \label{eq:177}
      \ppd{\ctentropy}{\ctvolume}(\ctenergy, \ctvolume) &< 0.
    \end{align}
  \end{subequations}
  These restrictions are however difficult to interpret in terms of some directly accessible quantities. Ideally we would like to see the restrictions in the form where the independent variables are the temperature and the volume. This can be done by a tedious manipulation, see, for example, \cite{muller.i:thermodynamics}, which is essentially a diagonalisation procedure for the matrix of second derivatives via a suitable choice of variables. Another possibility is to formulate the concavity condition for a different thermodynamic potential.

  We want the entropy $\ctentropy$ to be a concave function of $\ctenergy$ and $\ctvolume$, and we want to reformulate this condition for the energy~$\ctenergy$ as a function of entropy $\ctentropy$ and volume $\ctvolume$. (What follows is classical material from convex analysis---convexity/concavity of inverse function to a convex/concave function and convexity/concavity of the Legendre transform.) We know that
  \begin{equation}
    \label{eq:178}
    \pd{\ctentropy}{\ctenergy}(\ctenergy, \ctvolume) = \frac{1}{\cttemperature} > 0,
  \end{equation}
  and the differentiation of $\ctentropy \left( \ctenergy(\ctentropy, \ctvolume), \ctvolume \right) = \ctentropy$ with respect to $\ctentropy$ yields
  \begin{equation}
    \label{eq:179}
    \pd{\ctentropy}{\ctenergy}(\ctenergy, \ctvolume) \pd{\ctenergy}{\ctentropy}(\ctentropy, \ctvolume) = 1,
  \end{equation}
  hence
  \begin{equation}
    \label{eq:180}
    \pd{\ctentropy}{\ctenergy}(\ctenergy, \ctvolume) = \frac{1}{\pd{\ctenergy}{\ctentropy}(\ctentropy, \ctvolume)}.
  \end{equation}
  The second derivative thus reads
  \begin{equation}
    \label{eq:181}
    \ppd{\ctentropy}{\ctenergy}(\ctenergy, \ctvolume) = -\frac{\ppd{\ctenergy}{\ctentropy}(\ctentropy, \ctvolume) \pd{\ctentropy}{\ctenergy}(\ctenergy, \ctvolume) }{\left( \pd{\ctenergy}{\ctentropy}(\ctentropy, \ctvolume) \right)^2},
  \end{equation}
  which due to the positivity of temperature~\eqref{eq:178} and the requirement~\eqref{eq:176} implies that we must have
  \begin{equation}
    \label{eq:182}
    \ppd{\ctenergy}{\ctentropy}(\ctentropy, \ctvolume) > 0.
  \end{equation}
  The differentiation of  $\ctentropy \left( \ctenergy(\ctentropy, \ctvolume), \ctvolume \right) = \ctentropy$ with respect to volume $\ctvolume$ then yields
  \begin{equation}
    \label{eq:183}
     \pd{\ctentropy}{\ctenergy}(\ctenergy, \ctvolume) \pd{\ctenergy}{\ctvolume}(\ctentropy, \ctvolume) +  \pd{\ctentropy}{\ctvolume}(\ctenergy, \ctvolume) = 0
   \end{equation}
   and consequently
   \begin{equation}
     \label{eq:184}
     \pd{\ctentropy}{\ctvolume}(\ctenergy, \ctvolume) = - \pd{\ctentropy}{\ctenergy}(\ctenergy, \ctvolume) \pd{\ctenergy}{\ctvolume}(\ctentropy, \ctvolume) = - \frac{\pd{\ctenergy}{\ctvolume}(\ctentropy, \ctvolume)}{\pd{\ctenergy}{\ctentropy}(\ctentropy, \ctvolume)}.
   \end{equation}
   Concerning the first derivatives we thus have
   \begin{subequations}
     \label{eq:185}
     \begin{align}
       \label{eq:186}
       \pd{\ctentropy}{\ctenergy}(\ctenergy, \ctvolume) &= \frac{1}{\pd{\ctenergy}{\ctentropy}(\ctentropy, \ctvolume)}, \\
       \pd{\ctentropy}{\ctvolume}(\ctenergy, \ctvolume) &= - \frac{\pd{\ctenergy}{\ctvolume}(\ctentropy, \ctvolume)}{\pd{\ctenergy}{\ctentropy}(\ctentropy, \ctvolume)}.
     \end{align}
   \end{subequations}
   Next we calculate the mixed derivative
   \begin{equation}
     \label{eq:187}
     \pd{^2\ctentropy}{\ctvolume \partial \ctenergy}(\ctenergy, \ctvolume)
     =
     \pd{}{\ctvolume} \left( \frac{1}{\pd{\ctenergy}{\ctentropy}(\ctentropy(\ctenergy, \ctvolume), \ctvolume)} \right)
     =
     -
     \frac{\ppd{\ctenergy}{\ctentropy}(\ctentropy, \ctvolume) \pd{\ctentropy}{\ctvolume}(\ctenergy, \ctvolume) + \pd{^2\ctenergy}{\ctvolume \partial \ctentropy}(\ctentropy, \ctvolume) }{\left( \pd{\ctenergy}{\ctentropy}(\ctentropy, \ctvolume) \right)^2}
   \end{equation}
   and the second partial derivative with respect to $\ctvolume$,
   \begin{multline}
     \label{eq:188}
     \ppd{\ctentropy}{\ctvolume}(\ctenergy, \ctvolume)
     =
     -
     \pd{}{\ctvolume} \left( \frac{\pd{\ctenergy}{\ctvolume}(\ctentropy(\ctenergy, \ctvolume), \ctvolume)}{\pd{\ctenergy}{\ctentropy}(\ctentropy(\ctenergy, \ctvolume), \ctvolume)} \right)
     =
     -
     \frac{
       \pd{}{\ctvolume} \left( \pd{\ctenergy}{\ctvolume}(\ctentropy(\ctenergy, \ctvolume), \ctvolume) \right) \pd{\ctenergy}{\ctentropy} (\ctentropy, \ctvolume)
       -
       \pd{\ctenergy}{\ctvolume}(\ctentropy, \ctvolume) \pd{}{\ctvolume} \left( \pd{\ctenergy}{\ctentropy}(\ctentropy(\ctenergy, \ctvolume), \ctvolume) \right)
   }
   {\left( \pd{\ctenergy}{\ctentropy}(\ctentropy, \ctvolume) \right)^2}
   \\
   =
   \frac{
     \left(
       \pd{^2\ctenergy}{\ctentropy \partial \ctvolume}(\ctentropy, \ctvolume)
       \pd{\ctentropy}{\ctvolume}(\ctenergy, \ctvolume)
       +
       \ppd{\ctenergy}{\ctvolume}(\ctentropy, \ctvolume)
     \right)
     \pd{\ctenergy}{\ctentropy} (\ctentropy, \ctvolume)
     -
     \pd{\ctenergy}{\ctvolume}(\ctentropy, \ctvolume)
     \left(
       \ppd{\ctenergy}{\ctentropy}(\ctentropy , \ctvolume)
       \pd{\ctentropy}{\ctvolume}(\ctenergy, \ctvolume)
       +
       \pd{^2\ctenergy}{\ctvolume \partial \ctentropy}(\ctentropy , \ctvolume) 
     \right)
   }
   {
     \left( \pd{\ctenergy}{\ctentropy}(\ctentropy, \ctvolume) \right)^2
   }
   \\
   =
   \frac{
     -
     2
     \pd{^2\ctenergy}{\ctvolume \partial \ctentropy}(\ctentropy , \ctvolume)
     \pd{\ctenergy}{\ctvolume}(\ctentropy, \ctvolume)
     +
     \ppd{\ctenergy}{\ctvolume}(\ctentropy, \ctvolume)
     \pd{\ctenergy}{\ctentropy} (\ctentropy, \ctvolume)
     -
     \ppd{\ctenergy}{\ctentropy}(\ctentropy , \ctvolume)
     \pd{\ctentropy}{\ctvolume}(\ctenergy, \ctvolume)
     \pd{\ctenergy}{\ctvolume}(\ctentropy, \ctvolume)
   }
   {
     \left( \pd{\ctenergy}{\ctentropy}(\ctentropy, \ctvolume) \right)^2
   }
 \end{multline}
 where we have made use of identities following from~\eqref{eq:185}.

 Now we are ready to reinterpret the stability condition~\eqref{eq:174} in terms of the derivatives of the energy. We get
 \begin{multline}
   \label{eq:189}
   \ppd{\ctentropy}{\ctenergy}(\ctenergy, \ctvolume)\ppd{\ctentropy}{\ctvolume}(\ctenergy, \ctvolume) - \left( \pd{^2\ctentropy}{\ctenergy \partial \ctvolume}(\ctenergy, \ctvolume) \right)^2
   \\
   =
   \frac{1}{\left( \pd{\ctenergy}{\ctentropy}(\ctentropy, \ctvolume) \right)^4}
   \left(
     \left[
       \ppd{\ctenergy}{\ctentropy}(\ctentropy, \ctvolume) \pd{\ctentropy}{\ctenergy}(\ctenergy, \ctvolume)
     \right]
     \left[
       -
       2
       \pd{^2\ctenergy}{\ctvolume \partial \ctentropy}(\ctentropy , \ctvolume)
       \pd{\ctenergy}{\ctvolume}(\ctentropy, \ctvolume)
       +
       \ppd{\ctenergy}{\ctvolume}(\ctentropy, \ctvolume)
       \pd{\ctenergy}{\ctentropy} (\ctentropy, \ctvolume)
       -
       \ppd{\ctenergy}{\ctentropy}(\ctentropy , \ctvolume)
       \pd{\ctentropy}{\ctvolume}(\ctenergy, \ctvolume)
       \pd{\ctenergy}{\ctvolume}(\ctentropy, \ctvolume)
     \right]
   \right.
   \\
     -
     \left.
     \left[
       \ppd{\ctenergy}{\ctentropy}(\ctentropy, \ctvolume) \pd{\ctentropy}{\ctvolume}(\ctenergy, \ctvolume) + \pd{^2\ctenergy}{\ctvolume \partial \ctentropy}(\ctentropy, \ctvolume)
     \right]^2
   \right)
   \\
   =
   \frac{1}{\left( \pd{\ctenergy}{\ctentropy}(\ctentropy, \ctvolume) \right)^4}
   \left(
     \left[
     2
     \ppd{\ctenergy}{\ctentropy}(\ctentropy, \ctvolume)
     \pd{^2\ctenergy}{\ctvolume \partial \ctentropy}(\ctentropy , \ctvolume)
     \pd{\ctentropy}{\ctvolume}(\ctenergy, \ctvolume)
     +
     \ppd{\ctenergy}{\ctentropy}(\ctentropy, \ctvolume)
     \ppd{\ctenergy}{\ctvolume}(\ctentropy, \ctvolume)
     +
     \left( \ppd{\ctenergy}{\ctentropy}(\ctentropy , \ctvolume) \right)^2
     \left( \pd{\ctentropy}{\ctvolume}(\ctenergy, \ctvolume) \right)^2
   \right]
 \right.
 \\
   -
   \left.
   \left[
     \ppd{\ctenergy}{\ctentropy}(\ctentropy, \ctvolume) \pd{\ctentropy}{\ctvolume}(\ctenergy, \ctvolume) + \pd{^2\ctenergy}{\ctvolume \partial \ctentropy}(\ctentropy, \ctvolume)
   \right]^2
 \right)
 \\
 =
 \frac{1}{\left( \pd{\ctenergy}{\ctentropy}(\ctentropy, \ctvolume) \right)^4}
 \left(
   \ppd{\ctenergy}{\ctentropy}(\ctentropy, \ctvolume)
   \ppd{\ctenergy}{\ctvolume}(\ctentropy, \ctvolume)
   -
   \left(
     \pd{^2\ctenergy}{\ctvolume \partial \ctentropy}(\ctentropy, \ctvolume)
   \right)^2
 \right)
 \\
 =
 \frac{1}{\cttemperature^4}
 \left(
   \ppd{\ctenergy}{\ctentropy}(\ctentropy, \ctvolume)
   \ppd{\ctenergy}{\ctvolume}(\ctentropy, \ctvolume)
   -
   \left(
     \pd{^2\ctenergy}{\ctvolume \partial \ctentropy}(\ctentropy, \ctvolume)
   \right)^2
 \right)
 .
\end{multline}
All together we have
\begin{subequations}
  \label{eq:190}
  \begin{align}
    \label{eq:191}
    \ppd{\ctentropy}{\ctenergy}(\ctenergy, \ctvolume) &= -\frac{1}{\cttemperature}\frac{\ppd{\ctenergy}{\ctentropy}(\ctentropy, \ctvolume)  }{\left( \pd{\ctenergy}{\ctentropy}(\ctentropy, \ctvolume) \right)^2}, \\
    \label{eq:192}
    \ppd{\ctentropy}{\ctenergy}(\ctenergy, \ctvolume)\ppd{\ctentropy}{\ctvolume}(\ctenergy, \ctvolume) - \left( \pd{^2\ctentropy}{\ctenergy \partial \ctvolume}(\ctenergy, \ctvolume) \right)^2
                                                      &=
                                                         \frac{1}{\cttemperature^4}
 \left(
   \ppd{\ctenergy}{\ctentropy}(\ctentropy, \ctvolume)
   \ppd{\ctenergy}{\ctvolume}(\ctentropy, \ctvolume)
   -
   \left(
     \pd{^2\ctenergy}{\ctvolume \partial \ctentropy}(\ctentropy, \ctvolume)
   \right)^2
 \right)
                                                        ,
  \end{align}
\end{subequations}
see~\eqref{eq:181}. The stability conditions, or the concavity conditions/negative definiteness of the second derivatives matrix for the entropy function thus translate as conditions
\begin{subequations}
  \label{eq:193}
  \begin{align}
    \label{eq:194}
    \ppd{\ctenergy}{\ctentropy}(\ctentropy, \ctvolume)  &> 0, \\
    \label{eq:195}
   \ppd{\ctenergy}{\ctentropy}(\ctentropy, \ctvolume)
   \ppd{\ctenergy}{\ctvolume}(\ctentropy, \ctvolume)
   -
   \left(
     \pd{^2\ctenergy}{\ctvolume \partial \ctentropy}(\ctentropy, \ctvolume)
   \right)^2
                                                        &>
                                                          0,
  \end{align}
\end{subequations}
which are the conditions for the convexity/positive definiteness of the second derivatives matrix for the energy function. The concavity of the entropy function $\ctentropy (\ctenergy, \ctvolume)$ with respect to its natural variables is thus equivalent to the convexity of the energy function $\ctenergy (\ctentropy, \ctvolume)$ with respect to its natural variables. Once we have translated the convexity/concavity conditions from the entropic equation of state to the energetic equation of state, we can easily proceed with other thermodynamic potentials derived from the energy, namely with the Helmholtz free energy, the Gibbs free energy and the enthalpy---it suffices to exploit the known convexity/concavity properties of the Legendre transform.

Now it is time to interpret the conditions~\eqref{eq:193}. First, we see that
\begin{equation}
  \label{eq:196}
  \ppd{\ctenergy}{\ctentropy}(\ctentropy, \ctvolume)
  =
  \pd{}{\ctentropy}
  \left(
    \pd{\ctenergy}{\ctentropy}(\ctentropy, \ctvolume)
  \right)
  =
  \pd{\cttemperature}{\ctentropy}(\ctentropy, \ctvolume)
  =
  \frac{1}{\pd{\ctentropy}{\cttemperature}(\cttemperature, \ctvolume)}
  =
  \frac{\cttemperature}{C_\ctvolume} > 0,
\end{equation}
which with the positivity of temperature implies that the heat capacity at constant volume must be positive. Furthermore, if we define the Helmholtz free energy as
\begin{equation}
  \label{eq:197}
  \cthelmholtz(\cttemperature, \ctvolume) =_{\bydefinition} \left. \left( \ctenergy(\ctentropy, \ctvolume) - \cttemperature \ctentropy \right) \right|_{\ctentropy = \ctentropy (\cttemperature, \ctvolume)},
\end{equation}
then we see that
\begin{equation}
  \label{eq:198}
  \pd{\cthelmholtz}{\ctvolume}(\cttemperature, \ctvolume) = \left. \pd{\ctenergy}{\ctvolume}(\ctentropy, \ctvolume) \right|_{\ctentropy = \ctentropy (\cttemperature, \ctvolume)}
  =
  - \ctpressure (\cttemperature, \ctvolume),
\end{equation}
hence
\begin{equation}
  \label{eq:199}
   \ppd{\cthelmholtz}{\ctvolume}(\cttemperature, \ctvolume) = - \pd{\ctpressure}{\ctvolume} (\cttemperature, \ctvolume).
 \end{equation}
 On the other hand, the direct differentiation yields, with the usual abuse of notation, the following formula
 \begin{multline}
   \label{eq:200}
   \ppd{\cthelmholtz}{\ctvolume}(\cttemperature, \ctvolume)
   =
   \pd{}{\ctvolume}
   \left(
     \pd{\ctenergy}{\ctvolume}(\ctentropy(\cttemperature, \ctvolume), \ctvolume)
   \right)
   =
   \pd{^2\ctenergy}{\ctentropy \partial \ctvolume}(\ctentropy, \ctvolume)
   \pd{\ctentropy}{\ctvolume}(\cttemperature, \ctvolume)
   +
   \ppd{\ctenergy}{\ctvolume}(\ctentropy, \ctvolume)
   \\
   =
   \frac{
     \pd{^2\ctenergy}{\ctentropy \partial \ctvolume}(\ctentropy, \ctvolume)
     \pd{\ctentropy}{\ctvolume}(\cttemperature, \ctvolume)
     \ppd{\ctenergy}{\ctentropy}(\ctentropy, \ctvolume)
     +
     \ppd{\ctenergy}{\ctvolume}(\ctentropy, \ctvolume)
     \ppd{\ctenergy}{\ctentropy}(\ctentropy, \ctvolume)
   }
   {
     \ppd{\ctenergy}{\ctentropy}(\ctentropy, \ctvolume)
   }
   \\
   =
   \frac{
     -
     \left(
       \pd{^2\ctenergy}{\ctentropy \partial \ctvolume}(\ctentropy, \ctvolume)
     \right)^2
     +
     \ppd{\ctenergy}{\ctvolume}(\ctentropy, \ctvolume)
     \ppd{\ctenergy}{\ctentropy}(\ctentropy, \ctvolume)
   }
   {
     \ppd{\ctenergy}{\ctentropy}(\ctentropy, \ctvolume)
   }
   ,
 \end{multline}
 where we have used the subtle identity
 \begin{equation}
   \label{eq:201}
   0
   =
   \pd{\cttemperature}{\ctvolume}
   =
   \pd{}{\ctvolume}
   \left(
     \pd{\ctenergy}{\ctentropy}(\ctentropy, \ctvolume)
   \right)
   =
   \ppd{\ctenergy}{\ctentropy}(\ctentropy, \ctvolume)
   \pd{\ctentropy}{\ctvolume} (\cttemperature, \ctvolume)
   +
   \pd{^2\ctenergy}{\ctvolume \partial \ctentropy}(\ctentropy, \ctvolume)
   .
 \end{equation}
 (The first partial derivative is meant as the derivative with respect to $\ctvolume$ at constant temperature $\cttemperature$. Furthermore, we have also used the definition of thermodynamic temperature.) In virtue of equalities~\eqref{eq:196}, \eqref{eq:199} and \eqref{eq:200} we thus see that the stability conditions~\eqref{eq:193} formulated for the energy reduce to requirements
 \begin{subequations}
   \label{eq:202}
   \begin{align}
     \label{eq:203}
     C_\ctvolume &>  0, \\
     \label{eq:204}
     - \pd{\ctpressure}{\ctvolume} (\cttemperature, \ctvolume) &>0,
   \end{align}
 \end{subequations}
 which can be phrased in plain language as requirements ``heat supplied to a box held at constant volume increases the temperature'' and ``pressure in a box kept at fixed temperature decreases with increasing volume of the box''. Conditions~\eqref{eq:202} are in fact the same conditions as~\eqref{eq:96}.

 \subsection{Stability conditions and convexity of a potential}
 \label{sec:conv-requ}
 The previous section showed us that in classical thermodynamics we have a convex thermodynamic potential. The same is true in the continuum mechanics setting, but we must carefully work with the density. It turns out that the following modification of the internal energy, 
\begin{equation}
  \label{eq:modified-internal-energy}
  \volienergy(\volentropy, \rho) =_{\bydefinition} \rho  \left. \ienergy (\entropy, \rho) \right|_{\entropy = \frac{\volentropy}{\rho}}
  ,
\end{equation}
is the right choice concerning the search for a convex function in the continuum mechanics setting. This strange combination is motivated by the effort to write the internal energy and the entropy as functions normalised to unit volume rather than to unit mass.

We investigate the convexity of the modified internal energy function~\eqref{eq:modified-internal-energy} via the behaviour of the matrix of second derivatives. The matrix of second derivatives reads
\begin{equation}
  \label{eq:206}
  \begin{bmatrix}
    \ppd{\volienergy}{\volentropy} & \pd{^2 \volienergy}{\volentropy \partial \rho} \\
    \pd{^2 \volienergy}{\volentropy \partial \rho} & \ppd{\volienergy}{\rho}
  \end{bmatrix}
  ,
\end{equation}
and for convexity of $\volienergy(\volentropy, \rho)$ we need the matrix to be positive definite. The positive definiteness is equivalent to the positivity of the leading principal minors, see Sylvester criterion, \cite{meyer.c:matrix}. We thus need to prove that
\begin{subequations}
  \label{eq:207}
  \begin{align}
    \label{eq:208}
    \ppd{\volienergy}{\volentropy} &> 0, \\
    \label{eq:209}
    \ppd{\volienergy}{\volentropy} \ppd{\volienergy}{\rho} - \left( \pd{^2 \volienergy}{\volentropy \partial \rho} \right)^2 &> 0,
  \end{align}
\end{subequations}
which, as we shall show in a moment, turns out to be equivalent to the stability conditions~\eqref{eq:128}.

\subsubsection{Proof of convexity of modified internal energy function}
\label{sec:proof-convexity}
Let us start with the calculations. Using the chain rule, we see that
   \begin{equation}
     \label{eq:210}
     \pd{\volienergy}{\volentropy}(\volentropy, \rho)
     =
     \pd{}{\volentropy}
     \left(
       \rho
       \ienergy
       \left(
         \frac{\volentropy}{\rho}
         ,
         \rho
       \right)
     \right)
     =
     \rho
     \left. \pd{\ienergy}{\entropy}(\entropy, \rho) \right|_{\entropy = \frac{\volentropy}{\rho}} \pd{\entropy}{\volentropy}(\volentropy, \rho)
     =
     \rho
     \left.
       \pd{\ienergy}{\entropy}(\entropy, \rho)
     \right|_{\entropy = \frac{\volentropy}{\rho}} 
     \frac{1}{\rho}
     =
     \left. \pd{\ienergy}{\entropy}(\entropy, \rho) \right|_{\entropy = \frac{\volentropy}{\rho}}.
   \end{equation}
   A similar calculation also reveals that
   \begin{equation}
     \label{eq:211}
     \ppd{\volienergy}{\volentropy}(\volentropy, \rho)
     =
     \frac{1}{\rho}
     \left.
       \ppd{\ienergy}{\entropy}(\entropy, \rho)
     \right|_{\entropy = \frac{\volentropy}{\rho}}
     .
   \end{equation}
   We already know that the specific heat capacity at constant volume $\cheatvol$ is given by the formula
 \begin{equation}
   \label{eq:212}
   \cheatvol (\temp, \rho)
   =
   \temp
   \pd{\entropy}{\temp}(\temp, \rho)
   ,
 \end{equation}
 see~\eqref{eq:77},  which further implies that
 \begin{equation}
   \label{eq:213}
   \cheatvol (\temp, \rho)
   =
   \temp
   \frac{
     1
   }
   {
     \pd{\temp}{\entropy}(\entropy, \rho)
   }
   =
   \temp
   \frac{
     1
   }
   {
     \pd{}{\entropy} \left(\pd{\ienergy}{\entropy}(\entropy, \rho) \right)
   }
   =
   \frac{
     \temp
   }
   {
     \ppd{\ienergy}{\entropy}(\entropy, \rho)
   }
   ,
 \end{equation}
 where we have used formula $\temp = \pd{\ienergy}{\entropy}(\entropy, \rho)$ for the temperature. Using~\eqref{eq:213} in~\eqref{eq:211} then yields
 \begin{equation}
   \label{eq:214}
   \ppd{\volienergy}{\volentropy}(\volentropy, \rho)
   =
   \left.
     \frac{\temp}{\rho \cheatvol(\temp, \rho)}
   \right|_{\temp = \temp (\volentropy, \rho)}. 
 \end{equation}
 Consequently, since the density and the temperature are positive, we see that~\eqref{eq:214} implies that the positivity of the first principal minor $ \ppd{\volienergy}{\volentropy}(\volentropy, \rho)$ is equivalent to the positivity of the specific heat capacity at constant volume~$\cheatvol$.

 Let us now focus on the pressure growth condition~\eqref{eq:130}. Employing the chain rule we see that
 \begin{multline}
   \label{eq:215}
   \pd{\ienergy}{\rho}(\entropy, \rho)
   =
   \pd{}{\rho}
   \left(
     \left.
       \frac{\volienergy(\volentropy, \rho)}{\rho}
     \right|_{\volentropy = \rho \entropy}
   \right)
   =
   -
   \frac{1}{\rho^2}
   \left.
     \volienergy(\volentropy, \rho)
   \right|_{\volentropy = \rho \entropy}
   +
   \frac{1}{\rho}
   \left.
     \pd{\volienergy}{\volentropy}(\volentropy, \rho)
   \right|_{\volentropy = \rho \entropy}
   \pd{\volentropy}{\rho}(\entropy, \rho)
   +
   \frac{1}{\rho}
   \left.
     \pd{\volienergy}{\rho}(\volentropy, \rho)
   \right|_{\volentropy = \rho \entropy}
   \\
   =
   -
   \frac{1}{\rho^2}
   \left.
     \volienergy(\volentropy, \rho)
   \right|_{\volentropy = \rho \entropy}
   +
   \frac{1}{\rho^2}
   \temp(\entropy, \rho)
   \volentropy(\entropy, \rho)
   +
   \frac{1}{\rho}
   \left.
     \pd{\volienergy}{\rho}(\volentropy, \rho)
   \right|_{\volentropy = \rho \entropy},
 \end{multline}
 where we have used the temperature formula~$\temp = \pd{\ienergy}{\entropy}(\entropy, \rho)$. Using the pressure formula~$
 \thpressure(\temp, \rho) = \rho^2
 \left.
   \pd{\ienergy}{\rho}(\entropy, \rho)
 \right|_{\entropy = \entropy(\temp, \rho)}
 $
 and the just derived identity for the derivative of energy with respect to the density, we see that the pressure growth condition can be rewritten as
 \begin{multline}
   \label{eq:216}
   \pd{\thpressure}{\rho}(\temp, \rho)
   =
   \pd{}{\rho}
   \left(
     \rho^2
     \left.
       \pd{\ienergy}{\rho}(\entropy, \rho)
     \right|_{\entropy = \entropy(\temp, \rho)}
   \right)
   =
   \pd{}{\rho}
   \left(
     -
     \left.
       \volienergy(\volentropy, \rho)
     \right|_{\volentropy = \rho \entropy}
     +
     \temp
     \volentropy(\temp, \rho)
     +
     \rho
     \left.
       \pd{\volienergy}{\rho}(\volentropy, \rho)
     \right|_{\volentropy = \rho \entropy}
   \right)
   \\
   =
   -
   \left.
     \pd{\volienergy}{\rho}(\volentropy, \rho)
   \right|_{\volentropy = \rho \entropy}
   -
   \left.
     \pd{\volienergy}{\volentropy}(\volentropy, \rho)
   \right|_{\volentropy = \rho \entropy}
   \pd{\volentropy}{\rho}(\temp, \rho)
   +
   \temp
   \pd{\volentropy}{\rho}(\temp, \rho)
   +
   \left.
     \pd{\volienergy}{\rho}(\volentropy, \rho)
   \right|_{\volentropy = \rho \entropy}
   +
   \rho
   \left.
     \ppd{\volienergy}{\rho}(\volentropy, \rho)
   \right|_{\volentropy = \rho \entropy}
   +
   \rho
   \left.
     \pd{^2\volienergy}{\rho \partial \volentropy}(\volentropy, \rho)
   \right|_{\volentropy = \rho \entropy}
   \pd{\volentropy}{\rho}(\temp, \rho)
   \\
   =
   \rho
   \left(
     \left.
       \ppd{\volienergy}{\rho}(\volentropy, \rho)
     \right|_{\volentropy = \rho \entropy}
     +
     \left.
       \pd{^2\volienergy}{\rho \partial \volentropy}(\volentropy, \rho)
     \right|_{\volentropy = \rho \entropy}
     \pd{\volentropy}{\rho}(\temp, \rho)
   \right)
   .
 \end{multline}
 Now we make a subtle observation based on the differentiation of the temperature with respect to the density at \emph{constant temperature}. This clearly yields zero, but we also see that
 \begin{equation}
   \label{eq:217}
   0
   =
   \pd{\temp}{\rho}(\temp, \rho)
   =
   \pd{}{\rho}
   \left(
     \left. \pd{\volienergy}{\volentropy}(\volentropy, \rho) \right|_{\volentropy = \volentropy(\temp, \rho)}
   \right)
   =
   \left. \pd{^2\volienergy}{\volentropy \partial \rho}(\volentropy, \rho) \right|_{\volentropy = \volentropy(\temp, \rho)}
   +
   \left. \ppd{\volienergy}{\volentropy}(\volentropy, \rho) \right|_{\volentropy = \volentropy(\temp, \rho)}
   \pd{\volentropy}{\rho}(\temp, \rho)
   ,
 \end{equation}
 where we have used the temperature formula~$\temp = \pd{\ienergy}{\entropy}(\entropy, \rho)$ and identity~\eqref{eq:210}. Equality~\eqref{eq:217} yields
 \begin{equation}
   \label{eq:218}
   \pd{\volentropy}{\rho}(\temp, \rho)
   =
   -
   \frac{
     \left. \pd{^2\volienergy}{\volentropy \partial \rho}(\volentropy, \rho) \right|_{\volentropy = \volentropy(\temp, \rho)}
   }
   {
     \left. \ppd{\volienergy}{\volentropy}(\volentropy, \rho) \right|_{\volentropy = \volentropy(\temp, \rho)}
   }
   ,
 \end{equation}
 which upon substitution into~\eqref{eq:216} reveals that
 \begin{multline}
   \label{eq:219}
   \pd{\thpressure}{\rho}(\temp, \rho)
   =
   \frac{
     \rho
     \left(
       \left.
         \ppd{\volienergy}{\rho}(\volentropy, \rho)
       \right|_{\volentropy = \rho \entropy}
       \left.
         \ppd{\volienergy}{\volentropy}(\volentropy, \rho)
       \right|_{\volentropy = \volentropy(\temp, \rho)}
       -
       \left(
         \left.
           \pd{^2\volienergy}{\rho \partial \volentropy}(\volentropy, \rho)
         \right|_{\volentropy = \rho \entropy}
       \right)^2
     \right)
   }
   {
     \left.
       \ppd{\volienergy}{\volentropy}(\volentropy, \rho)
     \right|_{\volentropy = \volentropy(\temp, \rho)}
   }
   \\
   =
   \frac{\rho^2 \cheatvol}{\temp}
   \left.
     \left(
       \ppd{\volienergy}{\rho}(\volentropy, \rho)
       \ppd{\volienergy}{\volentropy}(\volentropy, \rho)
       -
       \left(
         \pd{^2\volienergy}{\rho \partial \volentropy}(\volentropy, \rho)
       \right)^2
     \right)
   \right|_{\volentropy = \volentropy(\temp, \rho)}
   ,
 \end{multline}
 where the expression on the right-hand side is the same as the expression in the condition on the positivity of the second principal minor, see~\eqref{eq:209}, of the second derivatives matrix. We thus see that the pressure growth condition~\eqref{eq:130} is tantamount to the positivity of the second principal minor. Overall we thus have
 \begin{subequations}
   \label{eq:220}
   \begin{align}
     \label{eq:221}
     \ppd{\volienergy}{\volentropy}(\volentropy, \rho)
     =
       \left.
       \frac{\temp}{\rho \cheatvol(\temp, \rho)}
       \right|_{\temp = \temp (\volentropy, \rho)}
       &> 0
       , \\
     \label{eq:222}
     \pd{\thpressure}{\rho}(\temp, \rho)
     =
     \frac{\rho^2 \cheatvol}{\temp}
     \left.
     \left(
     \ppd{\volienergy}{\rho}(\volentropy, \rho)
     \ppd{\volienergy}{\volentropy}(\volentropy, \rho)
     -
     \left(
     \pd{^2\volienergy}{\rho \partial \volentropy}(\volentropy, \rho)
     \right)^2
     \right)
     \right|_{\volentropy = \volentropy(\temp, \rho)}
     &>
       0
       ,
   \end{align}
 \end{subequations}
 which concludes the proof of the relation between the convexity of $\volienergy(\volentropy, \rho)$ and the stability conditions~\eqref{eq:128}. Our findings are summarised in Summary~\ref{summary:convexity}.

 \begin{summary}[Convexity of modified internal energy function $\volienergy(\volentropy, \rho) =_{\bydefinition} \rho  \left. \ienergy (\entropy, \rho) \right|_{\entropy = \frac{\volentropy}{\rho}}$]
   \label{summary:convexity}
   Let us consider a fluid with the Helmholtz free energy $\fenergy$ that leads to thermodynamic stability conditions
   \begin{subequations}
     \label{eq:131}
     \begin{align}
       \label{eq:132}
       \cheatvol \left(\temp, \rho\right) &> 0, \\
       \label{eq:133}
       \pd{\thpressure}{\rho}(\temp, \rho) &> 0,
     \end{align}
   \end{subequations}
   and let us define the modified internal energy as
   \begin{equation}
     \label{eq:134}
     \volienergy(\volentropy, \rho) =_{\bydefinition} \rho  \left. \ienergy (\entropy, \rho) \right|_{\entropy = \frac{\volentropy}{\rho}},
   \end{equation}
   where $\ienergy (\entropy, \rho)$ is the internal energy that corresponds to the Helmholtz free energy $\fenergy$. Then the matrix of second derivatives (Hessian matrix) of the modified internal energy function reads 
   \begin{equation}
     \label{eq:135}
     \hessian{\volienergy}(\volentropy, \rho)
     =
     \begin{bmatrix}
       \ppd{\volienergy}{\volentropy} & \pd{^2 \volienergy}{\volentropy \partial \rho} \\
       \pd{^2 \volienergy}{\volentropy \partial \rho} & \ppd{\volienergy}{\rho}
     \end{bmatrix}
     ,
   \end{equation}
   and the matrix is positive definite as its principal minors are given by the formulae
   \begin{subequations}
     \label{eq:137}
     \begin{align}
       \label{eq:280}
       \ppd{\volienergy}{\volentropy}(\volentropy, \rho)
       &=
       \left.
       \frac{\temp}{\rho \cheatvol(\temp, \rho)}
       \right|_{\temp = \temp (\volentropy, \rho)}
         , \\
       \label{eq:284}
       \left(
       \ppd{\volienergy}{\rho}(\volentropy, \rho)
       \ppd{\volienergy}{\volentropy}(\volentropy, \rho)
       -
       \left(
       \pd{^2\volienergy}{\rho \partial \volentropy}(\volentropy, \rho)
       \right)^2
       \right)
       &=
       \left.
       \frac{\temp}{\rho^2 \cheatvol (\temp, \rho)}
       \pd{\thpressure}{\rho}(\temp, \rho)
       \right|_{\temp = \temp (\volentropy, \rho)}
         ,
     \end{align}
   \end{subequations}
   and hence positive in the virtue of conditions~\eqref{eq:131}. The modified internal energy $\volienergy(\volentropy, \rho)$ is thus a convex function of its natural variables $\volentropy$ and $\rho$.
 \end{summary}

 \begin{remark}[Convexity of the modified internal energy and the perspective of a convex function]
   \label{rm:perspective}
   The convexity of the modified internal energy $\volienergy(\volentropy, \rho)$ can be also seen without the calculation of the principal minors. Let $\upsilon =_{\bydefinition} \frac{1}{\rho}$ denote the specific volume, and let
   \begin{equation}
     \label{eq:381}
     \check{\ienergy}(\entropy, \upsilon)
     =_{\bydefinition}
     \ienergy \left( \entropy, \frac{1}{\upsilon} \right)
   \end{equation}
   denote the specific internal energy as a function of the specific entropy $\entropy$ and the specific volume $\upsilon$. These are the natural variables of classical thermodynamics, see Section~\ref{sec:revis-constr-from}, just normalised to unit mass. The analysis in Section~\ref{sec:revis-constr-from}, wherein we formally replace $\ctenergy$, $\ctentropy$ and $\ctvolume$ by $\check{\ienergy}$, $\entropy$ and $\upsilon$, shows that $\check{\ienergy}(\entropy, \upsilon)$ is a strictly convex function provided that $\cheatvol > 0$ and $- \pd{\thpressure}{\upsilon} > 0$, see~\eqref{eq:193} and~\eqref{eq:202}. Since $\pd{\thpressure}{\upsilon} = - \rho^2 \pd{\thpressure}{\rho}$, these are again the thermodynamic stability conditions~\eqref{eq:131}.

   The definition of the modified internal energy, see~\eqref{eq:modified-internal-energy}, can be rewritten in terms of $\check{\ienergy}$ as
   \begin{equation}
     \label{eq:382}
     \volienergy(\volentropy, \rho)
     =
     \rho
     \check{\ienergy} \left( \frac{\volentropy}{\rho}, \frac{1}{\rho} \right)
     .
   \end{equation}
   This is a special instance of the \emph{perspective} of a function, see, for example, \cite[Section~3.2.6]{boyd.s.vandenberghe.l:convex}.

   The \emph{perspective of a function} is the following construction. If $g: \R^n \mapsto \R$ is a function, then its perspective is the function of $n + 1$ variables defined as
   \begin{equation}
     \label{eq:383}
     P_g(t, \vec{x})
     =_{\bydefinition}
     t g \left( \frac{\vec{x}}{t} \right),
     \qquad
     t > 0
     .
   \end{equation}
   If $g$ is convex, then its perspective $P_g$ is convex as a function of the pair $\left[t, \vec{x}\right]$. Indeed, let
   \begin{equation}
     \label{eq:386}
     \left[t, \vec{x}\right] = \lambda \left[t_1, \vec{x}_1\right] + (1 - \lambda) \left[t_2, \vec{x}_2\right]
   \end{equation}
   with $\lambda \in (0, 1)$, and let us denote $\mu =_{\bydefinition} \frac{\lambda t_1}{t} \in (0, 1)$. Then $1 - \mu = \frac{(1 - \lambda) t_2}{t}$, and
   \begin{equation}
     \label{eq:384}
     \frac{\vec{x}}{t}
     =
     \mu
     \frac{\vec{x}_1}{t_1}
     +
     (1 - \mu)
     \frac{\vec{x}_2}{t_2}
   \end{equation}
   is a convex combination of the points $\frac{\vec{x}_1}{t_1}$ and $\frac{\vec{x}_2}{t_2}$. The convexity of $g$ then yields
   \begin{equation}
     \label{eq:385}
     P_g(t, \vec{x})
     =
     t  g \left( \frac{\vec{x}}{t} \right)
     \leq
     t
     \left[
       \mu  g \left( \frac{\vec{x}_1}{t_1} \right)
       +
       (1 - \mu)  g \left( \frac{\vec{x}_2}{t_2} \right)
     \right]
     =
     \lambda P_g(t_1, \vec{x}_1)
     +
     (1 - \lambda) P_g(t_2, \vec{x}_2)
     ,
   \end{equation}
   which is the definition of convexity of $P_g$.
   
   The specific internal energy $\check{\ienergy}$ is a function of two variables, hence its perspective is a function of three variables $t$ and $\vec{x} = [\xi_1, \xi_2]$,
   \begin{equation}
     \label{eq:389}
     P_{\check{\ienergy}}(t, [\xi_1, \xi_2])
     =
     t \, \check{\ienergy} \left( \frac{\xi_1}{t}, \frac{\xi_2}{t} \right)
     .
   \end{equation}
   Formula~\eqref{eq:382} shows that
   \begin{equation}
     \label{eq:388}
     \volienergy(\volentropy, \rho)
     =
     \left.
       P_{\check{\ienergy}}(t, [\xi_1, \xi_2])
     \right|_{\left[t, \xi_1, \xi_2\right] = \left[\rho, \volentropy, 1\right]}
     ,
   \end{equation}
   that is the modified internal energy is the perspective of the specific internal energy $\check{\ienergy}$ composed with the affine map $\left[\volentropy, \rho\right] \mapsto \left[t, \xi_1, \xi_2\right] = \left[\rho, \volentropy, 1\right]$. The density plays a double role here. It is the prefactor $t$, and it is also the quantity by which the components $\xi_1 = \volentropy$ and $\xi_2 = 1$ are divided, which gives $\left[ \frac{\xi_1}{t}, \frac{\xi_2}{t} \right] = \left[ \frac{\volentropy}{\rho}, \frac{1}{\rho} \right] = \left[\entropy, \upsilon\right]$. In particular, the constant $1$ in the second component turns into the specific volume $\upsilon$ upon the division by $t = \rho$. Since the composition of a convex function with an affine map is a convex function, $\volienergy(\volentropy, \rho)$ is convex. Moreover, it is \emph{strictly} convex. If the equality holds in~\eqref{eq:385}, then the strict convexity of $\check{\ienergy}$ implies $\frac{\vec{x}_1}{t_1} = \frac{\vec{x}_2}{t_2}$. In our case the second component of $\vec{x}$ is fixed to one, $\xi_2 = 1$, hence $\frac{1}{\rho_1} = \frac{1}{\rho_2}$, which further implies $\rho_1 = \rho_2$ and $\volentropy_1 = \volentropy_2$.

   The perspective is in fact the mathematical counterpart of the passage from the quantities normalised to unit mass to the quantities normalised to unit volume---we divide the variables by the density, and we multiply the function by the density. The same construction is behind the kinetic energy per unit volume
   \begin{equation}
     \label{eq:387}
     \frac{\absnorm{\linmomentum}^2}{2 \rho} = \rho  g\left(\frac{\linmomentum}{\rho}\right),
   \end{equation}
   which is the perspective of the convex function $g(\vec{v}) = \frac{1}{2} \absnorm{\vec{v}}^2$, see Section~\ref{sec:kinet-energy-compl}.
 \end{remark}

 \subsubsection{Convexity and Bregman distance/divergence}
 \label{sec:conv-bregm-dist}
Convexity of $\volienergy(\volentropy, \rho)$ can be exploited in the construction of the Bregman distance/divergence, see~\cite{bregman.lm:relaxation}. Bregman distance/divergence is a construction that allows one to introduce a concept of distance/divergence generated by a convex function.

\begin{definition}[Bregman divergence]
  \label{dfn:bregman-divergence}
   Let $f: \R^n \mapsto \R$ be a differentiable and strictly convex function. Then
   \begin{equation}
     \label{eq:285}
     \bregmandivergenceop{f}{\vec{x}}{\vec{y}} =_{\bydefinition} f(\vec{x}) - f(\vec{y}) - \vectordot{\left( \left. \pd{f}{\vec{z}} \right|_{\vec{z} = \vec{y}} \right)}{\left( \vec{x} - \vec{y} \right)}
   \end{equation}
   is referred to as the Bregman distance/divergence induced by the function $f$.
 \end{definition}
 It turns out that the Bregman distance/divergence has convenient properties, namely
 \begin{equation}
   \label{eq:286}
   \forall \vec{x}, \vec{y} \in \R^n: \bregmandivergenceop{f}{\vec{x}}{\vec{y}} \geq 0,
 \end{equation}
 while
 \begin{equation}
   \label{eq:226}
   \bregmandivergenceop{f}{\vec{x}}{\vec{y}} = 0
 \end{equation}
 if and only if $\vec{x} = \vec{y}$. These properties are simple consequences of convexity of $f$. Despite its name, the Bregman distance/divergence does not satisfy the properties of a true distance function. In particular, it is not necessarily true that $\bregmandivergenceop{f}{\vec{x}}{\vec{y}} = \bregmandivergenceop{f}{\vec{y}}{\vec{x}}$, which would hold if we were dealing with the genuine concept of distance. (It is therefore preferred to use the terminology Bregman \emph{divergence}.) For further discussion of Bregman divergence see~\cite{amari.s:information} and references therein. 

 We note that Bregman distance/divergence is strictly convex in the first variable. Consequently, if we construct the Bregman divergence out of a convex function $f(\vec{z})$,
 \begin{equation}
   \label{eq:6}
   \bregmandivergenceop{f}{\vec{x}}{\vec{y}} =_{\bydefinition} f(\vec{x}) - f(\vec{y}) - \vectordot{\left( \left. \pd{f}{\vec{z}} \right|_{\vec{z} = \vec{y}} \right)}{\left( \vec{x} - \vec{y} \right)},
 \end{equation}
 then we get, for fixed $\vec{y}$, yet another convex function
 \begin{equation}
   \label{eq:299}
   g(\vec{z}) =_{\bydefinition} \bregmandivergenceop{f}{\vec{z}}{\vec{y}}.
 \end{equation}
 Now we can apply the same Bregman distance/divergence construction to the new convex function $g$,
 \begin{equation}
   \label{eq:290}
   \bregmandivergenceop{g}{\vec{u}}{\vec{v}}
   =_{\bydefinition}
   \bregmandivergenceop{f}{\vec{u}}{\vec{y}}
   -
   \bregmandivergenceop{f}{\vec{v}}{\vec{y}}
   -
   \vectordot{
     \left(
       \left.
         \pd{}{\vec{x}}
         \bregmandivergenceop{f}{\vec{x}}{\vec{y}}
       \right|_{\vec{x} = \vec{v}}
     \right)
   }
   {
     \left(
       \vec{u} - \vec{v}
     \right)
   }
   .
 \end{equation}
 We use the definition~\eqref{eq:6} of $\bregmandivergenceop{f}{\cdot}{\vec{y}}$, and we work out the formula~\eqref{eq:290}
 \begin{multline}
   \label{eq:300}
   \bregmandivergenceop{g}{\vec{u}}{\vec{v}}
   =
   \left[
     f(\vec{u}) - f(\vec{y}) - \vectordot{\left( \left. \pd{f}{\vec{z}} \right|_{\vec{z} = \vec{y}} \right)}{\left( \vec{u} - \vec{y} \right)}
   \right]
   -
   \left[
     f(\vec{v}) - f(\vec{y}) - \vectordot{\left( \left. \pd{f}{\vec{z}} \right|_{\vec{z} = \vec{y}} \right)}{\left( \vec{v} - \vec{y} \right)}
   \right]
   \\
   -
   \vectordot{
     \left[
       \left.
         \pd{}{\vec{x}}
         \left(
           f(\vec{x})
           -
           f(\vec{y})
           -
           \vectordot{\left( \left. \pd{f}{\vec{z}} \right|_{\vec{z} = \vec{y}} \right)}{\left( \vec{x} - \vec{y} \right)}
         \right)
       \right|_{\vec{x} = \vec{v}}
     \right]
   }
   {
     \left(
       \vec{u} - \vec{v}
     \right)
   }
   \\
   =
   f(\vec{u})
   -
   f(\vec{v})
   -
   \vectordot{\left( \left. \pd{f}{\vec{z}} \right|_{\vec{z} = \vec{y}} \right)}{\left( \vec{u} - \vec{y} \right)}
   +
   \vectordot{\left( \left. \pd{f}{\vec{z}} \right|_{\vec{z} = \vec{y}} \right)}{\left( \vec{v} - \vec{y} \right)}
   -
   \vectordot{\left( \left. \pd{f}{\vec{x}} \right|_{\vec{x} = \vec{v}} \right)}{\left( \vec{u} - \vec{v}\right)}
   +
   \vectordot{\left( \left. \pd{f}{\vec{z}} \right|_{\vec{z} = \vec{y}} \right)}{\left( \vec{u} - \vec{v} \right)}
   \\
   =
   f(\vec{u})
   -
   f(\vec{v})
   -
   \vectordot{\left( \left. \pd{f}{\vec{x}} \right|_{\vec{x} = \vec{v}} \right)}{\left( \vec{u} - \vec{v}\right)}
   ,
 \end{multline}
 and we see that
 \begin{equation}
   \label{eq:301}
   \bregmandivergenceop{g}{\vec{u}}{\vec{v}}
   =
   \bregmandivergenceop{f}{\vec{u}}{\vec{v}}.
 \end{equation}
 We can summarise our findings as a lemma.
 
  \begin{lemma}[Bregman divergence constructed out of Bregman divergence with a fixed reference point]
   \label{lm:1}
   Let $f: \R^n \mapsto \R$ be a differentiable and strictly convex function, and let $\bregmandivergenceop{f}{\vec{x}}{\vec{y}}$ be the corresponding Bregman divergence,
   \begin{equation}
     \label{eq:223}
     \bregmandivergenceop{f}{\vec{x}}{\vec{y}} =_{\bydefinition} f(\vec{x}) - f(\vec{y}) - \vectordot{\left( \left. \pd{f}{\vec{z}} \right|_{\vec{z} = \vec{y}} \right)}{\left( \vec{x} - \vec{y} \right)}.
   \end{equation}
   Let us fix $\vec{y}$ in $\R^n$ and let us denote $g(\vec{z}) = _{\bydefinition} \bregmandivergenceop{f}{\vec{z}}{\vec{y}}$. Function $g(\vec{z})$ is a strictly convex function, and we can construct the corresponding Bregman divergence out of function $g$, that is 
   \begin{equation}
     \label{eq:224}
     \bregmandivergenceop{g}{\vec{u}}{\vec{v}} =_{\bydefinition} g(\vec{u}) - g(\vec{v}) - \vectordot{\left( \left. \pd{g}{\vec{z}} \right|_{\vec{z} = \vec{v}} \right)}{\left( \vec{u} - \vec{v} \right)}.
   \end{equation}
   It holds
   \begin{equation}
     \label{eq:225}
     \bregmandivergenceop{g}{\vec{u}}{\vec{v}} = \bregmandivergenceop{f}{\vec{u}}{\vec{v}}.
   \end{equation}
 \end{lemma}

 The lemma gives us a hint how to introduce a divergence between two states $\vec{u}$ and $\vec{v}$ provided that we have at our disposal a divergence that measures the proximity between the given states $\vec{u}$ and $\vec{v}$ and a given reference state $\widehat{\vec{w}}_{\mathrm{eq}}$. If we start with Bregman divergence $\bregmandivergenceop{f}{\vec{w}}{\widehat{\vec{w}}_{\mathrm{eq}}}$ characterising the proximity of a generic state~$\vec{w}$ and a \emph{reference} state~$\widehat{\vec{w}}_{\mathrm{eq}}$, then the \emph{new} Bregman divergence $\bregmandivergenceop{g}{\cdot}{\cdot}$ constructed with the \emph{intention to measure the proximity between arbitrary states}~$\vec{u}$ and~$\vec{v}$ might be defined as the Bregman divergence associated to the strictly convex function $\bregmandivergenceop{f}{\cdot}{\widehat{\vec{w}}_{\mathrm{eq}}}$. According to Lemma~\ref{lm:1} we then have
 \begin{equation}
   \label{eq:227}
   \bregmandivergenceop{g}{\vec{u}}{\vec{v}} = \bregmandivergenceop{f}{\vec{u}}{\vec{v}},
 \end{equation}
 which means that the reference state $\widehat{\vec{w}}_{\mathrm{eq}}$ is no longer necessary for measuring the proximity between the given states, and that the new Bregman divergence has the same functional form as the original Bregman divergence; in other words we have just swapped $\widehat{\vec{w}}_{\mathrm{eq}}$ for $\vec{v}$.

  We note that the construction of Bregman divergence out of a Bregman divergence is tantamount to the \emph{affine correction trick} for a convex function, see~\cite{bulcek.m.malek.j.ea:thermodynamics} and references therein.

  \subsubsection{Quadratic part of Bregman distance/divergence}
  \label{sec:quadr-part-bregm}
  Let $f: \R^n \mapsto \R$ be a twice continuously differentiable strictly convex function, and let $\hessian{f}(\vec{z})$ denote its matrix of second derivatives (Hessian matrix),
  \begin{equation}
    \label{eq:378}
  \left[ \hessian{f}(\vec{z}) \right]_{ij}
  =_{\bydefinition}
  \pd{^2 f}{z_i \partial z_j}(\vec{z})
  .
\end{equation}
If we apply the Taylor formula with the remainder in the integral form, $\phi(1) = \phi(0) + \left. \dd{\phi}{\tau} \right|_{\tau = 0} + \int_0^1 (1 - \tau) \ddd{\phi}{\tau} \, \diff \tau$, to the function $\phi(\tau) =_{\bydefinition} f(\vec{y} + \tau (\vec{x} - \vec{y}))$, we get
\begin{equation}
  \label{eq:379}
  f(\vec{x})
  =
  f(\vec{y})
  +
  \vectordot{\left( \left. \pd{f}{\vec{z}} \right|_{\vec{z} = \vec{y}} \right)}{\left( \vec{x} - \vec{y} \right)}
  +
  \int_0^1
  (1 - \tau)
  \vectordot{\left( \vec{x} - \vec{y} \right)}{\hessian{f}\left(\vec{y} + \tau (\vec{x} - \vec{y})\right) \left( \vec{x} - \vec{y} \right)}
  \,
  \diff \tau
  .
\end{equation}
Comparing~\eqref{eq:379} with the definition of the Bregman distance/divergence, see Definition~\ref{dfn:bregman-divergence}, we see that the Bregman distance/divergence is \emph{exactly} the remainder in the Taylor formula,
\begin{equation}
  \label{eq:380}
  \bregmandivergenceop{f}{\vec{x}}{\vec{y}}
  =
  \int_0^1
  (1 - \tau)
  \vectordot{\left( \vec{x} - \vec{y} \right)}{\hessian{f}\left(\vec{y} + \tau (\vec{x} - \vec{y})\right) \left( \vec{x} - \vec{y} \right)}
  \,
  \diff \tau
  .
\end{equation}
(Formula~\eqref{eq:380} also provides yet another proof of the nonnegativity property~\eqref{eq:286}. If the Hessian matrix is positive definite everywhere, then the integrand is positive for $\vec{x} \neq \vec{y}$.) Since the Hessian matrix is continuous, we can replace $\hessian{f}\left(\vec{y} + \tau (\vec{x} - \vec{y})\right)$ by $\hessian{f}(\vec{y})$ at the expense of a term of order higher than two in $\vec{x} - \vec{y}$, and the integration with respect to $\tau$ yields the factor $\int_0^1 (1 - \tau) \, \diff \tau = \frac{1}{2}$. We can summarise this observation as a lemma.

\begin{lemma}[Quadratic part of Bregman distance/divergence]
  \label{lm:3}
  Let $f: \R^n \mapsto \R$ be a twice continuously differentiable strictly convex function, and let $\bregmandivergenceop{f}{\vec{x}}{\vec{y}}$ be the corresponding Bregman distance/divergence. Then
  \begin{equation}
    \label{eq:lin-bregman-quadratic-part}
    \bregmandivergenceop{f}{\vec{x}}{\vec{y}}
    =
    \frac{1}{2}
    \vectordot{\left( \vec{x} - \vec{y} \right)}{\hessian{f}(\vec{y}) \left( \vec{x} - \vec{y} \right)}
    +
    \smallo{\absnorm{\vec{x} - \vec{y}}^2}
  \end{equation}
  as $\vec{x} \to \vec{y}$, where $\hessian{f}(\vec{y})$ denotes the Hessian matrix of $f$ evaluated at the point $\vec{y}$. In other words, the leading order term in the Bregman distance/divergence $\bregmandivergenceop{f}{\vec{x}}{\vec{y}}$ is the quadratic form generated by the Hessian matrix of $f$ at the reference point $\vec{y}$.
\end{lemma}

 \subsubsection{Bregman distance/divergence generated by the modified internal energy function}
\label{sec:bregm-dist-gener}

 We now find the Bregman distance/divergence induced by the convex function $\volienergy(\volentropy, \rho)$. Using the definition of Bregman distance/divergence we get
 \begin{equation}
   \label{eq:287}
   \bregmandivergenceop{\volienergy}{\vec{A}}{\widehat{\vec{A}}}
   =_{\bydefinition}
   \volienergy \left( \vec{A} \right)
   -
   \volienergy \left( \widehat{\vec{A}} \right)
   -
   \vectordot{\left( \left. \pd{\volienergy(\vec{Z})}{\vec{Z}} \right|_{\vec{Z} = \widehat{\vec{A}}} \right)}{\left( \vec{A} - \widehat{\vec{A}} \right)},
 \end{equation}
 where we denote $\vec{A} = \left[ \volentropy, \rho \right]$ and $\widehat{\vec{A}} = \left[ \widehat{\volentropy}, \widehat{\rho} \right]$. Substituting into~\eqref{eq:287} we get
 \begin{equation}
   \label{eq:288}
   \bregmandivergenceop{\volienergy}{\vec{A}}{\widehat{\vec{A}}}
   =
   \volienergy (\volentropy, \rho)
   -
   \volienergy (\widehat{\volentropy}, \widehat{\rho})
   -
   \left.
     \pd{\volienergy}{\volentropy}(\volentropy, \rho)
   \right|_{\left[\volentropy, \rho\right] = \left[\widehat{\volentropy}, \widehat{\rho}\right]}
   \left(
     \volentropy
     -
     \widehat{\volentropy}
   \right)
   -
   \left.
     \pd{\volienergy}{\rho}(\volentropy, \rho)
   \right|_{\left[\volentropy, \rho\right] = \left[\widehat{\volentropy}, \widehat{\rho}\right]}
   \left(
     \rho
     -
     \widehat{\rho}
   \right)
   ,
 \end{equation}
 which we try to rewrite using the unmodified internal energy $\ienergy(\entropy, \rho)$ and the density and temperature. We already know that
 \begin{subequations}
   \label{eq:292}
   \begin{align}
     \label{eq:293}
     \pd{\volienergy}{\volentropy}(\volentropy, \rho)
     &=
       \left.
       \left[
       \pd{\ienergy}{\entropy}(\entropy, \rho)
       \right]
       \right|_{\entropy = \frac{\volentropy}{\rho}}, \\
     \label{eq:294}
     \pd{\volienergy}{\rho}(\volentropy, \rho)
     &=
       \left.
       \left[
       \ienergy(\entropy, \rho)
       -
       \temp(\entropy, \rho)
       \entropy
       +
       \rho \pd{\ienergy}{\rho}(\entropy, \rho)
       \right]
       \right|_{\entropy = \frac{\volentropy}{\rho}},
   \end{align}
 \end{subequations}
 see~\eqref{eq:210} and~\eqref{eq:215}. Using the definition of thermodynamic temperature and the thermodynamic pressure and the definition of the modified entropy $\volentropy = \rho \entropy$, we thus see that~\eqref{eq:288} reduces to
 \begin{equation}
   \label{eq:291}
   \bregmandivergenceop{\volienergy}{\vec{A}}{\widehat{\vec{A}}}
   =
   \rho
   \ienergy
   -
   \widehat{\rho}
   \widehat{\ienergy}
   -
   \widehat{\temp}
   \left(
     \rho
     \entropy
     -
     \widehat{\rho}
     \widehat{\entropy}
   \right)
   -
   \left(
     \widehat{\ienergy}
     -
     \widehat{\temp}
     \widehat{\entropy}
     +
     \frac{\widehat{\thpressure}}{\widehat{\rho}}
   \right)
   \left(
     \rho
     -
     \widehat{\rho}
   \right)
   .
 \end{equation}
 This formula can be further rewritten as
 \begin{equation}
   \label{eq:295}
   \bregmandivergenceop{\volienergy}{\vec{A}}{\widehat{\vec{A}}}
   =
   -
   \widehat{\temp}
   \left\{
     \rho
     \entropy
     -
     \widehat{\rho} \widehat{\entropy}
   \right\}
   +
   \left\{
     \rho \ienergy
     -
     \widehat{\rho}
     \widehat{\ienergy}
   \right\}
   -
   \left\{
     \frac{
       \widehat{\thpressure}
     }
     {
       \widehat{\rho}
     }
     +
     \left(
       \widehat{\ienergy}
       -
       \widehat{\temp} \widehat{\entropy}
     \right)
   \right\}
   \left(
     \rho - \widehat{\rho}
   \right)
   ,
 \end{equation}
 which is the same expression we have seen as the integrand in the Lyapunov type functional $\mathcal{V}_{\mathrm{meq}, \, \widehat{\temp}, \, \widehat{\rho}}(\rho, \vec{v}, \temp)$, see~\eqref{eq:68}, that is
 \begin{equation}
   \label{eq:22}
   \mathcal{V}_{\mathrm{meq}, \, \widehat{\temp}, \, \widehat{\rho}}(\rho, \vec{v}, \temp)
   =
   \int_{\Omega} \frac{1}{2} \rho \absnorm{\vec{v}}^2 \, \cvolumee
   +
   \int_{\Omega}
   \left[
     -
     \widehat{\temp} \left\{ \rho \entropy(\temp, \rho) - \widehat{\rho} \entropy(\widehat{\temp}, \widehat{\rho}) \right\}
     +
     \left\{ \rho \ienergy(\temp, \rho) - \widehat{\rho} \ienergy(\widehat{\temp}, \widehat{\rho}) \right\}
     -
     \left( \frac{\thpressure(\widehat{\temp}, \widehat{\rho})}{\widehat{\rho}} + \fenergy(\widehat{\temp}, \widehat{\rho}) \right) \left(\rho  - \widehat{\rho} \right)
   \right]
   \, \cvolumee
   .
 \end{equation}
 (Recall that $\widehat{\fenergy} = \widehat{\ienergy} - \widehat{\temp} \widehat{\entropy}$.) We can thus reformulate our findings concerning the decay equation in the nonlinear setting, see Summary~\ref{summary:decay-equation-nonlinear-rest-state} using alternative formula for the functional $\mathcal{V}_{\mathrm{meq}, \, \widehat{\temp}, \, \widehat{\rho}}(\rho, \vec{v}, \temp)$, in particular we have characterisation shown in Summary~\ref{summary:bregman-reformulation}.

 \begin{summary}[Functional $\mathcal{V}_{\mathrm{meq}, \, \widehat{\temp}, \, \widehat{\rho}}(\rho, \vec{v}, \temp)$ rewritten in terms of Bregman distance/divergence induced by the modified internal energy function $\volienergy(\volentropy, \rho) =_{\bydefinition} \rho  \left. \ienergy (\entropy, \rho) \right|_{\entropy = \frac{\volentropy}{\rho}}$]
   \label{summary:bregman-reformulation}
   Let $\widehat{\vec{A}} = \left[ \widehat{\volentropy}, \widehat{\rho} \right] = \left[ \widehat{\rho} \entropy(\widehat{\temp}, \widehat{\rho}), \widehat{\rho} \right]$ be the pair of the entropy--density variables corresponding to the \emph{spatially homogeneous rest state}~\eqref{eq:spatially-homogenoeus-rest-state}, and let $\vec{A} = \left[\volentropy, \rho \right] = \left[ \rho \entropy(\temp, \rho), \rho \right]$ be the pair corresponding to a general state $\left[ \rho, \vec{v}, \temp \right]$. Functional~$\mathcal{V}_{\mathrm{meq}, \, \widehat{\temp}, \, \widehat{\rho}}(\rho, \vec{v}, \temp)$ introduced in Summary~\ref{summary:lagrange-multipliers}, formula~\eqref{eq:270}, can be rewritten as
 \begin{equation}
   \label{eq:298}
   \mathcal{V}_{\mathrm{meq}, \, \widehat{\temp}, \, \widehat{\rho}}(\rho, \vec{v}, \temp)
   =
   \int_{\Omega}
   \frac{1}{2} \rho \absnorm{\vec{v}}^2
   \,
   \cvolumee
   +
   \int_{\Omega}
   \bregmandivergenceop{\volienergy}{\vec{A}}{\widehat{\vec{A}}}
   \,
   \cvolumee
   ,
 \end{equation}
 where $\bregmandivergenceop{\volienergy}{\vec{A}}{\widehat{\vec{A}}}$ denotes the Bregman distance/divergence induced by $\volienergy(\volentropy, \rho)$. (The Bregman distance/divergence must be computed in its natural variables $\volentropy$ and $\rho$, and then evaluated at $\temp$ and $\rho$.)
 \end{summary}

 \subsubsection{Bregman distance/divergence and Legendre transform}
\label{sec:bregm-dist-legendre}
 The Bregman distance/divergence is closely related to the Legendre transform, or, in the language of convex analysis, to the convex conjugate. Let $f: \R^n \mapsto \R$ be a differentiable and strictly convex function. Its convex conjugate $\conjugate{f}$ is defined as
 \begin{equation}
   \label{eq:451}
   \conjugate{f}(\vec{y})
   =_{\bydefinition}
   \sup_{\vec{x} \in \R^n}
   \left(
     \vectordot{\vec{x}}{\vec{y}}
     -
     f(\vec{x})
   \right)
   .
 \end{equation}
 If $\vec{y} = \left. \pd{f}{\vec{z}} \right|_{\vec{z} = \vec{x}}$, then the supremum in~\eqref{eq:451} is attained at the point $\vec{x}$, and we get the standard Fenchel--Young equality
 \begin{equation}
   \label{eq:452}
   f(\vec{x})
   +
   \conjugate{f}(\vec{y})
   =
   \vectordot{\vec{x}}{\vec{y}}
   ,
   \qquad
   \vec{y} = \left. \pd{f}{\vec{z}} \right|_{\vec{z} = \vec{x}}
   .
 \end{equation}
 Furthermore, the gradient maps of $f$ and $\conjugate{f}$ are mutually inverse,
 \begin{equation}
   \label{eq:453}
   \vec{y} = \left. \pd{f}{\vec{z}} \right|_{\vec{z} = \vec{x}}
   \quad
   \Longleftrightarrow
   \quad
   \vec{x} = \left. \pd{\conjugate{f}}{\vec{z}} \right|_{\vec{z} = \vec{y}}
   .
 \end{equation}
 (In the present context the variables $\vec{x}$ and $\vec{y}$ are the conjugate pairs such as the entropy and the temperature, and the switch from $f$ to $\conjugate{f}$ is the switch from one thermodynamic potential to another one.) The Bregman distance/divergence induced by $f$ and the Bregman distance/divergence induced by $\conjugate{f}$ are related as follows.

 \begin{lemma}[Bregman distance/divergence and convex conjugate]
   \label{lm:5}
   Let $f: \R^n \mapsto \R$ be a differentiable and strictly convex function, and let $\conjugate{f}$ be its convex conjugate~\eqref{eq:451}. Let $\vec{x}, \widehat{\vec{x}} \in \R^n$ and let $\vec{y} =_{\bydefinition} \left. \pd{f}{\vec{z}} \right|_{\vec{z} = \vec{x}}$ and $\widehat{\vec{y}} =_{\bydefinition} \left. \pd{f}{\vec{z}} \right|_{\vec{z} = \widehat{\vec{x}}}$
   be the corresponding conjugate variables. Then
   \begin{equation}
     \label{eq:455}
     \bregmandivergenceop{f}{\vec{x}}{\widehat{\vec{x}}}
     =
     \bregmandivergenceop{\conjugate{f}}{\widehat{\vec{y}}}{\vec{y}}
     .
   \end{equation}
 \end{lemma}
 Note the order of the arguments in~\eqref{eq:455}. The Bregman distance/divergence induced by the convex conjugate $\conjugate{f}$ is evaluated with the swapped roles of the state and the reference state.

 The lemma is straightforward to prove by direct computation. By the definition of the Bregman distance/divergence, see~\eqref{eq:285}, and by~\eqref{eq:453} we have
 \begin{equation}
   \label{eq:456}
   \bregmandivergenceop{\conjugate{f}}{\widehat{\vec{y}}}{\vec{y}}
   =
   \conjugate{f}(\widehat{\vec{y}})
   -
   \conjugate{f}(\vec{y})
   -
   \vectordot{\left( \left. \pd{\conjugate{f}}{\vec{z}} \right|_{\vec{z} = \vec{y}} \right)}{\left( \widehat{\vec{y}} - \vec{y} \right)}
   =
   \conjugate{f}(\widehat{\vec{y}})
   -
   \conjugate{f}(\vec{y})
   -
   \vectordot{\vec{x}}{\left( \widehat{\vec{y}} - \vec{y} \right)}
   .
 \end{equation}
 The Fenchel--Young equality~\eqref{eq:452} implies $\conjugate{f}(\vec{y}) = \vectordot{\vec{x}}{\vec{y}} - f(\vec{x})$ and $\conjugate{f}(\widehat{\vec{y}}) = \vectordot{\widehat{\vec{x}}}{\widehat{\vec{y}}} - f(\widehat{\vec{x}})$, and the substitution into~\eqref{eq:456} yields
 \begin{equation}
   \label{eq:457}
   \bregmandivergenceop{\conjugate{f}}{\widehat{\vec{y}}}{\vec{y}}
   =
   \vectordot{\widehat{\vec{x}}}{\widehat{\vec{y}}}
   -
   f(\widehat{\vec{x}})
   -
   \vectordot{\vec{x}}{\vec{y}}
   +
   f(\vec{x})
   -
   \vectordot{\vec{x}}{\widehat{\vec{y}}}
   +
   \vectordot{\vec{x}}{\vec{y}}
   =
   f(\vec{x})
   -
   f(\widehat{\vec{x}})
   -
   \vectordot{\widehat{\vec{y}}}{\left( \vec{x} - \widehat{\vec{x}} \right)}
   =
   \bregmandivergenceop{f}{\vec{x}}{\widehat{\vec{x}}}
   ,
 \end{equation}
 where we have used $\widehat{\vec{y}} = \left. \pd{f}{\vec{z}} \right|_{\vec{z} = \widehat{\vec{x}}}$ in the last equality. This concludes the proof.
 
 \subsection{Answer to stability problem for an isolated system}
 The analysis in this section provides an answer to Question~\ref{q:isolated}. Our findings are the following\footnote{We slightly modify the notation in order to keep it in line with the notation used for Bregman distance/divergence and for later use in analysis of open systems. We thus systematically write $\widehat{\vec{W}}_{\equilibrium}$ instead of mere $\widehat{\vec{W}}$, and similarly for $\widehat{\rho}_{\equilibrium}$, $\widehat{\temp}_{\equilibrium}$ instead of mere $\widehat{\rho}$, $\widehat{\temp}$. This helps us to indicate that we are working in the thermodynamically isolated setting as we shall later need to distinguish whether the target state $\widehat{\rho}$, $\widehat{\temp}$ is the spatially homogeneous equilibrium rest state~\eqref{eq:spatially-homogenoeus-rest-state} (thermodynamically isolated setting) or the spatially inhomogeneous steady state~\eqref{eq:spatially-inhomogenoeus-steady-state-dirichlet} (thermodynamically open setting).}. 
 We have constructed the functional
 \begin{multline}
   \label{eq:302}
   \lyapunovfn{\mequilibrium}{\vec{W}}{\widehat{\vec{W}}_{\equilibrium}}
   \equiv
   \mathcal{V}_{\mathrm{meq}, \, \widehat{\temp}_{\equilibrium}, \, \widehat{\rho}_{\equilibrium}}(\rho, \vec{v}, \temp)
   =_{\bydefinition}
   \int_{\Omega} \frac{1}{2} \rho \absnorm{\vec{v}}^2 \, \cvolumee
   \\
   +
   \int_{\Omega}
   \left[
     -
     \widehat{\temp}_{\equilibrium} \left\{ \rho \entropy(\temp, \rho) - \widehat{\rho}_{\equilibrium} \entropy(\widehat{\temp}_{\equilibrium}, \widehat{\rho}_{\equilibrium}) \right\}
     +
     \left\{ \rho \ienergy(\temp, \rho) - \widehat{\rho}_{\equilibrium} \ienergy(\widehat{\temp}_{\equilibrium}, \widehat{\rho}_{\equilibrium}) \right\}
     -
     \left( \frac{\thpressure(\widehat{\temp}_{\equilibrium}, \widehat{\rho}_{\equilibrium})}{\widehat{\rho}_{\equilibrium}} + \fenergy(\widehat{\temp}_{\equilibrium}, \widehat{\rho}_{\equilibrium}) \right) \left(\rho  - \widehat{\rho}_{\equilibrium} \right)
   \right]
   \, \cvolumee
   \\
   =
   \int_{\Omega} \frac{1}{2} \rho \absnorm{\vec{v}}^2 \, \cvolumee
   +
   \int_{\Omega}
   \bregmandivergenceop{\volienergy}{\vec{A}}{\widehat{\vec{A}}_{\equilibrium}}
   \,
   \cvolumee
   ,
 \end{multline}
 where $\vec{A} = \left[ \volentropy, \rho \right] = \left[ \rho \entropy(\temp, \rho), \rho \right]$ and $\widehat{\vec{A}}_{\equilibrium} = \left[ \widehat{\rho}_{\equilibrium} \entropy(\widehat{\temp}_{\equilibrium}, \widehat{\rho}_{\equilibrium}), \widehat{\rho}_{\equilibrium} \right]$ denote the entropy--density variables corresponding to the state $\vec{W}$ and to the rest state $\widehat{\vec{W}}_{\equilibrium}$, see Summary~\ref{summary:bregman-reformulation}. The functional is nonnegative
 \begin{equation}
   \label{eq:303}
   \lyapunovfn{\mequilibrium}{\vec{W}}{\widehat{\vec{W}}_{\equilibrium}} \geq 0
 \end{equation}
 and it vanishes if and only if the state
 $
 \vec{W} = \left[ \rho, \vec{v}, \temp \right]
 $
 coincides with the spatially homogeneous equilibrium rest state
 $
 \widehat{\vec{W}}_{\equilibrium} = \left[ \rhoeq, \vec{0}, \tempeq \right],
 $
 that is
 \begin{equation}
   \label{eq:304}
   \lyapunovfn{\mequilibrium}{\vec{W}}{\widehat{\vec{W}}_{\equilibrium}} = 0 \Leftrightarrow  \vec{W} = \widehat{\vec{W}}_{\equilibrium}.
 \end{equation}
 These properties are proved in Section~\ref{sec:nonnegativity} and they are consequences of thermodynamic stability conditions regarding the structure of the Helmholtz free energy~\eqref{eq:thermodynamic-stability-conditions-conjecture}.  Furthermore, the time derivative of the functional is nonpositive,
 \begin{equation}
   \label{eq:305}
   \dd{}{t} \lyapunovfn{\mequilibrium}{\vec{W}}{\widehat{\vec{W}}_{\equilibrium}} =
   -
   \widehat{\temp}_{\equilibrium}
   \int_{\Omega}
   \left(
     \frac{
       \tilde{\lambda} \left(\divergence \vec{v}\right)^2
       +
       2 \nu \tensordot{\traceless{\gradsym}}{\traceless{\gradsym}}
     }
     {
       \temp
     }
     +
     \frac{
       \kappa
       \vectordot{\nabla \temp}{\nabla{\temp}}
     }
     {
       \temp^2
     }
   \right)
   \, \cvolumee
   .
  \end{equation}
  The nonnegativity of the functional and the nonpositivity of its time derivative make the functional ideal for monitoring the progress of state $\vec{W}$ towards the spatially homogeneous rest state $\widehat{\vec{W}}_{\equilibrium}$.

  The functional $\lyapunovfn{\mequilibrium}{\vec{W}}{\widehat{\vec{W}}_{\equilibrium}}$ defined in~\eqref{eq:302} is generated by the maximisation procedure discussed in Section~\ref{sec:ident-lagr-mult}, and the formula~\eqref{eq:302} can be also rewritten as
  \begin{subequations}
    \label{eq:307}
    \begin{equation}
      \label{eq:306}
      \lyapunovfn{\mequilibrium}{\vec{W}}{\widehat{\vec{W}}_{\equilibrium}}
      =
      -
      \left(
        \netentropy_{\widehat{\temp}_{\equilibrium}} \left(\vec{W}\right)
        -
        \netentropy_{\widehat{\temp}_{\equilibrium}} \left(\widehat{\vec{W}}_{\equilibrium}\right)
      \right)
      +
      \left(\nettenergy(\vec{W}) - \nettenergy \left(\widehat{\vec{W}}_{\equilibrium}\right) \right)
      -
      \int_{\Omega}
      \left(
        \frac{\thpressure(\widehat{\temp}_{\equilibrium}, \widehat{\rho}_{\equilibrium})}{\widehat{\rho}_{\equilibrium}}
        +
        \fenergy(\widehat{\temp}_{\equilibrium}, \widehat{\rho}_{\equilibrium})
      \right)
      \left(
        \rho
        -
        \widehat{\rho}_{\equilibrium}
      \right)
      \,
      \cvolumee,
    \end{equation}
    where the symbols $\netentropy_{\widehat{\temp}_{\equilibrium}} \left(\vec{W}\right)$ and $\nettenergy\left(\vec{W}\right)$ denote the (modified) net total entropy and the net total energy, 
    \begin{align}
      \label{eq:310}
      \netentropy_{\widehat{\temp}_{\equilibrium}} \left(\vec{W}\right)
      &=
        _{\bydefinition}
        \int_{\Omega}
        \rho
        \widehat{\temp}_{\equilibrium}
        \entropy(\vec{W})
        \,
        \cvolumee
        ,
      \\
      \label{eq:311}
      \nettenergy\left(\vec{W}\right)
      &=
        _{\bydefinition}
        \int_{\Omega}
        \left(
          \frac{1}{2} \rho \absnorm{\vec{v}}^2
          +
          \rho
          \ienergy(\vec{W})
        \right)
        \,
        \cvolumee
        .
    \end{align}
  \end{subequations}

  \begin{figure}[h]
    \centering
    \includegraphics[width=0.9\textwidth]{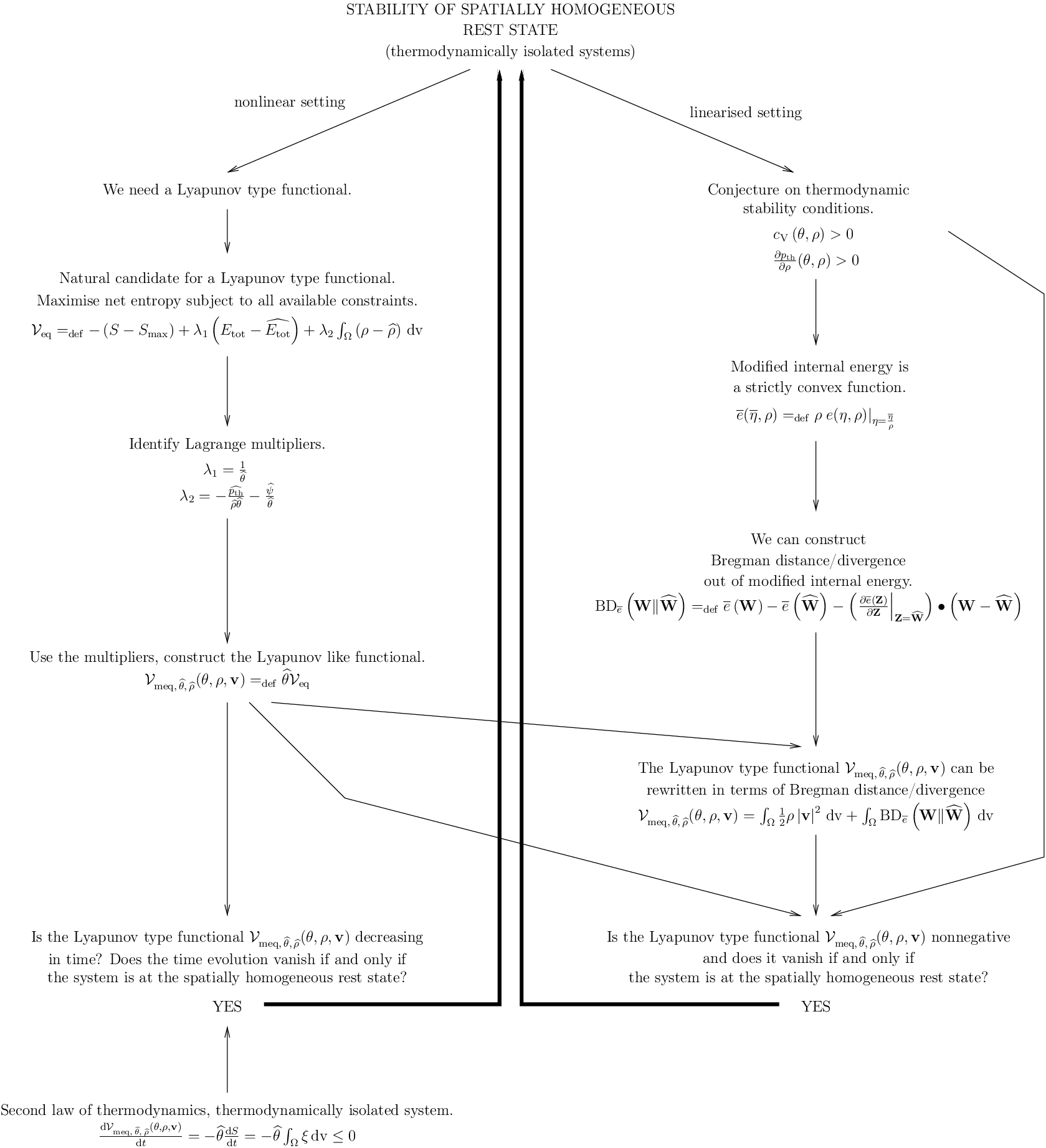}
    \caption{Graphical summary---stability of spatially homogeneous rest state in thermodynamically isolated systems.}
    \label{fig:thermodynamically-isolated-mindmap}
  \end{figure}

  \begin{summary}[Answer to Question~\ref{q:isolated}---stability of spatially homogeneous rest state in a thermodynamically isolated system]
    \label{summary:anwswer-to-stability-problem-in-isolated-system}
    Let us consider the thermodynamically isolated system introduced in Question~\ref{q:isolated}, and let $\vec{W} = \left[\rho, \vec{v}, \temp \right]$ be a generic admissible state of the system, and let $\widehat{\vec{W}}_{\equilibrium} = \left[\widehat{\rho}_{\equilibrium}, \vec{0}, \widehat{\temp}_{\equilibrium} \right] \equiv  \left[\rhoeq, \vec{0}, \tempeq \right] $ denote the spatially homogeneous equilibrium rest state~\eqref{eq:spatially-homogenoeus-rest-state}. The functional
    \begin{equation}
      \label{eq:228}
      \lyapunovfn{\mequilibrium}{\vec{W}}{\widehat{\vec{W}}_{\equilibrium}}
      =
      -
      \left(
        \netentropy_{\widehat{\temp}_{\equilibrium}} \left(\vec{W}\right)
        -
        \netentropy_{\widehat{\temp}_{\equilibrium}} \left(\widehat{\vec{W}}_{\equilibrium}\right)
      \right)
      +
      \left(\nettenergy(\vec{W}) - \nettenergy \left(\widehat{\vec{W}}_{\equilibrium}\right) \right)
      -
      \int_{\Omega}
      \left(
        \frac{\thpressure(\widehat{\temp}_{\equilibrium}, \widehat{\rho}_{\equilibrium})}{\widehat{\rho}_{\equilibrium}}
        +
        \fenergy(\widehat{\temp}_{\equilibrium}, \widehat{\rho}_{\equilibrium})
      \right)
      \left(
        \rho
        -
        \widehat{\rho}_{\equilibrium}
      \right)
      \,
      \cvolumee,
    \end{equation}
    see~\eqref{eq:307}, which has an alternative representation in terms of Bregman distance/divergence associated to the modified internal energy function
    \begin{equation}
      \label{eq:229}
      \lyapunovfn{\mequilibrium}{\vec{W}}{\widehat{\vec{W}}_{\equilibrium}}
      =
      \int_{\Omega}
      \frac{1}{2} \rho \absnorm{\vec{v}}^2
      \,
      \cvolumee
      +
      \int_{\Omega}
      \bregmandivergenceop{\volienergy}{\vec{A}}{\widehat{\vec{A}}_{\equilibrium}}
      \,
      \cvolumee
      ,
    \end{equation}
    see~\eqref{eq:298}, where $\vec{A} = \left[ \volentropy, \rho \right]$ and $\widehat{\vec{A}}_{\equilibrium} = \left[ \widehat{\volentropy}_{\equilibrium}, \widehat{\rho}_{\equilibrium} \right]$ denote the entropy--density variables corresponding to the states $\vec{W}$ and $\widehat{\vec{W}}_{\equilibrium}$, has the following properties. The functional is nonnegative
    \begin{subequations}
      \label{eq:231}
      \begin{equation}
        \label{eq:230}
        \lyapunovfn{\mequilibrium}{\vec{W}}{\widehat{\vec{W}}_{\equilibrium}} \geq 0,
      \end{equation}
      it vanishes if and only if the state $\vec{W}$ coincides with the spatially homogeneous equilibrium rest state $\widehat{\vec{W}}_{\equilibrium}$,
      \begin{equation}
        \label{eq:233}
        \lyapunovfn{\mequilibrium}{\vec{W}}{\widehat{\vec{W}}_{\equilibrium}} = 0 \Leftrightarrow  \vec{W} = \widehat{\vec{W}}_{\equilibrium},
      \end{equation}
      and the time derivative of the functional is nonpositive,
      \begin{equation}
        \label{eq:232}
        \dd{}{t} \lyapunovfn{\mequilibrium}{\vec{W}}{\widehat{\vec{W}}_{\equilibrium}} =
        -
        \widehat{\temp}_{\equilibrium}
        \int_{\Omega}
        \left(
          \frac{
            \tilde{\lambda} \left(\divergence \vec{v}\right)^2
            +
            2 \nu \tensordot{\traceless{\gradsym}}{\traceless{\gradsym}}
          }
          {
            \temp
          }
          +
          \frac{
            \kappa
            \vectordot{\nabla \temp}{\nabla{\temp}}
          }
          {
            \temp^2
          }
        \right)
        \, \cvolumee
        ,
      \end{equation}
      and the time derivative vanishes if and only if the velocity field vanishes and the temperature field is spatially homogeneous,
      \begin{equation}
        \label{eq:234}
        \dd{}{t} \lyapunovfn{\mequilibrium}{\vec{W}}{\widehat{\vec{W}}_{\equilibrium}} = 0
        \Leftrightarrow
        \vec{v} = \vec{0} \text{ and } \nabla \temp = \vec{0}.
      \end{equation}
    \end{subequations}
    The functional $\lyapunovfn{\mequilibrium}{\vec{W}}{\widehat{\vec{W}}_{\equilibrium}}$ is thus a Lyapunov type functional.
  \end{summary}

\subsection{Lyapunov type functional in various representations}
\label{sec:lyap-type-funct-2}

The reformulation of Lyapunov type functional $\mathcal{V}_{\mathrm{meq}, \, \widehat{\temp}, \, \widehat{\rho}}(\rho, \vec{v}, \temp)$ in terms of the Bregman distance/divergence, that is the formula
\begin{equation}
  \label{eq:391}
  \mathcal{V}_{\mathrm{meq}, \, \widehat{\temp}, \, \widehat{\rho}}(\rho, \vec{v}, \temp)
   =
   \int_{\Omega}
   \frac{1}{2} \rho \absnorm{\vec{v}}^2
   \,
   \cvolumee
   +
   \int_{\Omega}
   \bregmandivergenceop{\volienergy}{\vec{A}}{\widehat{\vec{A}}}
   \,
   \cvolumee
   ,
\end{equation}
see Summary~\ref{summary:bregman-reformulation}, is particularly appealing for its simplicity and relation to convex analysis. It would be nice to rewrite the kinetic energy term in terms of a Bregman distance/divergence as well. Furthermore, it would be interesting to see the functional rewritten in terms of other thermodynamics variables/potentials. We proceed with the corresponding calculations in the following sections. 

\subsubsection{Lyapunov type functional as Bregman distance/divergence induced by the total energy}
\label{sec:kinet-energy-compl}
The kinetic energy term can be incorporated into the Bregman distance/divergence framework provided that we use the momentum
\begin{equation}
  \label{eq:linear-momentum}
  \linmomentum =_{\bydefinition} \rho \vec{v}
\end{equation}
instead of the velocity, see Section~\ref{sec:suit-choice-vari}. Let us consider the variables $\vec{V} = \left[ \rho, \linmomentum, \volentropy \right]$, see~\eqref{eq:297}, and the total energy per unit volume expressed in these variables,
\begin{equation}
  \label{eq:modified-total-energy}
  \voltenergy(\rho, \linmomentum, \volentropy)
  =_{\bydefinition}
  \frac{1}{2}
  \frac{\absnorm{\linmomentum}^2}{\rho}
  +
  \volienergy(\volentropy, \rho)
  .
\end{equation}

The function $\voltenergy$ is strictly convex,\footnote{%
  The function $\voltenergy(\rho, \linmomentum, \volentropy)$ is convex with respect to all its variables simultaneously, that is as a function of the triple $\left[\rho, \linmomentum, \volentropy\right]$, $\rho > 0$. Indeed, let us denote
  \begin{equation}
    \label{eq:390}
    k(\rho, \linmomentum) =_{\bydefinition} \frac{\absnorm{\linmomentum}^2}{2 \rho},
  \end{equation}
  and let us investigate the second G\^ateaux derivative of $\voltenergy$ at a point $\widehat{\vec{V}} = \left[\widehat{\rho}, \widehat{\linmomentum}, \widehat{\volentropy}\right]$, $\widehat{\rho} > 0$, in an arbitrary direction (perturbation) $\widetilde{\vec{V}} = \left[\widetilde{\rho}, \widetilde{\linmomentum}, \widetilde{\volentropy}\right]$. The second derivatives of $k$ read
  \begin{align*}
    \ppd{k}{\rho}(\rho, \linmomentum)
    &=
      \frac{\absnorm{\linmomentum}^2}{\rho^3}
      ,
    \\
    \pd{^2 k}{\rho \partial \linmomentum}(\rho, \linmomentum)
    &=
      -
      \frac{\linmomentum}{\rho^2}
      ,
    \\
    \pd{^2 k}{\linmomentum \partial \linmomentum}(\rho, \linmomentum)
    &=
      \frac{1}{\rho} \identity
      ,
  \end{align*}
  hence the second G\^ateaux derivative of $k$ at the point $\left[\widehat{\rho}, \widehat{\linmomentum}\right]$ in the direction $\left[\widetilde{\rho}, \widetilde{\linmomentum}\right]$ reads
  \begin{equation*}
    \Diff^2 k(\widehat{\rho}, \widehat{\linmomentum})
    \left[
      \begin{bmatrix}
        \widetilde{\rho} \\
        \widetilde{\linmomentum}
      \end{bmatrix}
      ,
      \begin{bmatrix}
        \widetilde{\rho} \\
        \widetilde{\linmomentum}
      \end{bmatrix}
    \right]
    =
    \vectordot{
      \begin{bmatrix}
        \widetilde{\rho} \\
        \widetilde{\linmomentum}
      \end{bmatrix}
    }
    {
    \begin{bmatrix}
        \frac{\absnorm{\widehat{\linmomentum}}^2}{\widehat{\rho}^3} & - \transpose{\left( \frac{\widehat{\linmomentum}}{\widehat{\rho}^2} \right)} \\
      - \frac{\widehat{\linmomentum}}{\widehat{\rho}^2} & \frac{1}{\widehat{\rho}} \identity
    \end{bmatrix}
    \begin{bmatrix}
      \widetilde{\rho} \\
      \widetilde{\linmomentum}
    \end{bmatrix}
    }
    =
    \frac{\absnorm{\widehat{\linmomentum}}^2}{\widehat{\rho}^3} \widetilde{\rho}^2
    -
    2 \widetilde{\rho} \frac{\vectordot{\widehat{\linmomentum}}{\widetilde{\linmomentum}}}{\widehat{\rho}^2}
    +
    \frac{\absnorm{\widetilde{\linmomentum}}^2}{\widehat{\rho}}
    =
    \frac{1}{\widehat{\rho}}
    \absnorm{\widetilde{\linmomentum} - \widetilde{\rho} \frac{\widehat{\linmomentum}}{\widehat{\rho}}}^2
    \geq
    0
    ,
  \end{equation*}
  and it is positive semidefinite. The function~$k(\rho, \linmomentum)$ is thus convex as a function of the pair $\left[\rho, \linmomentum\right]$, but it is not strictly convex---the quadratic form vanishes in the directions~$\widetilde{\linmomentum} = \widetilde{\rho} \frac{\widehat{\linmomentum}}{\widehat{\rho}}$.

  Equivalently, we could have concluded on the convexity of~$k$ by an abstract argument. The function~$k(\rho, \linmomentum)$ can be rewritten as
  \begin{equation*}
    k(\rho, \linmomentum) = \rho \, g\left(\frac{\linmomentum}{\rho}\right)
  \end{equation*}
  with $g(\vec{u}) = \frac{1}{2} \absnorm{\vec{u}}^2$ hence $k(\rho, \linmomentum)$ is the perspective of a convex function, and perspectives of convex functions are convex, see Remark~\ref{rm:perspective}.

  Having analysed the second G\^ateaux derivative of the kinetic energy term, we can proceed with the total modified energy $\voltenergy$. The second G\^ateaux derivative of $\voltenergy$ at the point $\widehat{\vec{V}}$ in the direction $\widetilde{\vec{V}}$ reads
  \begin{equation*}
    \Diff^2 \voltenergy(\widehat{\vec{V}})
    \left[
      \widetilde{\vec{V}}, \widetilde{\vec{V}}
    \right]
    =
    \frac{1}{\widehat{\rho}}
    \absnorm{\widetilde{\linmomentum} - \widetilde{\rho} \frac{\widehat{\linmomentum}}{\widehat{\rho}}}^2
    +
    \vectordot{
      \begin{bmatrix}
        \widetilde{\volentropy} \\
        \widetilde{\rho}
      \end{bmatrix}
    }
    {
      \hessian{\volienergy}(\widehat{\volentropy}, \widehat{\rho})
      \begin{bmatrix}
        \widetilde{\volentropy} \\
        \widetilde{\rho}
      \end{bmatrix}
    }
    ,
  \end{equation*}
  where $\hessian{\volienergy}$ is the Hessian matrix~\eqref{eq:206} of the modified internal energy with respect to the variables $\volentropy$ and $\rho$. The first term is nonnegative, and the second term is positive unless $\widetilde{\volentropy} = 0$ and $\widetilde{\rho} = 0$, since $\hessian{\volienergy}$ is positive definite, see Summary~\ref{summary:convexity}. If $\widetilde{\volentropy} = 0$ and $\widetilde{\rho} = 0$, then the first term reduces to~$\frac{1}{\widehat{\rho}} \absnorm{\widetilde{\linmomentum}}^2$, which is positive unless $\widetilde{\linmomentum} = \vec{0}$. The second G\^ateaux derivative is thus positive in every nonzero direction $\widetilde{\vec{V}}$, which means that the Hessian matrix of $\voltenergy$ is positive definite, and $\voltenergy$ is a strictly convex function of $\left[\rho, \linmomentum, \volentropy\right]$. In other words, the strict convexity of the modified internal energy $\volienergy(\volentropy, \rho)$ compensates for the missing strict convexity of the kinetic energy term.
} 
since the kinetic energy $\frac{\absnorm{\linmomentum}^2}{2 \rho}$ is a convex function of the density--momentum pair $\left[\rho, \linmomentum\right]$, $\rho > 0$, and the modified internal energy $\volienergy(\volentropy, \rho)$ is strictly convex by Summary~\ref{summary:convexity}. At the spatially homogeneous rest state $\widehat{\vec{V}} = \left[ \widehat{\rho}, \vec{0}, \widehat{\volentropy} \right]$ the gradient of the kinetic energy term vanishes, hence
\begin{equation}
  \label{eq:355}
  \bregmandivergenceop{\voltenergy}{\vec{V}}{\widehat{\vec{V}}}
  =
  \frac{1}{2}
  \frac{\absnorm{\linmomentum}^2}{\rho}
  +
  \bregmandivergenceop{\volienergy}{\vec{A}}{\widehat{\vec{A}}}
  =
  \frac{1}{2}
  \rho
  \absnorm{\vec{v}}^2
  +
  \bregmandivergenceop{\volienergy}{\vec{A}}{\widehat{\vec{A}}}
  ,
\end{equation}
where $\vec{A} = \left[ \volentropy, \rho \right]$ is the pair of the entropy--density variables, and the Lyapunov type functional~\eqref{eq:229} can be compactly written as
\begin{equation}
  \label{eq:356}
  \lyapunovfn{\mequilibrium}{\vec{W}}{\widehat{\vec{W}}_{\equilibrium}}
  =
  \int_{\Omega}
  \bregmandivergenceop{\voltenergy}{\vec{V}}{\widehat{\vec{V}}_{\equilibrium}}
  \,
  \cvolumee
  .
\end{equation}
Note that the natural variables on the right-hand side in~\eqref{eq:356} are $\vec{V} = \left[ \rho, \linmomentum, \volentropy \right]$ and $\widehat{\vec{V}}_{\equilibrium} = \left[ \widehat{\rho}_{\equilibrium}, \widehat{\linmomentum}_{\equilibrium}, \widehat{\volentropy}_{\equilibrium} \right]$, meaning that the Bregman distance/divergence must be computed in its natural variables $\linmomentum$, $\volentropy$ and $\rho$, and then evaluated at $\vec{v}$, $\temp$ and $\rho$, provided that we want the functional on the left-hand side written as $\mathcal{V}_{\mathrm{meq}, \, \widehat{\temp}, \, \widehat{\rho}}(\rho, \vec{v}, \temp)$.

\subsubsection{Lyapunov type functional as Bregman distance/divergence induced by the generalised entropy}
\label{sec:lyap-type-funct-3}
Having introduced the total modified energy, see~\eqref{eq:modified-total-energy}, we can formally solve the equation
\begin{equation}
  \label{eq:338}
  \voltenergy = \voltenergy(\rho, \linmomentum, \volentropy)
\end{equation}
for $\volentropy$ and we call this function the generalised entropy $\volentropygen$. In this way we can now interpret $\volentropygen$ as a function of the \emph{full} triple~$\left[\rho, \linmomentum, \voltenergy \right]$, that is
\begin{equation}
  \label{eq:generalised-entropy}
  \volentropygen = _{\bydefinition} \volentropygen \left( \rho, \linmomentum, \voltenergy \right).
\end{equation}
This is an analogue of the standard thermodynamic manipulation that inverts $\volienergy = \volienergy(\volentropy, \rho)$ as $\volentropy = \volentropy (\volienergy, \rho)$. In particular, as
\begin{equation}
  \label{eq:393}
  \voltenergy(\rho, \linmomentum, \volentropy)
  =
  \frac{1}{2}
  \frac{\absnorm{\linmomentum}^2}{\rho}
  +
  \volienergy(\volentropy, \rho),
\end{equation}
we see that
\begin{equation}
  \label{eq:394}
  \volentropygen \left( \rho, \linmomentum, \voltenergy \right)
  =
  \volentropy\left( \voltenergy - \frac{1}{2}\frac{\absnorm{\linmomentum}^2}{\rho}, \rho\right)
  .
\end{equation}
The generalised entropy $\volentropygen$ is a strictly concave function of the triple $\left[\rho, \linmomentum, \voltenergy \right]$. Equivalently, $- \volentropygen$ is a strictly convex function of the triple $\left[\rho, \linmomentum, \voltenergy \right]$.

The fact that the generalised entropy $\volentropygen$ is a strictly concave function of the triple $\left[\rho, \linmomentum, \voltenergy \right]$ can be shown as follows. Let $\widehat{\vec{U}} = \left[\widehat{\rho}, \widehat{\linmomentum}, \widehat{\voltenergy}\right]$ be a point and let $\widetilde{\vec{U}} = \left[\widetilde{\rho}, \widetilde{\linmomentum}, \widetilde{\voltenergy}\right]$ be a perturbation. Formula~\eqref{eq:394} shows that $\volentropygen$ is the composition
  \begin{equation}
    \label{eq:428}
    \volentropygen(\vec{U})
    =
    \volentropy\left( \volienergy(\vec{U}), \rho \right)
    ,
    \qquad
    \volienergy(\vec{U})
    =_{\bydefinition}
    \voltenergy - \frac{1}{2}\frac{\absnorm{\linmomentum}^2}{\rho}
    ,
  \end{equation}
  where $\volentropy(\volienergy, \rho)$ is the modified entropy as a function of its natural variables $\volienergy$ and $\rho$. The chain rule for the second derivative of a composite function yields
  \begin{equation}
    \label{eq:429}
    \Diff^2 \volentropygen(\widehat{\vec{U}})
    \left[
      \widetilde{\vec{U}}, \widetilde{\vec{U}}
    \right]
    =
    \vectordot{
      \begin{bmatrix}
        \widetilde{\volienergy} \\
        \widetilde{\rho}
      \end{bmatrix}
    }
    {
      \hessian{\volentropy}(\widehat{\volienergy}, \widehat{\rho})
      \begin{bmatrix}
        \widetilde{\volienergy} \\
        \widetilde{\rho}
      \end{bmatrix}
    }
    +
    \pd{\volentropy}{\volienergy}(\widehat{\volienergy}, \widehat{\rho})
    \,
    \Diff^2 \volienergy(\widehat{\vec{U}})
    \left[
      \widetilde{\vec{U}}, \widetilde{\vec{U}}
    \right]
    ,
  \end{equation}
  where $\hessian{\volentropy}(\volienergy, \rho)$ is the Hessian matrix of $\volentropy(\volienergy, \rho)$ with respect to $\volienergy$ and $\rho$, and where
  \begin{equation}
    \label{eq:430}
    \widetilde{\volienergy}
    =_{\bydefinition}
    \Diff \volienergy(\widehat{\vec{U}})
    \left[
      \widetilde{\vec{U}}
    \right]
    =
    \begin{bmatrix}
      \frac{1}{2}
      \frac{\absnorm{\widehat{\linmomentum}}^2}{\widehat{\rho}^2}
      &
      - \transpose{\left( \frac{\widehat{\linmomentum}}{\widehat{\rho}} \right)}
      &
      1
    \end{bmatrix}
    \begin{bmatrix}
      \widetilde{\rho} \\
      \widetilde{\linmomentum} \\
      \widetilde{\voltenergy}
    \end{bmatrix}
    =
    \widetilde{\voltenergy}
    -
    \vectordot{\frac{\widehat{\linmomentum}}{\widehat{\rho}}}{\widetilde{\linmomentum}}
    +
    \frac{1}{2}
    \frac{\absnorm{\widehat{\linmomentum}}^2}{\widehat{\rho}^2}
    \widetilde{\rho}
  \end{equation}
  is the (first order) perturbation of the modified internal energy induced by the perturbation $\widetilde{\vec{U}}$. The second derivative of~$\volienergy(\vec{U})$ is, up to the sign, the second derivative of the kinetic energy term $k(\rho, \linmomentum) = \frac{1}{2}\frac{\absnorm{\linmomentum}^2}{\rho}$, see the discussion after~\eqref{eq:modified-total-energy},
  \begin{equation}
    \label{eq:431}
    \Diff^2 \volienergy(\widehat{\vec{U}})
    \left[
      \widetilde{\vec{U}}, \widetilde{\vec{U}}
    \right]
    =
    -
    \vectordot{
      \begin{bmatrix}
        \widetilde{\rho} \\
        \widetilde{\linmomentum} \\
        \widetilde{\voltenergy}
      \end{bmatrix}
    }
    {
      \begin{bmatrix}
        \frac{\absnorm{\widehat{\linmomentum}}^2}{\widehat{\rho}^3} & - \transpose{\left( \frac{\widehat{\linmomentum}}{\widehat{\rho}^2} \right)} & 0 \\
        - \frac{\widehat{\linmomentum}}{\widehat{\rho}^2} & \frac{1}{\widehat{\rho}} \identity & \vec{0} \\
        0 & \transpose{\vec{0}} & 0
      \end{bmatrix}
      \begin{bmatrix}
        \widetilde{\rho} \\
        \widetilde{\linmomentum} \\
        \widetilde{\voltenergy}
      \end{bmatrix}
    }
    =
    -
    \frac{1}{\widehat{\rho}}
    \absnorm{\widetilde{\linmomentum} - \widetilde{\rho} \frac{\widehat{\linmomentum}}{\widehat{\rho}}}^2
    ,
  \end{equation}
  and we also know that
  $
  \pd{\volentropy}{\volienergy}(\widehat{\volienergy}, \widehat{\rho})
  =
  \frac{1}{\widehat{\temp}}
  $.
  Consequently, we see that
  \begin{equation}
    \label{eq:432}
    \Diff^2 \volentropygen(\widehat{\vec{U}})
    \left[
      \widetilde{\vec{U}}, \widetilde{\vec{U}}
    \right]
    =
    \vectordot{
      \begin{bmatrix}
        \widetilde{\volienergy} \\
        \widetilde{\rho}
      \end{bmatrix}
    }
    {
      \hessian{\volentropy}(\widehat{\volienergy}, \widehat{\rho})
      \begin{bmatrix}
        \widetilde{\volienergy} \\
        \widetilde{\rho}
      \end{bmatrix}
    }
    -
    \frac{1}{\widehat{\temp} \widehat{\rho}}
    \absnorm{\widetilde{\linmomentum} - \widetilde{\rho} \frac{\widehat{\linmomentum}}{\widehat{\rho}}}^2
    .
  \end{equation}
  Both terms are nonpositive as the second derivative matrix $\hessian{\volentropy}(\widehat{\volienergy}, \widehat{\rho})$ is negative definite. (This follows from Lemma~\ref{lm:4} and the positive definiteness of Hessian $\hessian{\volienergy}(\widehat{\volentropy}, \widehat{\rho})$.)  The sum in~\eqref{eq:432} vanishes only if $\widetilde{\volienergy} = 0$, $\widetilde{\rho} = 0$ and $\widetilde{\linmomentum} - \widetilde{\rho} \frac{\widehat{\linmomentum}}{\widehat{\rho}} = \vec{0}$, that is if $\widetilde{\rho} = 0$ and $\widetilde{\linmomentum} = \vec{0}$, and then the formula for $\widetilde{\volienergy}$ reduces to $\widetilde{\volienergy} = \widetilde{\voltenergy}$, hence also $\widetilde{\voltenergy} = 0$. The Hessian matrix of $\volentropygen$ is thus negative definite, and $\volentropygen$ is a strictly concave function of $\left[\rho, \linmomentum, \voltenergy\right]$.

  Since the (negative) generalised entropy $- \volentropygen$ is a strictly convex function of the triple $\left[\rho, \linmomentum, \voltenergy \right]$, \emph{we can construct the Bregman distance/divergence out of it}. We show that the Bregman distance/divergence induced by $- \volentropygen$ coincides, up to a constant factor, with the Bregman distance/divergence induced by the total modified energy $\voltenergy$, and that this holds \emph{pointwise}, that is on the level of integrands.

  Let $\vec{V} = \left[ \rho, \linmomentum, \volentropy \right]$ be a state and let $\widehat{\vec{V}} = \left[ \widehat{\rho}, \widehat{\linmomentum}, \widehat{\volentropy} \right]$ be an arbitrary reference state. We denote by
  \begin{equation}
    \label{eq:435}
    \vec{U}(\vec{V})
    =_{\bydefinition}
    \left[ \rho, \linmomentum, \voltenergy(\rho, \linmomentum, \volentropy) \right]
  \end{equation}
  the change of variables from $\vec{V} = \left[ \rho, \linmomentum, \volentropy \right]$ to the variables $\vec{U} = \left[ \rho, \linmomentum, \voltenergy \right]$, and we write $\vec{U} = \vec{U}(\vec{V})$, $\widehat{\vec{U}} = \vec{U}(\widehat{\vec{V}}) = \left[ \widehat{\rho}, \widehat{\linmomentum}, \widehat{\voltenergy} \right]$ with $\widehat{\voltenergy} = \voltenergy(\widehat{\rho}, \widehat{\linmomentum}, \widehat{\volentropy})$. By the definition of the generalised entropy we have $\volentropygen(\vec{U}(\vec{V})) = \volentropy$, in particular $\volentropygen(\widehat{\vec{U}}) = \widehat{\volentropy}$. Furthermore, we denote $\widehat{\vec{v}} =_{\bydefinition} \frac{\widehat{\linmomentum}}{\widehat{\rho}}$, and $\widehat{\temp}$ and $\widehat{\gibbs}$ denote the temperature and the specific Gibbs free energy evaluated at $\left[ \widehat{\volentropy}, \widehat{\rho} \right]$.

  First, we compute the Bregman distance/divergence induced by the total modified energy $\voltenergy(\rho, \linmomentum, \volentropy)$. The partial derivatives of~\eqref{eq:modified-total-energy} are, see~\eqref{eq:341},
  \begin{subequations}
    \label{eq:436}
    \begin{align}
      \label{eq:441}
      \pd{\voltenergy}{\rho}(\vec{V})
      &=
        \gibbs
        -
        \frac{1}{2}
        \frac{\absnorm{\linmomentum}^2}{\rho^2}
        ,
      \\
      \label{eq:442}
      \pd{\voltenergy}{\linmomentum}(\vec{V})
      &=
        \frac{\linmomentum}{\rho}
        ,
      \\
      \label{eq:443}
      \pd{\voltenergy}{\volentropy}(\vec{V})
      &=
        \temp
        ,
    \end{align}
  \end{subequations}
  hence the Bregman distance/divergence induced by $\voltenergy$, see Definition~\ref{dfn:bregman-divergence}, reads
  \begin{equation}
    \label{eq:437}
    \bregmandivergenceop{\voltenergy}{\vec{V}}{\widehat{\vec{V}}}
    =
    \left( \voltenergy - \widehat{\voltenergy} \right)
    -
    \left( \widehat{\gibbs} - \frac{1}{2} \absnorm{\widehat{\vec{v}}}^2 \right)
    \left( \rho - \widehat{\rho} \right)
    -
    \vectordot{\widehat{\vec{v}}}{\left( \linmomentum - \widehat{\linmomentum} \right)}
    -
    \widehat{\temp}
    \left( \volentropy - \widehat{\volentropy} \right)
    ,
  \end{equation}
  where $\voltenergy = \voltenergy(\rho, \linmomentum, \volentropy)$.

  Second, we compute the Bregman distance/divergence induced by $- \volentropygen(\rho, \linmomentum, \voltenergy)$. The partial derivatives of $\volentropygen$ are obtained by the differentiation of the identity $\voltenergy \left( \rho, \linmomentum, \volentropygen(\rho, \linmomentum, \voltenergy) \right) = \voltenergy$, and~\eqref{eq:436} yields
  \begin{subequations}
    \label{eq:438}
    \begin{align}
      \label{eq:444}
      \pd{\volentropygen}{\rho}(\vec{U})
      &=
        -
        \frac{1}{\temp}
        \left(
          \gibbs
          -
          \frac{1}{2}
          \frac{\absnorm{\linmomentum}^2}{\rho^2}
        \right)
        ,
      \\
      \label{eq:445}
      \pd{\volentropygen}{\linmomentum}(\vec{U})
      &=
        -
        \frac{1}{\temp}
        \frac{\linmomentum}{\rho}
        ,
      \\
      \label{eq:446}
      \pd{\volentropygen}{\voltenergy}(\vec{U})
      &=
        \frac{1}{\temp}
        .
    \end{align}
  \end{subequations}
  Consequently, the Bregman distance/divergence induced by $- \volentropygen$ reads
  \begin{multline}
    \label{eq:439}
    \bregmandivergenceop{- \volentropygen}{\vec{U}}{\widehat{\vec{U}}}
    =
    -
    \volentropygen(\vec{U})
    +
    \volentropygen(\widehat{\vec{U}})
    +
    \vectordot{
      \left(
        \left.
          \pd{\volentropygen}{\vec{U}}
        \right|_{\vec{U} = \widehat{\vec{U}}}
      \right)
    }
    {
      \left( \vec{U} - \widehat{\vec{U}} \right)
    }
    \\
    =
    -
    \left( \volentropy - \widehat{\volentropy} \right)
    +
    \frac{1}{\widehat{\temp}}
    \left(
      \left( \voltenergy - \widehat{\voltenergy} \right)
      -
      \left( \widehat{\gibbs} - \frac{1}{2} \absnorm{\widehat{\vec{v}}}^2 \right)
      \left( \rho - \widehat{\rho} \right)
      -
      \vectordot{\widehat{\vec{v}}}{\left( \linmomentum - \widehat{\linmomentum} \right)}
    \right)
    ,
  \end{multline}
  where we have used $\volentropygen(\vec{U}) = \volentropy$ and $\volentropygen(\widehat{\vec{U}}) = \widehat{\volentropy}$.

  Finally, if we multiply~\eqref{eq:439} by $\widehat{\temp}$ and compare the result with~\eqref{eq:437}, we see that
  \begin{equation}
    \label{eq:440}
    \bregmandivergenceop{\voltenergy}{\vec{V}}{\widehat{\vec{V}}}
    =
    \widehat{\temp}
    \,
    \bregmandivergenceop{- \volentropygen}{\vec{U}(\vec{V})}{\vec{U}(\widehat{\vec{V}})}
    .
  \end{equation}
  The Bregman distance/divergence induced by the total modified energy and the Bregman distance/divergence induced by the negative of the generalised entropy thus coincide pointwise, up to the factor $\widehat{\temp}$ given by the temperature of the reference state. Note that~\eqref{eq:440} holds for an arbitrary reference state $\widehat{\vec{V}}$, in particular we do not need $\widehat{\linmomentum} = \vec{0}$. If the reference state is the spatially homogeneous equilibrium rest state $\widehat{\vec{V}}_{\equilibrium}$, then $\widehat{\vec{v}} = \vec{0}$, the factor $\widehat{\temp} = \widehat{\temp}_{\equilibrium}$ is a constant, and the integration of~\eqref{eq:440} over the domain $\Omega$ yields the equality of the two representations of the Lyapunov type functional $\lyapunovfn{\mequilibrium}{\vec{W}}{\widehat{\vec{W}}_{\equilibrium}}$ stated in Summary~\ref{summary:lyapunov-bregman-representations}.

  \begin{summary}[Lyapunov type functional as Bregman distance/divergence in different representations]
  \label{summary:lyapunov-bregman-representations}
  Let $\widehat{\vec{V}}_{\equilibrium}$ be the spatially homogeneous equilibrium rest state~\eqref{eq:spatially-homogenoeus-rest-state}. The Lyapunov type functional $\lyapunovfn{\mequilibrium}{\vec{W}}{\widehat{\vec{W}}_{\equilibrium}}$, see~\eqref{eq:302}, has two equivalent representations in terms of Bregman distance/divergence,
  \begin{subequations}
    \label{eq:433}
    \begin{align}
      \label{eq:449}
      \lyapunovfn{\mequilibrium}{\vec{W}}{\widehat{\vec{W}}_{\equilibrium}}
      &=
      \int_{\Omega}
      \bregmandivergenceop{\voltenergy}{\vec{V}}{\widehat{\vec{V}}_{\equilibrium}}
      \,
      \cvolumee
      ,
      \\
      \label{eq:450}
      \lyapunovfn{\mequilibrium}{\vec{W}}{\widehat{\vec{W}}_{\equilibrium}}
      &=
      \widehat{\temp}_{\equilibrium}
      \int_{\Omega}
      \bregmandivergenceop{- \volentropygen}{\vec{U}}{\widehat{\vec{U}}_{\equilibrium}}
      \,
      \cvolumee
      ,
    \end{align}
  \end{subequations}
  where
  \begin{equation}
    \label{eq:423}
    \voltenergy(\rho, \linmomentum, \volentropy)
    =_{\bydefinition}
    \frac{1}{2}
    \frac{\absnorm{\linmomentum}^2}{\rho} + \volienergy(\volentropy, \rho)
  \end{equation}
  denotes the total modified energy $\voltenergy$ as a function of variables $\vec{V} = \left[ \rho, \linmomentum, \volentropy \right]$, and
  \begin{equation}
    \label{eq:434}
    \volentropygen \left( \rho, \linmomentum, \voltenergy \right)
    =_{\bydefinition}
    \volentropy\left( \voltenergy - \frac{1}{2}\frac{\absnorm{\linmomentum}^2}{\rho}, \rho\right)
  \end{equation}
  denotes the generalised entropy $\volentropygen$ as a function of variables $\vec{U} = \left[ \rho, \linmomentum, \voltenergy \right]$. 
\end{summary}

\subsection{Leading order expansion of Lyapunov type functional}
\label{sec:lead-order-expans}
In Section~\ref{sec:gateaux-derivatives} we have found, by a direct but tedious calculation in the temperature--density representation, that the leading order term in the expansion of the Lyapunov type functional $\lyapunovfn{\mequilibrium}{\vec{W}}{\widehat{\vec{W}}_{\equilibrium}}$ in the neighbourhood of the spatially homogeneous rest state is the quadratic functional~\eqref{eq:127}. Up to the kinetic energy term, this is the functional known from the stability analysis of the linearised governing equations, see~\eqref{eq:57}. The representation of the Lyapunov type functional in terms of the Bregman distance/divergence, see~\eqref{eq:229}, allows us to obtain this result almost immediately, and, more importantly, it allows us to give another interpretation to the quadratic functional~\eqref{eq:127}.

In particular, we show that the leading order (quadratic) term in the Bregman distance/divergence induced by the modified internal energy $\volienergy(\volentropy, \rho)$ is given by the Hessian matrix of $\volienergy(\volentropy, \rho)$, and that the very same quadratic term is, up to the factor $- \frac{\widehat{\temp}_{\equilibrium}}{2}$, the Hessian matrix of $\volentropy(\volienergy, \rho)$ with respect to its natural variables. (The same holds for the Bregman distance/divergence induced by the total modified energy $\voltenergy(\rho, \linmomentum, \volentropy)$ and the generalised entropy $\volentropygen(\rho, \linmomentum, \voltenergy)$, which is the most abstract framework including also the linear momentum.) This can be further generalised to other thermodynamic potentials such as the Gibbs potential; the key tool would be Lemma~\ref{lm:5}.

\subsubsection{Transformation of Hessian with respect to a change of variables}
\label{sec:transf-hess-with}
Before we proceed with the computations, let us prove a simple yet useful lemma.

\begin{lemma}[Transformation of Hessian with respect to a change of variables]
  \label{lm:4}
  Let us consider the function
$
    \volienergy = \volienergy \left(\volentropy, \vec{\xi}\right)
$,
  and the corresponding inverse function with respect to the first variable
$
    \volentropy = \volentropy(\volienergy, \vec{\xi})
$,
  where $\vec{\xi}$ is an arbitrary vector. The inverse relation is understood in the sense that
  \begin{equation}
    \label{eq:400}
    \volienergy = \volienergy \left(\volentropy(\volienergy, \vec{\xi}), \vec{\xi}\right).
  \end{equation}
  Then the relation between the Hessian matrix for the original function $\hessian{\volienergy}(\volentropy, \vec{\xi})$ and the Hessian matrix  $\hessian{\volentropy}(\volienergy, \vec{\xi})$ for the inverse function reads
  \begin{equation}
    \label{eq:396}
    \hessian{\volentropy}(\volienergy, \vec{\xi})
    =
    \left.
      \left(
        -
        \frac{1}{\pd{\volienergy}{\volentropy}(\volentropy, \vec{\xi})}
        \transpose{\tensorq{J}} (\volentropy, \vec{\xi})
        \hessian{\volienergy}(\volentropy, \vec{\xi})
        \tensorq{J} (\volentropy, \vec{\xi})
      \right)
    \right|_{\left[\volentropy, \vec{\xi}\right] = \left[\volentropy(\volienergy, \vec{\xi}), \vec{\xi}\right]}
    ,
  \end{equation}
  where
  \begin{equation}
    \label{eq:399}
    \tensorq{J}
    =
    \begin{bmatrix}
      \frac{1}{\pd{\volienergy}{\volentropy}(\volentropy, \vec{\xi})} & -\frac{\transpose{\pd{\volienergy}{\vec{\xi}}(\volentropy, \vec{\xi})}}{\pd{\volienergy}{\volentropy}(\volentropy, \vec{\xi})}   \\
      \vec{0} & \identity
    \end{bmatrix}
    .
  \end{equation}
\end{lemma}

We are in particular interested in the setting where $\vec{\xi}$ has just one component, $\vec{\xi} = \rho$. This setting corresponds to the change of variables from the energy--density representation to the entropy--density representation.

The lemma is straightforward to prove by direct computation. Let us consider the transformation of variables in the form
  \begin{equation}
  \label{eq:395}
  \vec{\varphi} :
  \vec{x}
  =
  \begin{bmatrix}
    \volienergy \\
    \vec{\xi}
  \end{bmatrix}
  \mapsto
  \vec{y}
  =
  \vec{\varphi}(\vec{x})
  =
  \begin{bmatrix}
    \volentropy(\volienergy, \vec{\xi}) \\
    \vec{\xi}
  \end{bmatrix}.
\end{equation}
Direct differentiation of the inverse function relation
$
  \volienergy = \volienergy \left(\volentropy(\volienergy, \vec{\xi}), \vec{\xi}\right)
$,
with respect to $\volienergy$ and $\vec{\xi}$ yields
\begin{subequations}
  \label{eq:402}
  \begin{align}
    \label{eq:403}
    1 &= \left. \pd{\volienergy}{\volentropy}(\volentropy, \vec{\xi}) \right|_{\volentropy = \volentropy(\volienergy, \vec{\xi})} \pd{\volentropy}{\volienergy}(\volienergy, \vec{\xi}), \\
    \label{eq:404}
    \vec{0} &= \left. \pd{\volienergy}{\volentropy}(\volentropy, \vec{\xi}) \right|_{\volentropy = \volentropy(\volienergy, \vec{\xi})} \pd{\volentropy}{\vec{\xi}}(\volienergy, \vec{\xi}) +  \left. \pd{\volienergy}{\vec{\xi}}(\volentropy, \vec{\xi}) \right|_{\volentropy = \volentropy(\volienergy, \vec{\xi})}
    .
  \end{align}
\end{subequations}
We can write the outcome of this computation as a substitution rule
\begin{equation}
  \label{eq:405}
  \begin{bmatrix}
    \pd{\volienergy}{\volentropy} \\
    \pd{\volienergy}{\vec{\xi}}
  \end{bmatrix}
  =
  \begin{bmatrix}
    \frac{1}{\pd{\volentropy}{\volienergy}} \\
    -
    \frac{
    \pd{\volentropy}{\vec{\xi}}
    }
    {
    \pd{\volentropy}{\volienergy}
    }
  \end{bmatrix}
\end{equation}
or vice versa
\begin{equation}
  \label{eq:406}
    \begin{bmatrix}
    \pd{\volentropy}{\volienergy} \\
    \pd{\volentropy}{\vec{\xi}}
  \end{bmatrix}
  =
  \begin{bmatrix}
    \frac{1}{\pd{\volienergy}{\volentropy}} \\
    -
    \frac{
    \pd{\volienergy}{\vec{\xi}}
    }
    {
    \pd{\volienergy}{\volentropy}
    }
  \end{bmatrix}
  ,
\end{equation}
where we, for the sake of notational simplicity, tacitly assume that all functions are functions of their natural variables, and that the left-hand side is evaluated at the point corresponding to the point at the right-hand side as in~\eqref{eq:402}. 

The chain rule for a generic function $f(\vec{y})$ yields
\begin{equation}
  \label{eq:407}
  \pd{}{\vec{x}} \left( f(\vec{y}\left(\vec{x}\right)) \right)
  =
  \transpose{
    \left(
      \pd{\vec{y}}{\vec{x}}
    \right)
  }
  \left.
    \pd{f}{\vec{y}}(\vec{y})
  \right|_{\vec{y} = \vec{y}(\vec{x})}
  ,
\end{equation}
where the transformation matrix $\pd{\vec{y}}{\vec{x}}$ is
\begin{equation}
  \label{eq:408}
  \pd{\vec{y}}{\vec{x}}
  =
  \pd{\vec{\varphi}(\vec{x})}{\vec{x}}
  =
  \begin{bmatrix}
    \pd{\volentropy}{\volienergy} & \transpose{\left( \pd{\volentropy}{\vec{\xi}} \right)} \\
    \vec{0} & \identity
  \end{bmatrix}
  =
  \begin{bmatrix}
    \frac{1}{\pd{\volienergy}{\volentropy}} &
                                              -
                                              \frac{
                                              \transpose{ \left( \pd{\volienergy}{\vec{\xi}} \right) }
                                              }
                                              {
                                              \pd{\volienergy}{\volentropy}
                                              }
    \\
    \vec{0} & \identity
  \end{bmatrix}
  =
  \tensorq{J}
  ,
\end{equation}
where the last equality follows from~\eqref{eq:406}. Now we can take the second derivative of a generic function $f(\vec{y})$, and we get
\begin{equation}
  \label{eq:354}
  \pd{^2 f}{\vec{x} \partial \vec{x}}(\vec{x})
  =
  \pd{^2}{\vec{x} \partial \vec{x}} \left( f(\vec{y}\left(\vec{x}\right)) \right)
  =
  \pd{}{\vec{x}}
  \left(
    \left.
      \pd{f}{\vec{y}}(\vec{y})
    \right|_{\vec{y} = \vec{y}(\vec{x})}
    \pd{\vec{y}}{\vec{x}}
  \right)
  =
  \transpose{
    \left(
      \pd{\vec{y}}{\vec{x}}
    \right)
  }
  \left.
    \pd{^2f}{\vec{y} \partial \vec{y}}(\vec{y})
  \right|_{\vec{y} = \vec{y}(\vec{x})}
  \pd{\vec{y}}{\vec{x}}
  +
  \pd{f}{\vec{y}}(\vec{y})
  \pd{^2 \vec{y}}{\vec{x} \partial \vec{x}}
  .
\end{equation}
Now let us choose $f(\vec{y}) = \volienergy(\volentropy, \vec{\xi})$, that is $f(\vec{y}(\vec{x})) = x_1 = \volienergy$. Then
\begin{subequations}
  \label{eq:409}
  \begin{align}
    \label{eq:392}
    \pd{^2 f}{\vec{x} \partial \vec{x}}(\vec{x})
    &=
      \tensorzero,
    \\
    \label{eq:410}
    \pd{^2f}{\vec{y} \partial \vec{y}}(\vec{y})
    &=
      \pd{^2\volienergy}{\vec{y} \partial \vec{y}}(\vec{y})
      =
      \hessian{\volienergy}(\volentropy, \vec{\xi}),
    \\
    \label{eq:411}
    \pd{f}{\vec{y}}(\vec{y})
    &=
    \begin{bmatrix}
      \pd{\volienergy}{\volentropy} \\
      \pd{\volienergy}{\vec{\xi}}
    \end{bmatrix},
    \\
    \label{eq:412}
    \pd{^2 \vec{y}}{\vec{x} \partial \vec{x}}
    &=
      \delta_{\volentropy {\vec{y}}}
      \pd{^2 \volentropy}{\vec{x} \partial \vec{x}}
      =
      \delta_{\volentropy {\vec{y}}}
      \hessian{\volentropy}(\volienergy, \vec{\xi}),
  \end{align}
\end{subequations}
where the last equality follows from the fact that the transformation $\vec{\varphi}$ is nonlinear only in its first component. All mixed derivatives of components $y_i$, $i=2, \dots, n$ thus vanish, and the only second differentiation that survives is the second differentiation of $y_1 \equiv \volentropy$. Substituting~\eqref{eq:409} into~\eqref{eq:354} yields
\begin{equation}
  \label{eq:413}
  \tensorzero
  =
  \transpose{
    \tensorq{J}
  }
  \hessian{\volienergy}(\volentropy, \vec{\xi})
  \tensorq{J}
  +
  \pd{\volienergy}{\volentropy}\hessian{\volentropy}(\volienergy, \vec{\xi}),
\end{equation}
which gives
\begin{equation}
  \label{eq:414}
  \hessian{\volentropy}(\volienergy, \vec{\xi})
  =
  -
  \frac{
    1
  }
  {
    \pd{\volienergy}{\volentropy}
  }
  \transpose{
    \tensorq{J}
  }
  \hessian{\volienergy}(\volentropy, \vec{\xi})
  \tensorq{J},
\end{equation}
which was to be proved.

\subsubsection{Relation between the G\^ateaux derivatives of the modified internal energy $\volienergy(\volentropy, \rho)$ and the modified entropy $\volentropy(\volienergy, \rho)$}
\label{sec:relat-betw-gate}
Using the transformation rule for the Hessian we can find the following relation between the G\^ateaux derivatives of the modified internal energy $\volienergy(\volentropy, \rho)$ and the modified entropy $\volentropy(\volienergy, \rho)$ functions with respect to their natural variables. Let us for the moment introduce the notation $\vec{A} = \left[ \volentropy, \rho \right]$ and $\widehat{\vec{A}} = \left[ \widehat{\volentropy}, \widehat{\rho} \right]$, and we denote the perturbation as
\begin{equation}
  \label{eq:416}
  \widetilde{\vec{A}}
  =_{\bydefinition}
  \vec{A} - \widehat{\vec{A}}
  =
  \left[ \volentropy - \widehat{\volentropy}, \rho - \widehat{\rho} \right]
  =
  \left[ \widetilde{\volentropy}, \widetilde{\rho} \right]
  ,
\end{equation}
and, similarly, for $\vec{B} = \left[ \volienergy, \rho \right]$ and $\widehat{\vec{B}} = \left[ \widehat{\volienergy}, \widehat{\rho} \right]$ we denote
\begin{equation}
  \label{eq:417}
  \widetilde{\vec{B}}
  =_{\bydefinition}
  \vec{B} - \widehat{\vec{B}}
  =
  \left[ \volienergy - \widehat{\volienergy}, \rho - \widehat{\rho} \right]
  =
  \left[ \widetilde{\volienergy}, \widetilde{\rho} \right]
  .
\end{equation}
We have
\begin{equation}
  \label{eq:415}
  \Diff^2 \volienergy (\volentropy, \rho)
  \left[
    \widetilde{\vec{A}}, \widetilde{\vec{A}}
  \right]
  =
  \vectordot{
    \begin{bmatrix}
      \widetilde{\volentropy} \\
      \widetilde{\rho}
    \end{bmatrix}
  }
  {
    \begin{bmatrix}
      \ppd{\volienergy}{\volentropy} & \pd{^2 \volienergy}{\volentropy \partial \rho} \\
      \pd{^2 \volienergy}{\volentropy \partial \rho} & \ppd{\volienergy}{\rho}
    \end{bmatrix}
    \begin{bmatrix}
      \widetilde{\volentropy} \\
      \widetilde{\rho}
    \end{bmatrix}
  }
  =
  \vectordot{
    \begin{bmatrix}
      \widetilde{\volentropy} \\
      \widetilde{\rho}
    \end{bmatrix}
  }
  {
    \hessian{\volienergy}(\volentropy, \rho)
    \begin{bmatrix}
      \widetilde{\volentropy} \\
      \widetilde{\rho}
    \end{bmatrix}
  }
  .
\end{equation}
We observe that the relation between the perturbations in the energy--density and entropy--density representation is
\begin{equation}
  \label{eq:418}
  \widetilde{\volienergy}
  =
  \volienergy - \widehat{\volienergy}
  =
  \volienergy(\widehat{\volentropy} + \widetilde{\volentropy}, \widehat{\rho} + \widetilde{\rho}) - \volienergy(\widehat{\volentropy}, \widehat{\rho})
  \approx
  \pd{\volienergy}{\volentropy} \widetilde{\volentropy} + \pd{\volienergy}{\rho} \widetilde{\rho}.
\end{equation}
We thus have the transformation rule
\begin{equation}
  \label{eq:420}
  \begin{bmatrix}
    \widetilde{\volienergy} \\
    \widetilde{\rho}
  \end{bmatrix}
  =
  \begin{bmatrix}
    \pd{\volienergy}{\volentropy} &  \pd{\volienergy}{\rho} \\
    0 & 1
  \end{bmatrix}
  \begin{bmatrix}
    \widetilde{\volentropy} \\
    \widetilde{\rho}
  \end{bmatrix}
\end{equation}
or vice versa
\begin{equation}
  \label{eq:421}
  \begin{bmatrix}
    \widetilde{\volentropy} \\
    \widetilde{\rho}
  \end{bmatrix}
  =
  \frac{1}{\pd{\volienergy}{\volentropy}}
  \begin{bmatrix}
    1 &  -\pd{\volienergy}{\rho} \\
    0 & \pd{\volienergy}{\volentropy}
  \end{bmatrix}
  \begin{bmatrix}
    \widetilde{\volienergy} \\
    \widetilde{\rho}
  \end{bmatrix}
  =
  \tensorq{J}
  \begin{bmatrix}
    \widetilde{\volienergy} \\
    \widetilde{\rho}
  \end{bmatrix}
  ,
\end{equation}
where we have used the notation from Lemma~\ref{lm:4}. We can now substitute~\eqref{eq:421} into~\eqref{eq:415}, and we get
\begin{multline}
  \label{eq:419}
  \Diff^2 \volienergy (\volentropy, \rho)
  \left[
    \widetilde{\vec{A}}, \widetilde{\vec{A}}
  \right]
  =
  \vectordot{
    \begin{bmatrix}
      \widetilde{\volentropy} \\
      \widetilde{\rho}
    \end{bmatrix}
  }
  {
    \hessian{\volienergy}(\volentropy, \rho)
    \begin{bmatrix}
      \widetilde{\volentropy} \\
      \widetilde{\rho}
    \end{bmatrix}
  }
  =
  \vectordot{
    \tensorq{J}
    \begin{bmatrix}
      \widetilde{\volienergy} \\
      \widetilde{\rho}
    \end{bmatrix}
  }
  {
    \hessian{\volienergy}(\volentropy, \rho)
    \tensorq{J}
    \begin{bmatrix}
      \widetilde{\volienergy} \\
      \widetilde{\rho}
    \end{bmatrix}
  }
  \\
  =
  -
  \pd{\volienergy}{\volentropy}
  \vectordot{
    \begin{bmatrix}
      \widetilde{\volienergy} \\
      \widetilde{\rho}
    \end{bmatrix}
  }
  {
    \left(
      -
      \frac{1}{\pd{\volienergy}{\volentropy}}
      \transpose{\tensorq{J}}
      \hessian{\volienergy}(\volentropy, \rho)
      \tensorq{J}
    \right)
      \begin{bmatrix}
        \widetilde{\volienergy} \\
        \widetilde{\rho}
      \end{bmatrix}
    }
  =
  -
  \pd{\volienergy}{\volentropy}
  \vectordot{
    \begin{bmatrix}
      \widetilde{\volienergy} \\
      \widetilde{\rho}
    \end{bmatrix}
  }
  {
    \hessian{\volentropy}(\volienergy, \rho)
    \begin{bmatrix}
      \widetilde{\volienergy} \\
      \widetilde{\rho}
    \end{bmatrix}
  }
  =
  -
  \temp
  \Diff^2 \volentropy (\volienergy, \rho)
  \left[
    \widetilde{\vec{B}}, \widetilde{\vec{B}}
  \right]
  .
\end{multline}
Consequently, if the modified energy $\volienergy(\volentropy, \rho)$ is a \emph{convex} function, then the modified entropy $\volentropy = \volentropy(\volienergy, \rho)$ is a \emph{concave function}, and~\eqref{eq:419} provides an explicit relation between the corresponding Hessian matrices. We summarise our findings in Summary~\ref{summary:second-derivative-in-various-representations}.

\begin{summary}[Relation between second G\^ateaux derivatives of $\volienergy (\volentropy, \rho)$ and $\volentropy (\volienergy, \rho)$]
  \label{summary:second-derivative-in-various-representations}
  Let $\vec{A} = \left[ \volentropy, \rho \right]$ and $\widehat{\vec{A}} = \left[ \widehat{\volentropy}, \widehat{\rho} \right]$, and let
$
  \widetilde{\vec{A}}
  =_{\bydefinition}
  \vec{A} - \widehat{\vec{A}}
  =
  \left[ \volentropy - \widehat{\volentropy}, \rho - \widehat{\rho} \right]
  =
  \left[ \widetilde{\volentropy}, \widetilde{\rho} \right]
  $.
  (The symbol $\vec{A}$ thus stands for the variables in the entropy--density representation.)
Similarly, let $\vec{B} = \left[ \volienergy, \rho \right]$ and $\widehat{\vec{B}} = \left[ \widehat{\volienergy}, \widehat{\rho} \right]$, and let
$
  \widetilde{\vec{B}}
  =_{\bydefinition}
  \vec{B} - \widehat{\vec{B}}
  =
  \left[ \volienergy - \widehat{\volienergy}, \rho - \widehat{\rho} \right]
  =
  \left[ \widetilde{\volienergy}, \widetilde{\rho} \right]
  $.
  (The symbol $\vec{B}$ thus stands for the variables in the energy--density representation.)
  Further, let~$\tensorq{J}$ denote the matrix
  \begin{equation}
    \tensorq{J}
    =
    \label{eq:427}
    \frac{1}{\pd{\volienergy}{\volentropy}}
    \begin{bmatrix}
      1 &  -\pd{\volienergy}{\rho} \\
      0 & \pd{\volienergy}{\volentropy}
    \end{bmatrix}
    .
  \end{equation}
Then
\begin{subequations}
  \label{eq:424}
  \begin{align}
    \label{eq:426}
    \widetilde{\vec{B}}
    &=
    \tensorq{J}
      \widetilde{\vec{A}}
      ,\\
    \label{eq:425}
    \Diff^2 \volienergy (\volentropy, \rho)
    \left[
    \widetilde{\vec{A}}, \widetilde{\vec{A}}
    \right]
    &=
      -
      \temp
      \Diff^2 \volentropy (\volienergy, \rho)
      \left[
      \widetilde{\vec{B}}, \widetilde{\vec{B}}
      \right]
      .
  \end{align}
\end{subequations}
Consequently, if $\volienergy (\volentropy, \rho)$ is a strictly convex function, then $\volentropy (\volienergy, \rho)$ is a strictly concave function.
\end{summary}

\subsubsection{Leading order approximation in entropy--density representation}
\label{sec:lead-order-appr-1}
Now we apply Lemma~\ref{lm:3} to the Bregman distance/divergence induced by the modified internal energy $\volienergy(\volentropy, \rho)$, see~\eqref{eq:287}. We again use the notation $\vec{A} = \left[ \volentropy, \rho \right]$ and $\widehat{\vec{A}} = \left[ \widehat{\volentropy}, \widehat{\rho} \right]$, and we denote the perturbation as
\begin{equation}
  \label{eq:339}
  \widetilde{\vec{A}}
  =_{\bydefinition}
  \vec{A} - \widehat{\vec{A}}
  =
  \left[ \volentropy - \widehat{\volentropy}, \rho - \widehat{\rho} \right]
  =
  \left[ \widetilde{\volentropy}, \widetilde{\rho} \right]
  .
\end{equation}
Lemma~\ref{lm:3} implies that
\begin{equation}
  \label{eq:340}
  \int_{\Omega}
  \bregmandivergenceop{\volienergy}{\vec{A}}{\widehat{\vec{A}}}
  \,
  \cvolumee
  =
  \frac{1}{2}
  \int_{\Omega}
  \vectordot{\widetilde{\vec{A}}}{\hessian{\volienergy}(\widehat{\vec{A}}) \widetilde{\vec{A}}}
  \,
  \cvolumee
  +
  \cdots
  ,
\end{equation}
where the dots denote the terms of order higher than two in the perturbation, and where the Hessian matrix $\hessian{\volienergy}(\widehat{\vec{A}})$ is the matrix~\eqref{eq:206} evaluated at the spatially homogeneous rest state. In other words, the leading order approximation of $\int_{\Omega} \bregmandivergenceop{\volienergy}{\vec{A}}{\widehat{\vec{A}}} \, \cvolumee$ reads
\begin{equation}
  \label{eq:bregman-divergence-hessian-approximation}
  \int_{\Omega}
  \bregmandivergenceop{\volienergy}{\vec{A}}{\widehat{\vec{A}}}
  \,
  \cvolumee
  \approx
  \frac{1}{2}
  \int_{\Omega}
  \vectordot{\widetilde{\vec{A}}}{\hessian{\volienergy}(\widehat{\vec{A}}) \widetilde{\vec{A}}}
  \,
  \cvolumee
  =
  \frac{1}{2}
  \int_{\Omega}
  \vectordot{
    \begin{bmatrix}
      \widetilde{\volentropy} \\
      \widetilde{\rho}
    \end{bmatrix}
  }
  {
    \begin{bmatrix}
      \ppd{\volienergy}{\volentropy} & \pd{^2 \volienergy}{\volentropy \partial \rho} \\
      \pd{^2 \volienergy}{\volentropy \partial \rho} & \ppd{\volienergy}{\rho}
    \end{bmatrix}
    \begin{bmatrix}
      \widetilde{\volentropy} \\
      \widetilde{\rho}
    \end{bmatrix}
  }
  \,
  \cvolumee
  .  
\end{equation}

We could now evaluate all entries of the Hessian matrix~\eqref{eq:206} one by one, but the quadratic form in~\eqref{eq:340} can be found more efficiently. We recall that the first derivatives of the modified internal energy are given by~\eqref{eq:292}. If we use the formula $\thpressure = \rho^2 \pd{\ienergy}{\rho}(\entropy, \rho)$ for the thermodynamic pressure in~\eqref{eq:294}, we see that
\begin{subequations}
  \label{eq:341}
  \begin{align}
    \label{eq:342}
    \pd{\volienergy}{\volentropy}(\volentropy, \rho)
    &=
      \temp,
    \\
    \label{eq:343}
    \pd{\volienergy}{\rho}(\volentropy, \rho)
    &=
      \ienergy - \temp \entropy + \frac{\thpressure}{\rho}
      =
      \fenergy + \frac{\thpressure}{\rho}
      =
      \gibbs
      ,
  \end{align}
\end{subequations}
where $\gibbs$ denotes the specific Gibbs free energy, and where the temperature $\temp$ and the specific Gibbs free energy $\gibbs$ are interpreted as functions of $\volentropy$ and $\rho$. The Hessian matrix is the derivative of the gradient, hence the quadratic form in~\eqref{eq:340} can be written as
\begin{equation}
  \label{eq:344}
  \vectordot{\widetilde{\vec{A}}}{\hessian{\volienergy}(\widehat{\vec{A}}) \widetilde{\vec{A}}}
  =
  \Diff \temp(\widehat{\vec{A}})[\widetilde{\vec{A}}] \, \widetilde{\volentropy}
  +
  \Diff \gibbs(\widehat{\vec{A}})[\widetilde{\vec{A}}] \, \widetilde{\rho}
  ,
\end{equation}
where $\Diff \temp(\widehat{\vec{A}})[\widetilde{\vec{A}}]$ and $\Diff \gibbs(\widehat{\vec{A}})[\widetilde{\vec{A}}]$ denote the derivatives of the temperature $\temp(\volentropy, \rho)$ and the specific Gibbs free energy $\gibbs(\volentropy, \rho)$ at the point $\widehat{\vec{A}}$ in the direction $\widetilde{\vec{A}}$. Formula~\eqref{eq:344} is the continuum mechanics counterpart of the classical expression $\delta \cttemperature \, \delta \ctentropy - \delta \ctpressure \, \delta \ctvolume$ for the second order variation of the energy, see, for example, \cite{callen.hb:thermodynamics}.

Now we express the quadratic form~\eqref{eq:344} in terms of the perturbations of the temperature and the density, $\temp = \widehat{\temp} + \widetilde{\temp}$ and $\rho = \widehat{\rho} + \widetilde{\rho}$, see~\eqref{eq:41}. The perturbation of the modified entropy $\volentropy = \rho \entropy(\temp, \rho)$ reads
\begin{equation}
  \label{eq:345}
  \widetilde{\volentropy}
  =
  \rho \entropy(\temp, \rho)
  -
  \widehat{\rho} \entropy(\widehat{\temp}, \widehat{\rho})
  =
  \left.
    \pd{\left(\rho \entropy\right)}{\temp}
  \right|_{\left[\temp, \rho\right] = \left[\widehat{\temp}, \widehat{\rho}\right]}
  \widetilde{\temp}
  +
  \left.
    \pd{\left(\rho \entropy\right)}{\rho}
  \right|_{\left[\temp, \rho\right] = \left[\widehat{\temp}, \widehat{\rho}\right]}
  \widetilde{\rho}
  +
  \cdots
  .
\end{equation}
The first partial derivative is given by~\eqref{eq:77}, that is $\pd{\left(\rho \entropy\right)}{\temp} = \rho \frac{\cheatvol}{\temp}$. Concerning the second partial derivative, we use the Maxwell relation
\begin{equation}
  \label{eq:346}
  \pd{\entropy}{\rho}(\temp, \rho)
  =
  -
  \pd{^2 \fenergy}{\temp \partial \rho}(\temp, \rho)
  =
  -
  \frac{1}{\rho^2}
  \pd{\thpressure}{\temp}(\temp, \rho)
  ,
\end{equation}
which follows from~\eqref{eq:205} and~\eqref{eq:5}, and we get $\pd{\left(\rho \entropy\right)}{\rho} = \entropy - \frac{1}{\rho} \pd{\thpressure}{\temp}$. Consequently,
\begin{equation}
  \label{eq:347}
  \widetilde{\volentropy}
  =
  \frac{\widehat{\rho} \cheatvol(\widehat{\temp}, \widehat{\rho})}{\widehat{\temp}}
  \widetilde{\temp}
  +
  \left(
    \entropy(\widehat{\temp}, \widehat{\rho})
    -
    \frac{1}{\widehat{\rho}}
    \pd{\thpressure}{\temp}(\widehat{\temp}, \widehat{\rho})
  \right)
  \widetilde{\rho}
  +
  \cdots
  .
\end{equation}
In particular, since $\cheatvol > 0$, the map $\left[\temp, \rho\right] \mapsto \left[\volentropy, \rho\right]$ is locally invertible, and the chain rule implies that the directional derivatives in~\eqref{eq:344} are, up to higher order terms, the increments of the temperature and the specific Gibbs free energy computed in the temperature--density representation. For the temperature we trivially get $\Diff \temp(\widehat{\vec{A}})[\widetilde{\vec{A}}] = \widetilde{\temp} + \cdots$. For the specific Gibbs free energy $\gibbs(\temp, \rho) = \fenergy(\temp, \rho) + \frac{\thpressure(\temp, \rho)}{\rho}$ we use~\eqref{eq:205} and~\eqref{eq:5}, which yields
\begin{subequations}
  \label{eq:348}
  \begin{align}
    \label{eq:349}
    \pd{\gibbs}{\temp}(\temp, \rho)
    &=
      \pd{\fenergy}{\temp}(\temp, \rho)
      +
      \frac{1}{\rho}
      \pd{\thpressure}{\temp}(\temp, \rho)
      =
      -
      \entropy(\temp, \rho)
      +
      \frac{1}{\rho}
      \pd{\thpressure}{\temp}(\temp, \rho)
      ,
    \\
    \label{eq:350}
    \pd{\gibbs}{\rho}(\temp, \rho)
    &=
      \pd{\fenergy}{\rho}(\temp, \rho)
      +
      \frac{1}{\rho}
      \pd{\thpressure}{\rho}(\temp, \rho)
      -
      \frac{\thpressure(\temp, \rho)}{\rho^2}
      =
      \frac{1}{\rho}
      \pd{\thpressure}{\rho}(\temp, \rho)
      ,
  \end{align}
\end{subequations}
hence
\begin{equation}
  \label{eq:351}
  \Diff \gibbs(\widehat{\vec{A}})[\widetilde{\vec{A}}]
  =
  \left(
    -
    \entropy(\widehat{\temp}, \widehat{\rho})
    +
    \frac{1}{\widehat{\rho}}
    \pd{\thpressure}{\temp}(\widehat{\temp}, \widehat{\rho})
  \right)
  \widetilde{\temp}
  +
  \frac{1}{\widehat{\rho}}
  \pd{\thpressure}{\rho}(\widehat{\temp}, \widehat{\rho})
  \widetilde{\rho}
  +
  \cdots
  .
\end{equation}
Substituting~\eqref{eq:347} and~\eqref{eq:351} into~\eqref{eq:344} we see that the mixed terms $\widetilde{\temp} \widetilde{\rho}$ cancel---this is the Maxwell relation~\eqref{eq:346} at work---and we are left with
\begin{equation}
  \label{eq:352}
  \vectordot{\widetilde{\vec{A}}}{\hessian{\volienergy}(\widehat{\vec{A}}) \widetilde{\vec{A}}}
  =
  \frac{\widehat{\rho} \cheatvol(\widehat{\temp}, \widehat{\rho})}{\widehat{\temp}}
  \widetilde{\temp}^2
  +
  \frac{1}{\widehat{\rho}}
  \pd{\thpressure}{\rho}(\widehat{\temp}, \widehat{\rho})
  \widetilde{\rho}^2
  +
  \cdots
  .
\end{equation}
In other words, the switch from the variables $\volentropy$, $\rho$ to the variables $\temp$, $\rho$ diagonalises the Hessian matrix of the modified internal energy. Formulae~\eqref{eq:340} and~\eqref{eq:352} finally give us the following approximation of the Lyapunov type functional
\begin{equation}
  \label{eq:353}
  \int_{\Omega}
  \bregmandivergenceop{\volienergy}{\vec{A}}{\widehat{\vec{A}}}
  \,
  \cvolumee
  \approx
  \int_{\Omega}
  \left(
    \frac{1}{2}
    \widehat{\rho}
    \frac{\cheatvol \left(\widehat{\temp}, \widehat{\rho}\right)}{\widehat{\temp}}
    \widetilde{\temp}^2
    +
    \frac{1}{2} \frac{1}{\widehat{\rho}} \pd{\thpressure}{\rho}(\widehat{\temp}, \widehat{\rho})
    \widetilde{\rho}^2
  \right)
  \,
  \cvolumee
  ,
\end{equation}
which is the quadratic functional~\eqref{eq:127} we have obtained in Section~\ref{sec:gateaux-derivatives} by a brute force calculation. We also see that the positivity of the coefficients in the quadratic functional is tantamount to the thermodynamic stability conditions~\eqref{eq:128}, which is in agreement with the convexity analysis in Section~\ref{sec:proof-convexity}.

The second derivatives of the kinetic energy term $\frac{\absnorm{\linmomentum}^2}{2 \rho}$ at the rest state are
\begin{subequations}
  \label{eq:357}
  \begin{align}
    \label{eq:358}
    \left.
    \pd{^2}{\linmomentum \partial \linmomentum}
    \left(
    \frac{\absnorm{\linmomentum}^2}{2 \rho}
    \right)
    \right|_{\left[\rho, \linmomentum\right] = \left[\widehat{\rho}, \widehat{\linmomentum}\right]}
    &=
      \frac{1}{\widehat{\rho}} \identity
      ,
    \\
    \label{eq:359}
    \left.
    \pd{^2}{\rho \partial \linmomentum}
    \left(
    \frac{\absnorm{\linmomentum}^2}{2 \rho}
    \right)
    \right|_{\left[\rho, \linmomentum\right] = \left[\widehat{\rho}, \widehat{\linmomentum}\right]}
    &=
      -
      \frac{\widehat{\linmomentum}}{\widehat{\rho}^2}
      =
      \vec{0}
      ,
    \\
    \label{eq:360}
    \left.
    \ppd{}{\rho}
    \left(
    \frac{\absnorm{\linmomentum}^2}{2 \rho}
    \right)
    \right|_{\left[\rho, \linmomentum\right] = \left[\widehat{\rho}, \widehat{\linmomentum}\right]}
    &=
      \frac{\absnorm{\widehat{\linmomentum}}^2}{\widehat{\rho}^3}
      =
      0
      ,
  \end{align}
\end{subequations}
where we have used $\widehat{\linmomentum} = \vec{0}$. The Hessian matrix of $\voltenergy$ at the rest state is thus block diagonal in the sense that the momentum perturbation decouples from the perturbation of the other variables, and Lemma~\ref{lm:3} yields
\begin{equation}
  \label{eq:361}
  \lyapunovfn{\mequilibrium}{\vec{W}}{\widehat{\vec{W}}_{\equilibrium}}
  =
  \int_{\Omega}
  \frac{1}{2}
  \frac{\absnorm{\widetilde{\linmomentum}}^2}{\widehat{\rho}}
  \,
  \cvolumee
  +
  \int_{\Omega}
  \bregmandivergenceop{\volienergy}{\vec{A}}{\widehat{\vec{A}}_{\equilibrium}}
  \,
  \cvolumee
  +
  \cdots
  .
\end{equation}
Since $\widetilde{\linmomentum} = \rho \vec{v} = \widehat{\rho} \widetilde{\vec{v}} + \widetilde{\rho} \widetilde{\vec{v}}$, we have $\frac{\absnorm{\widetilde{\linmomentum}}^2}{\widehat{\rho}} = \widehat{\rho} \absnorm{\widetilde{\vec{v}}}^2 + \cdots$, and using~\eqref{eq:353} we finally get
\begin{equation}
  \label{eq:362}
  \lyapunovfn{\mequilibrium}{\vec{W}}{\widehat{\vec{W}}_{\equilibrium}}
  =
  \frac{1}{2}
  \left(
    \widehat{\rho}
    \int_{\Omega} \absnorm{\widetilde{\vec{v}}} ^2 \, \cvolumee
    +
    \frac{1}{\widehat{\rho}}
    \pd{\thpressure}{\rho}(\widehat{\temp}, \widehat{\rho})
    \int_{\Omega} \widetilde{\rho}^2 \, \cvolumee
    +
    \widehat{\rho}
    \frac{\cheatvol\left(\widehat{\temp}, \widehat{\rho}\right)}{\widehat{\temp}}
    \int_{\Omega} \widetilde{\temp}^2 \, \cvolumee
  \right)
  +
  \cdots
  ,
\end{equation}
where we have written $\widehat{\rho}$ and $\widehat{\temp}$ instead of $\widehat{\rho}_{\equilibrium}$ and $\widehat{\temp}_{\equilibrium}$. The leading order term in~\eqref{eq:362} is \emph{exactly} the quadratic functional whose time derivative we have computed using the linearised governing equations, see~\eqref{eq:57}. Consequently, we again see that the Lyapunov type functional $\lyapunovfn{\mequilibrium}{\vec{W}}{\widehat{\vec{W}}_{\equilibrium}}$ is thus the nonlinear counterpart of the quadratic functional known from the linearised stability analysis.

Furthermore, we can even use Lemma~\ref{lm:3} directly to the representation based on the Bregman distance/divergence induced by the total modified energy $\voltenergy = \voltenergy(\rho, \linmomentum, \volentropy)$. Then the leading order approximation of the Lyapunov type functional has a particularly compact form
\begin{equation}
  \label{eq:401}
  \lyapunovfn{\mequilibrium}{\vec{W}}{\widehat{\vec{W}}_{\equilibrium}}
  =
  \int_{\Omega}
  \bregmandivergenceop{\voltenergy}{\vec{V}}{\widehat{\vec{V}}_{\equilibrium}}
  \,
  \cvolumee
  \approx
  \frac{1}{2}
  \int_{\Omega}
  \vectordot{\widetilde{\vec{V}}}{\hessian{\voltenergy}(\widehat{\vec{V}}_{\equilibrium}) \widetilde{\vec{V}}}
  \,
  \cvolumee
  ,
\end{equation}
which reduces to
\begin{equation}
  \label{eq:422}
  \frac{1}{2}
  \int_{\Omega}
  \vectordot{\widetilde{\vec{V}}}{\hessian{\voltenergy}(\widehat{\vec{V}}_{\equilibrium}) \widetilde{\vec{V}}}
  \,
  \cvolumee
  =
  \frac{1}{2}
  \int_{\Omega}
  \widehat{\rho}
  \absnorm{\widetilde{\vec{v}}}^2
  \,
  \cvolumee
  +
  \frac{1}{2}
  \int_{\Omega}
  \vectordot{\widetilde{\vec{A}}}{\hessian{\volienergy}(\widehat{\vec{A}}) \widetilde{\vec{A}}}
  \,
  \cvolumee
  \\
  =
  \frac{1}{2}
  \int_{\Omega}
  \widehat{\rho}
  \absnorm{\widetilde{\vec{v}}}^2
  \,
  \cvolumee
  +
  \frac{1}{2}
  \int_{\Omega}
  \vectordot{
    \begin{bmatrix}
      \widetilde{\volentropy} \\
      \widetilde{\rho}
    \end{bmatrix}
  }
  {
    \begin{bmatrix}
      \ppd{\volienergy}{\volentropy} & \pd{^2 \volienergy}{\volentropy \partial \rho} \\
      \pd{^2 \volienergy}{\volentropy \partial \rho} & \ppd{\volienergy}{\rho}
    \end{bmatrix}
    \begin{bmatrix}
      \widetilde{\volentropy} \\
      \widetilde{\rho}
    \end{bmatrix}
  }
  \,
  \cvolumee
  ,
\end{equation}
where $\widetilde{\vec{V}} = \vec{V} - \widehat{\vec{V}}_{\equilibrium}$ is the perturbation of the full triple $\vec{V} = \left[ \rho, \linmomentum, \volentropy \right]$, while $\widetilde{\vec{A}}$ is the perturbation of the pair $\vec{A} = \left[ \volentropy, \rho \right]$ of genuine thermodynamic variables.

\subsubsection{Leading order approximation in energy--density representation}
\label{sec:lead-order-appr}

Analogously, we can apply Lemma~\ref{lm:3} to the representation based on the Bregman distance/divergence induced by the negative generalised entropy $- \volentropygen = - \volentropygen(\rho, \linmomentum, \voltenergy)$, see Summary~\ref{summary:lyapunov-bregman-representations}. The leading order approximation of the Lyapunov type functional then reads
\begin{equation}
  \label{eq:447}
  \lyapunovfn{\mequilibrium}{\vec{W}}{\widehat{\vec{W}}_{\equilibrium}}
  =
  \widehat{\temp}_{\equilibrium}
  \int_{\Omega}
  \bregmandivergenceop{- \volentropygen}{\vec{U}}{\widehat{\vec{U}}_{\equilibrium}}
  \,
  \cvolumee
  \approx
  -
  \frac{\widehat{\temp}_{\equilibrium}}{2}
  \int_{\Omega}
  \vectordot{\widetilde{\vec{U}}}{\hessian{\volentropygen}(\widehat{\vec{U}}_{\equilibrium}) \widetilde{\vec{U}}}
  \,
  \cvolumee
  ,
\end{equation}
where $\vec{U} = \left[ \rho, \linmomentum, \voltenergy \right]$ stands for the full triple of variables in the energy--density representation, and $\widetilde{\vec{U}} = \vec{U} - \widehat{\vec{U}}_{\equilibrium}$. The quadratic form has been computed in~\eqref{eq:432}. At the rest state we have $\widehat{\linmomentum} = \vec{0}$, hence the perturbation of the modified internal energy~\eqref{eq:430} reduces to $\widetilde{\volienergy} = \widetilde{\voltenergy}$, and~\eqref{eq:432} yields
\begin{multline}
  \label{eq:448}
  -
  \frac{\widehat{\temp}}{2}
  \int_{\Omega}
  \vectordot{\widetilde{\vec{U}}}{\hessian{\volentropygen}(\widehat{\vec{U}}) \widetilde{\vec{U}}}
  \,
  \cvolumee
  =
  \frac{1}{2}
  \int_{\Omega}
  \frac{\absnorm{\widetilde{\linmomentum}}^2}{\widehat{\rho}}
  \,
  \cvolumee
  -
  \frac{\widehat{\temp}}{2}
  \int_{\Omega}
  \vectordot{\widetilde{\vec{B}}}{\hessian{\volentropy}(\widehat{\vec{B}}) \widetilde{\vec{B}}}
  \,
  \cvolumee
  \\
  =
  \frac{1}{2}
  \int_{\Omega}
  \frac{\absnorm{\widetilde{\linmomentum}}^2}{\widehat{\rho}}
  \,
  \cvolumee
  -
  \frac{\widehat{\temp}}{2}
  \int_{\Omega}
  \vectordot{
    \begin{bmatrix}
      \widetilde{\volienergy} \\
      \widetilde{\rho}
    \end{bmatrix}
  }
  {
    \begin{bmatrix}
      \ppd{\volentropy}{\volienergy} & \pd{^2 \volentropy}{\volienergy \partial \rho} \\
      \pd{^2 \volentropy}{\volienergy \partial \rho} & \ppd{\volentropy}{\rho}
    \end{bmatrix}
    \begin{bmatrix}
      \widetilde{\volienergy} \\
      \widetilde{\rho}
    \end{bmatrix}
  }
  \,
  \cvolumee
  ,
\end{multline}
where $\vec{B} = \left[ \volienergy, \rho \right]$ stands for the pair of genuine thermodynamic variables in the energy--density representation, and where we have again written $\widehat{\rho}$ and $\widehat{\temp}$ instead of $\widehat{\rho}_{\equilibrium}$ and $\widehat{\temp}_{\equilibrium}$. Since $\widetilde{\linmomentum} = \widehat{\rho} \widetilde{\vec{v}} + \cdots$ at the rest state, the first term in~\eqref{eq:448} is the kinetic energy term in~\eqref{eq:422}, while the second term coincides with the second term in~\eqref{eq:422} by virtue of~\eqref{eq:425}. Both representations thus lead to the same leading order approximation, which is in agreement with the pointwise identity~\eqref{eq:440} between the two Bregman distances/divergences.

The results obtained in this section are summarised in Summary~\ref{summary:leading-order-lyapunov}.

\begin{summary}[Leading order approximation of the Lyapunov type functional]
  \label{summary:leading-order-lyapunov}
  Let $\widehat{\vec{V}}_{\equilibrium} = \left[ \widehat{\rho}_{\equilibrium}, \vec{0}, \widehat{\volentropy}_{\equilibrium} \right]$ be the spatially homogeneous equilibrium rest state~\eqref{eq:spatially-homogenoeus-rest-state} written in the variables $\vec{V} = \left[ \rho, \linmomentum, \volentropy \right]$, and let $\widehat{\vec{U}}_{\equilibrium} = \left[ \widehat{\rho}_{\equilibrium}, \vec{0}, \widehat{\voltenergy}_{\equilibrium} \right]$ be the same state written in the variables $\vec{U} = \left[ \rho, \linmomentum, \voltenergy \right]$. Let $\widetilde{\vec{V}} = \vec{V} - \widehat{\vec{V}}_{\equilibrium}$ and $\widetilde{\vec{U}} = \vec{U} - \widehat{\vec{U}}_{\equilibrium}$ denote the perturbations. The leading order term of the Lyapunov type functional $\lyapunovfn{\mequilibrium}{\vec{W}}{\widehat{\vec{W}}_{\equilibrium}}$ is, in virtue of the representations in Summary~\ref{summary:lyapunov-bregman-representations}, the quadratic form generated by the Hessian matrix of the total modified energy $\voltenergy(\vec{V})$, or, equivalently, by the Hessian matrix of the negative generalised entropy $- \volentropygen(\vec{U})$,
  \begin{equation}
    \label{eq:376}
    \lyapunovfn{\mequilibrium}{\vec{W}}{\widehat{\vec{W}}_{\equilibrium}}
    \approx
    \frac{1}{2}
    \int_{\Omega}
    \vectordot{\widetilde{\vec{V}}}{\hessian{\voltenergy}(\widehat{\vec{V}}_{\equilibrium}) \widetilde{\vec{V}}}
    \,
    \cvolumee
    \approx
    -
    \frac{\widehat{\temp}_{\equilibrium}}{2}
    \int_{\Omega}
    \vectordot{\widetilde{\vec{U}}}{\hessian{\volentropygen}(\widehat{\vec{U}}_{\equilibrium}) \widetilde{\vec{U}}}
    \,
    \cvolumee
    .
  \end{equation}
  In terms of the net total energy $\nettenergy(\vec{V}) = \int_{\Omega} \voltenergy(\vec{V}) \, \cvolumee$ and the net entropy $\netentropy(\vec{U}) = \int_{\Omega} \volentropygen(\vec{U}) \, \cvolumee$, the leading order term is given by the second G\^ateaux derivatives of these functionals,
  \begin{equation}
    \label{eq:377}
    \lyapunovfn{\mequilibrium}{\vec{W}}{\widehat{\vec{W}}_{\equilibrium}}
    \approx
    \frac{1}{2}
    \Diff^2 \nettenergy(\widehat{\vec{V}}_{\equilibrium})[\widetilde{\vec{V}}, \widetilde{\vec{V}}]
    \approx
    -
    \frac{\widehat{\temp}_{\equilibrium}}{2}
    \Diff^2 \netentropy(\widehat{\vec{U}}_{\equilibrium})[\widetilde{\vec{U}}, \widetilde{\vec{U}}]
    .
  \end{equation}
  The leading order term is the quadratic functional obtained in the stability analysis of the linearised governing equations, see Summary~\ref{summary:decay-linearised}.
\end{summary}

  \section{Open systems---spatially inhomogeneous steady state}
\label{sec:open-syst-spat}
Having obtained an answer to Question~\ref{q:isolated} (thermodynamically isolated system) we can proceed with the more demanding analysis of the simple thermodynamically open system specified in Question~\ref{q:open}.

\subsection{Affine correction trick}
\label{sec:affine-corr-trick}
Concerning the stability of open systems, we might follow the \emph{affine correction trick}, and we can try to construct a Lyapunov type functional for a \emph{spatially inhomogeneous steady state} out of the known Lyapunov type functional $\lyapunovfn{\mequilibrium}{\vec{W}}{\widehat{\vec{W}}_{\equilibrium}}$ for the \emph{spatially homogeneous equilibrium rest state}.

In other words, we might conjecture that the functional $\lyapunovfn{\nequilibrium}{\vec{W}}{\widehat{\vec{W}}_{\steady}}$ measuring the proximity of the state $\vec{W}$ to the spatially inhomogeneous steady state~$\widehat{\vec{W}}_{\steady}$ might be constructed as 
\begin{equation}
  \label{eq:37}
  \lyapunovfn{\nequilibrium}{\vec{W}}{\widehat{\vec{W}}_{\steady}}
  =_{\bydefinition}
  \lyapunovfn{\mequilibrium}{\vec{W}}{\widehat{\vec{W}}_{\equilibrium}}
  -
  \lyapunovfn{\mequilibrium}{\widehat{\vec{W}}_{\steady}}{\widehat{\vec{W}}_{\equilibrium}}
  -
  \left.
    \Diff[\vec{W}]
    \lyapunovfn{\mequilibrium}{\vec{W}}{\widehat{\vec{W}}_{\equilibrium}}
  \right|_{\vec{W} = \widehat{\vec{W}}_{\steady}}
  \left[
    \vec{W} - \widehat{\vec{W}}_{\steady}
  \right].
\end{equation}
This is essentially the same construction as that behind the Bregman distance/divergence, see Definition~\ref{dfn:bregman-divergence}, provided that~$\lyapunovfn{\mequilibrium}{\cdot}{\widehat{\vec{W}}_{\equilibrium}}$ is a strictly convex function.

\subsection{Suitable choice of variables for affine correction trick}
\label{sec:suit-choice-vari}
The key question is how to describe the state $\vec{W}$, or, in other words, how to choose the variables suitable for this construction. Since we have already seen that the pair $\volentropy \equiv \rho \entropy$ and~$\rho$ is more suitable than the pair $\temp$, $\rho$, we might want to keep these variables in~\eqref{eq:37}. The clear advantage is that~\eqref{eq:37} in this case reduces to the \emph{construction of the (new) Bregman distance/divergence out of the (original) Bregman distance/divergence} as discussed in Section~\ref{sec:conv-bregm-dist}, see especially Lemma~\ref{lm:1}.

Indeed, if we for a moment fix
\begin{equation}
  \label{eq:83}
  \vec{A} = \left[ \volentropy, \rho \right],
\end{equation}
and if we temporarily restrict $\lyapunovfn{\mequilibrium}{\vec{A}}{\widehat{\vec{A}}_{\equilibrium}}$ just to the Bregman part,
\begin{equation}
  \label{eq:85}
  \lyapunovfnup{\mequilibrium}{\Bregman}{\vec{A}}{\widehat{\vec{A}}_{\equilibrium}}
  =_{\bydefinition}
  \int_{\Omega}
   \bregmandivergenceop{\volienergy}{\vec{A}}{\widehat{\vec{A}}_{\equilibrium}}
   \,
   \cvolumee
   ,
\end{equation}
then the \emph{affine correction trick} construction~\eqref{eq:37} or, for that matter, the \emph{Bregman distance/divergence out of Bregman distance/divergence trick} applied to the Bregman part of $\lyapunovfn{\mequilibrium}{\vec{A}}{\widehat{\vec{A}}_{\equilibrium}}$ reads
\begin{multline}
  \label{eq:84}
  \lyapunovfnup{\nequilibrium}{\Bregman}{\vec{A}}{\widehat{\vec{A}}_{\steady}}
  =
  \lyapunovfnup{\mequilibrium}{\Bregman}{\vec{A}}{\widehat{\vec{A}}_{\equilibrium}}
  -
  \lyapunovfnup{\mequilibrium}{\Bregman}{\widehat{\vec{A}}_{\steady}}{\widehat{\vec{A}}_{\equilibrium}}
  \\
  -
  \left.
    \Diff[\vec{A}]
    \lyapunovfnup{\mequilibrium}{\Bregman}{\vec{A}}{\widehat{\vec{A}}_{\equilibrium}}
  \right|_{\vec{A} = \widehat{\vec{A}}_{\steady}}
  \left[
    \vec{A} - \widehat{\vec{A}}_{\steady}
  \right]
  .
\end{multline}
Now we can use the explicit formula~\eqref{eq:85} for $\lyapunovfnup{\mequilibrium}{\Bregman}{\vec{A}}{\widehat{\vec{A}}_{\equilibrium}}$, and we see that~\eqref{eq:84} actually reads
\begin{equation}
  \label{eq:86}
  \lyapunovfnup{\nequilibrium}{\Bregman}{\vec{A}}{\widehat{\vec{A}}_{\steady}}
  =
  \int_{\Omega}
  \left(
    \bregmandivergenceop{\volienergy}{\vec{A}}{\widehat{\vec{A}}_{\equilibrium}}
    -
    \bregmandivergenceop{\volienergy}{\widehat{\vec{A}}_{\steady}}{\widehat{\vec{A}}_{\equilibrium}}
    -
    \left.
      \pd{}{\vec{A}}
      \bregmandivergenceop{\volienergy}{\vec{A}}{\widehat{\vec{A}}_{\equilibrium}}
    \right|_{\vec{A} = \widehat{\vec{A}}_{\steady}}
    \left[
      \vec{A} - \widehat{\vec{A}}_{\steady}
    \right]
  \right)
  \,
  \cvolumee
  ,
\end{equation}
where the integrand is indeed obtained via the \emph{Bregman distance/divergence out of the Bregman distance/divergence} construction discussed in Section~\ref{sec:conv-bregm-dist}. (See the discussion following~\eqref{eq:6} and Lemma~\ref{lm:1}.) Consequently, we have
\begin{equation}
  \label{eq:87}
  \lyapunovfnup{\nequilibrium}{\Bregman}{\vec{A}}{\widehat{\vec{A}}_{\steady}}
  =
  \int_{\Omega}
  \bregmandivergenceop{\volienergy}{\vec{A}}{\widehat{\vec{A}}_{\steady}}
  \,
  \cvolumee
  .
\end{equation}
This observation suggests that we should prefer the variables $\volentropy$ and $\rho$ in all variational constructions.

Now the question is how to apply the \emph{affine correction trick} to the kinetic energy term. It turns out that the suitable variable is not the \emph{velocity} $\vec{v}$ but the \emph{momentum} $\linmomentum$,
\begin{equation}
  \label{eq:241}
  \linmomentum =_{\bydefinition} \rho \vec{v}.
\end{equation}
(This is a well known observation applied typically in the theory of compressible fluids, see~\cite{dostalk.m:nonlinear} for further references.) If we use the momentum, then the kinetic energy term can be rewritten as
\begin{equation}
  \label{eq:88}
  \netkenergy\left(\rho, \linmomentum\right)
  =
  _{\bydefinition}
  \int_{\Omega}
  \frac{1}{2} \frac{\absnorm{\linmomentum}^2}{\rho}
  \,
  \cvolumee,
\end{equation}
and the affine correction trick with respect to the variables $\rho$ and $\linmomentum$ then yields
\begin{equation}
  \label{eq:90}
  \netkenergy\left(\widehat{\rho} + \widetilde{\rho}, \widehat{\linmomentum} + \widetilde{\linmomentum}\right)
  -
  \netkenergy\left(\widehat{\rho}, \widehat{\linmomentum}\right)
  -
  \Diff \netkenergy
  \left(
    \widehat{\rho}, \widehat{\linmomentum}
  \right)
  \left[
    \widetilde{\rho}, \widetilde{\linmomentum}
  \right]
  =
  \int_{\Omega}
  \left[
    \frac{1}{2} \frac{\absnorm{\widehat{\linmomentum} + \widetilde{\linmomentum}}^2}{\widehat{\rho} + \widetilde{\rho}}
    -
    \frac{1}{2} \frac{\absnorm{\widehat{\linmomentum}}^2}{\widehat{\rho}}
    -
    \left(
      \frac{\vectordot{\widehat{\linmomentum}}{\widetilde{\linmomentum}}}{\widehat{\rho}}
      -
      \frac{1}{2} \frac{\absnorm{\widehat{\linmomentum}}^2}{\widehat{\rho}^2} \widetilde{\rho}
    \right)
  \right]
  \,
  \cvolumee
  ,
\end{equation}
which can be, upon using the definition of linear momentum $\widehat{\linmomentum} + \widetilde{\linmomentum} = \left( \widehat{\rho} + \widetilde{\rho} \right) \left(\widehat{\vec{v}} + \widetilde{\vec{v}}\right)$, $\widehat{\linmomentum} = \widehat{\rho} \widehat{\vec{v}}$, rewritten as
\begin{equation}
  \label{eq:112}
  \netkenergy\left(\widehat{\rho} + \widetilde{\rho}, \widehat{\linmomentum} + \widetilde{\linmomentum}\right)
  -
  \netkenergy\left(\widehat{\rho}, \widehat{\linmomentum}\right)
  -
  \Diff \netkenergy
  \left(
    \widehat{\rho}, \widehat{\linmomentum}
  \right)
  \left[
    \widetilde{\rho}, \widetilde{\linmomentum}
  \right]
  =
  \int_{\Omega}
  \frac{1}{2} \left( \widehat{\rho} + \widetilde{\rho} \right)\absnorm{\widetilde{\vec{v}}}^2
  \,
  \cvolumee
  ,
\end{equation}
which is a nonnegative quantity. Consequently, we see that the momentum $\linmomentum$ is indeed the preferred choice of variable. Note however that for incompressible materials, that is for materials with constant density, one can use the \emph{velocity} as well, and the same holds in the setting where $\widehat{\vec{v}}_{\steady} = \vec{0}$, see \cite{dostalk.m.prusa.v.ea:unconditional*1} and \cite{dostalk.m.prusa.v.ea:unconditional} for applications in the incompressible setting.

\subsection{Lyapunov type functional}
\label{sec:lyap-type-funct}

Concerning the construction of a Lyapunov type functional for the open system investigated in Question~\ref{q:open}, we therefore propose to use the functional $\lyapunovfn{\nequilibrium}{\vec{W}}{\widehat{\vec{W}}_{\steady}}$ constructed out of the Lyapunov type functional for the spatially homogeneous rest state $\lyapunovfn{\mequilibrium}{\vec{W}}{\widehat{\vec{W}}_{\equilibrium}}$ by the \emph{affine correction trick}
\begin{equation}
  \label{eq:296}
  \lyapunovfn{\nequilibrium}{\vec{W}}{\widehat{\vec{W}}_{\steady}}
  =_{\bydefinition}
  \lyapunovfn{\mequilibrium}{\vec{W}}{\widehat{\vec{W}}_{\equilibrium}}
  -
  \lyapunovfn{\mequilibrium}{\widehat{\vec{W}}_{\steady}}{\widehat{\vec{W}}_{\equilibrium}}
  -
  \left.
    \Diff[\vec{V}]
    \lyapunovfn{\mequilibrium}{\vec{W}}{\widehat{\vec{W}}_{\equilibrium}}
  \right|_{\vec{V} = \widehat{\vec{V}}_{\steady}}
  \left[
    \vec{V} - \widehat{\vec{V}}_{\steady}
  \right]
  ,
\end{equation}
with respect to the variables
\begin{equation}
  \label{eq:297}
  \vec{V} = \left[ \rho, \linmomentum, \volentropy \right],
\end{equation}
where $\volentropy \equiv \rho \entropy$ and $\linmomentum \equiv \rho \vec{v}$, that is we interpret the state $\vec{W} = \left[ \rho, \vec{v}, \temp \right]$ via the variables $\vec{V}$ when computing the G\^ateaux derivative. (The steady state is then $\widehat{\vec{V}}_{\steady} = \left[ \widehat{\rho}_{\steady}, \vec{0}, \widehat{\volentropy}_{\steady} \right]$.) Using this construction, we arrive at the functional
\begin{equation}
  \label{eq:312}
  \lyapunovfn{\nequilibrium}{\vec{W}}{\widehat{\vec{W}}_{\steady}}
  =
  \int_{\Omega}
  \frac{1}{2} \rho \absnorm{\widetilde{\vec{v}}}^2
  \,
  \cvolumee
  +
  \int_{\Omega}
  \bregmandivergenceop{\volienergy}{\vec{A}}{\widehat{\vec{A}}_{\steady}}
  \,
  \cvolumee
  ,
\end{equation}
where $\vec{A} = \left[ \volentropy, \rho \right]$ and $\widehat{\vec{A}}_{\steady} = \left[ \widehat{\volentropy}_{\steady}, \widehat{\rho}_{\steady} \right]$ are the pairs of the entropy--density variables corresponding to $\vec{W}$ and $\widehat{\vec{W}}_{\steady}$, and where $\widetilde{\vec{v}}$ is the perturbation of the velocity field $\vec{v}$ with respect to the steady velocity field $\widehat{\vec{v}}_{\steady}$, which is in the case of Question~\ref{q:open} the zero velocity field. (Consequently we have $\vec{v} \equiv \widetilde{\vec{v}}$.) If we rewrite the functional in terms of our original variables $\rho$, $\vec{v}$ and $\temp$, we get
\begin{multline}
  \label{eq:325}
  \lyapunovfn{\nequilibrium}{\vec{W}}{\widehat{\vec{W}}_{\steady}}
  =
  \int_{\Omega} \frac{1}{2} \rho \absnorm{\vec{v}}^2 \, \cvolumee
  \\
  +
  \int_{\Omega}
  \left[
    -
    \widehat{\temp}_{\steady} \left\{ \rho \entropy(\temp, \rho) - \widehat{\rho}_{\steady} \entropy(\widehat{\temp}_{\steady}, \widehat{\rho}_{\steady}) \right\}
    +
    \left\{ \rho \ienergy(\temp, \rho) - \widehat{\rho}_{\steady} \ienergy(\widehat{\temp}_{\steady}, \widehat{\rho}_{\steady}) \right\}
  \right]
  \, \cvolumee
  \\
  -
  \int_{\Omega}
  \left[
    \left( \frac{\thpressure(\widehat{\temp}_{\steady}, \widehat{\rho}_{\steady})}{\widehat{\rho}_{\steady}} + \fenergy(\widehat{\temp}_{\steady}, \widehat{\rho}_{\steady}) \right) \left(\rho  - \widehat{\rho}_{\steady} \right)
  \right]
  \, \cvolumee
  .
\end{multline}
This is visually the same functional as in the thermodynamically isolated setting~\eqref{eq:302} in which we formally replace the equilibrium values $\widehat{\left( \cdot \right)}_{\equilibrium}$ with the steady state values $\widehat{\left( \cdot \right)}_{\steady}$. By design the functional has all desired properties---it is nonnegative and it vanishes if and only if the system is at the steady state $\widehat{\vec{W}}_{\steady}$. We note that the functional construction \emph{never} exploited the fact that $\widehat{\vec{W}}_{\steady}$ is a \emph{steady} state; this property is only of interest in evaluating the \emph{time} derivative of the functional.

\subsection{Time derivative of Lyapunov type functional}
\label{sec:time-deriv-lyap-1}
Having constructed the functional we can investigate its time derivative. Unlike in the thermodynamically isolated setting, see Section~\ref{sec:time-deriv-lyap}, \emph{we can not directly exploit the fact that the net entropy is an increasing function}. The Dirichlet boundary condition~\eqref{eq:boundary-conditions-dirichlet-temperature} does not imply that the heat flux vanishes on the boundary, hence we can not exploit equations such as the entropy evolution equation~\eqref{eq:72}. (In fact, we can not even expect energy conservation---we are dealing with an open system.) But let us try to compute the time derivative anyway.

\subsubsection{Spatially homogeneous Dirichlet boundary condition}
\label{sec:spat-homog-dirichl}
Let us first consider the simple case wherein the boundary temperature $\tempbdr$ in~\eqref{eq:31} is \emph{constant in space}. This implies that the steady state $\widehat{\temp}_{\steady}$ and $\widehat{\rho}_{\steady}$ are spatially homogeneous as well.

Next we observe that the system of interest in Question~\ref{q:open} is mechanically isolated, hence it can not exchange mass with its surroundings. This means that the time derivative of the last integral in~\eqref{eq:325} vanishes. Moreover all other terms in~\eqref{eq:325} that are constant in time vanish, and concerning the time differentiation we in fact only need to deal with a few terms,
\begin{equation}
  \label{eq:326}
  \dd{}{t}
  \lyapunovfn{\nequilibrium}{\vec{W}}{\widehat{\vec{W}}_{\steady}}
  =
  \dd{}{t}
  \left(
    \int_{\Omega} \left( \frac{1}{2} \rho \absnorm{\vec{v}}^2  + \rho \ienergy(\temp, \rho) \right) \, \cvolumee
    -
    \int_{\Omega}
    \widehat{\temp}_{\steady} \rho \entropy(\temp, \rho)
    \,
    \cvolumee
  \right)
  .
\end{equation}
The first integral is in fact the net total energy of the system, and the second integral is easy to differentiate. Consequently, we have
\begin{equation}
  \label{eq:327}
  \dd{}{t}
  \lyapunovfn{\nequilibrium}{\vec{W}}{\widehat{\vec{W}}_{\steady}}
  =
  \dd{\nettenergy(\vec{W})}{t}
  -
  \int_{\Omega}
  \rho \widehat{\temp}_{\steady} \dd{\entropy(\temp, \rho)}{t} 
  \,
  \cvolumee
  .
\end{equation}
(Since the steady state is spatially homogeneous, the \emph{material time derivative} $\dd{\widehat{\temp}_{\steady}}{t} = \pd{ \widehat{\temp}_{\steady}}{t} + \vectordot{\vec{v}}{\nabla  \widehat{\temp}_{\steady}} $ vanishes.) The time derivative of the net total energy is easy to find,
\begin{equation}
  \label{eq:328}
  \dd{\nettenergy(\vec{W})}{t}
  =
  \int_{\Omega}
  \left(
    \divergence \left(\cstress \vec{v}\right)
    -
    \divergence \hfluxc
  \right)
  \,
  \cvolumee
  ,
\end{equation}
where $\cstress$ denotes the Cauchy stress tensor and $\hfluxc$ denotes the heat flux, $\hfluxc = - \kappa \nabla \temp$. Using the Stokes theorem and the boundary condition~\eqref{eq:30} we see that the first term (mechanical energy exchange) vanishes, and that the time derivative of the net total energy reads
\begin{equation}
  \label{eq:329}
  \dd{\nettenergy(\vec{W})}{t}
  =
  \int_{\Omega}
  \divergence \left( \kappa \nabla \temp \right)
  \,
  \cvolumee
  .
\end{equation}

Furthermore, the \emph{pointwise} entropy evolution equation reads
\begin{equation}
  \label{eq:330}
  \rho
  \dd{\entropy}{t}
  =
  \frac{
    \tilde{\lambda} \left(\divergence \vec{v}\right)^2
    +
    2 \nu \tensordot{\traceless{\gradsym}}{\traceless{\gradsym}}
   }
   {
     \temp
   }
   +
   \frac{
     \divergence
     \left(
       \kappa
       \nabla \temp
     \right)
   }
   {
     \temp
   }
   .
 \end{equation}
 Substituting~\eqref{eq:329} and \eqref{eq:330} into~\eqref{eq:327} we get
 \begin{equation}
   \label{eq:331}
   \dd{}{t}
   \lyapunovfn{\nequilibrium}{\vec{W}}{\widehat{\vec{W}}_{\steady}}
   =
   \int_{\Omega}
   \left(
     \divergence \left( \kappa \nabla \temp \right)
     -
     \widehat{\temp}_{\steady}
     \frac{
       \divergence
       \left(
         \kappa
         \nabla \temp
       \right)
     }
     {
       \temp
     }
   \right)
   \, \cvolumee
   -
   \int_{\Omega}
   \widehat{\temp}_{\steady}
   \frac{
    \tilde{\lambda} \left(\divergence \vec{v}\right)^2
    +
    2 \nu \tensordot{\traceless{\gradsym}}{\traceless{\gradsym}}
   }
   {
     \temp
   }
   \, \cvolumee
   .
 \end{equation}
 The second integral on the right-hand side is nonpositive, which is what we want. Concerning the first integral on the right-hand side we can proceed as follows
 \begin{equation}
   \label{eq:332}
   \int_{\Omega}
   \left(
     \divergence \left( \kappa \nabla \temp \right)
     -
     \widehat{\temp}_{\steady}
     \frac{
       \divergence
       \left(
         \kappa
         \nabla \temp
       \right)
     }
     {
       \temp
     }
   \right)
   \, \cvolumee
   =
   \int_{\Omega}
   \left(
     \divergence \left( \kappa \nabla \temp \right)
     -
     \divergence
     \left(
       \frac{
         \kappa
         \nabla \temp
       }
       {
         \frac{\temp}{\widehat{\temp}_{\steady}}
       }
     \right)
     +
     \kappa
     \vectordot{
       \nabla \temp
     }
     {
       \nabla
       \left(
         \frac{1}{\frac{\temp}{\widehat{\temp}_{\steady}}}
       \right)
     }
   \right)
   \, \cvolumee.
 \end{equation}
 Now we group the first two terms, we use Stokes theorem and we find that
 \begin{equation}
   \label{eq:333}
   \int_{\Omega}
   \left(
     \divergence \left( \kappa \nabla \temp \right)
     -
     \divergence
     \left(
       \frac{
         \kappa
         \nabla \temp
       }
       {
         \frac{\temp}{\widehat{\temp}_{\steady}}
       }
     \right)
   \right)
   \,
   \cvolumee
   =
   \int_{\partial \Omega}
   \left(
     1
     -
     \frac{\widehat{\temp}_{\steady}}{\temp}
   \right)
   \vectordot{
     \kappa \nabla \temp
   }
   {
     \vec{n}
   }
   \,
   \csurfacees
   =
   0,
 \end{equation}
 where we have used the fact that the temperature fields $\temp$ and $\widehat{\temp}_{\steady}$ coincide on the \emph{boundary} of $\Omega$. Consequently, it remains to work only with the last term in~\eqref{eq:332}. We manipulate the integrand in the last term as follows,
 \begin{multline}
   \label{eq:334}
   \kappa
   \vectordot{
     \nabla \temp
   }
   {
     \nabla
     \left(
       \frac{1}{\frac{\temp}{\widehat{\temp}_{\steady}}}
     \right)
   }
   =
   -
   \kappa
   \frac{
     \vectordot{
       \nabla \temp
     }
     {
       \nabla
       \frac{\temp}{\widehat{\temp}_{\steady}}
     }
   }
   {
     \left(
       \frac{\temp}{\widehat{\temp}_{\steady}}
     \right)^2
   }
   =
   -
   \kappa
   \widehat{\temp}_{\steady}
   \frac{
     \vectordot{
       \frac{\nabla \temp}{\widehat{\temp}_{\steady}}
     }
     {
       \nabla
       \frac{\temp}{\widehat{\temp}_{\steady}}
     }
   }
   {
     \left(
       \frac{\temp}{\widehat{\temp}_{\steady}}
     \right)^2
   }
   =
   -
   \kappa
   \widehat{\temp}_{\steady}
   \frac{
     \vectordot{
       \nabla
       \frac{\temp}{\widehat{\temp}_{\steady}}
     }
     {
       \nabla
       \frac{\temp}{\widehat{\temp}_{\steady}}
     }
   }
   {
     \left(
       \frac{\temp}{\widehat{\temp}_{\steady}}
     \right)^2
   }
   +
   \kappa
   \widehat{\temp}_{\steady} \temp
   \frac{
     \vectordot{
       \nabla \frac{1}{\widehat{\temp}_{\steady}}
     }
     {
       \nabla
       \frac{\temp}{\widehat{\temp}_{\steady}}
     }
   }
   {
     \left(
       \frac{\temp}{\widehat{\temp}_{\steady}}
     \right)^2
   }
   \\
   =
   -
   \kappa
   \widehat{\temp}_{\steady}
   \vectordot{
     \nabla
     \left(
       \ln
       \frac{\temp}{\widehat{\temp}_{\steady}}
     \right)
   }
   {
     \nabla
     \left(
       \ln
       \frac{\temp}{\widehat{\temp}_{\steady}}
     \right)
   }
   -
   \kappa
   \frac{\temp}{\widehat{\temp}_{\steady}}
   \frac{
     \vectordot{
       \nabla \widehat{\temp}_{\steady}
     }
     {
       \nabla
       \frac{\temp}{\widehat{\temp}_{\steady}}
     }
   }
   {
     \left(
       \frac{\temp}{\widehat{\temp}_{\steady}}
     \right)^2
   }
   \\
   =
   -
   \kappa
   \widehat{\temp}_{\steady}
   \vectordot{
     \nabla
     \left(
       \ln
       \frac{\temp}{\widehat{\temp}_{\steady}}
     \right)
   }
   {
     \nabla
     \left(
       \ln
       \frac{\temp}{\widehat{\temp}_{\steady}}
     \right)
   }
   -
   \vectordot{
     \kappa
     \nabla \widehat{\temp}_{\steady}
   }
   {
     \nabla
     \left(
       \ln
       \frac{\temp}{\widehat{\temp}_{\steady}}
     \right)
   }
   \\
   =
   -
   \kappa
   \widehat{\temp}_{\steady}
   \vectordot{
     \nabla
     \left(
       \ln
       \frac{\temp}{\widehat{\temp}_{\steady}}
     \right)
   }
   {
     \nabla
     \left(
       \ln
       \frac{\temp}{\widehat{\temp}_{\steady}}
     \right)
   }
   -
   \divergence
   \left(
     \left(
       \kappa
       \nabla
       \widehat{\temp}_{\steady}
     \right)
     \left(
       \ln
       \frac{\temp}{\widehat{\temp}_{\steady}}
     \right)
   \right)
   +
   \left(
     \divergence
     \left(
       \kappa
       \nabla
       \widehat{\temp}_{\steady}
     \right)
   \right)
   \ln
   \frac{\temp}{\widehat{\temp}_{\steady}}
   .
 \end{multline}
 However, the steady state temperature $\widehat{\temp}_{\steady}$ is the solution to the steady heat equation~\eqref{eq:33}, hence the last term vanishes. The middle term vanishes after the application of Stokes theorem and due to the fact that the steady state temperature field $\widehat{\temp}_{\steady}$  shares the same boundary condition with the temperature field $\temp$. Consequently, we find that
 \begin{equation}
   \label{eq:335}
      \dd{}{t}
   \lyapunovfn{\nequilibrium}{\vec{W}}{\widehat{\vec{W}}_{\steady}}
   =
   -
   \int_{\Omega}
   \kappa
   \widehat{\temp}_{\steady}
   \vectordot{
     \nabla
     \left(
       \ln
       \frac{\temp}{\widehat{\temp}_{\steady}}
     \right)
   }
   {
     \nabla
     \left(
       \ln
       \frac{\temp}{\widehat{\temp}_{\steady}}
     \right)
   }
   \, \cvolumee
   -
   \int_{\Omega}
   \widehat{\temp}_{\steady}
   \frac{
    \tilde{\lambda} \left(\divergence \vec{v}\right)^2
    +
    2 \nu \tensordot{\traceless{\gradsym}}{\traceless{\gradsym}}
   }
   {
     \temp
   }
   \, \cvolumee
   ,
 \end{equation}
 hence the time derivative is clearly nonpositive as required. Our findings are summarised in Summary~\ref{summary:anwswer-to-stability-problem-in-open-system-restricted}. Note that we have formally gained nothing with respect to the isolated system setting in Summary~\ref{summary:anwswer-to-stability-problem-in-isolated-system}. The thermodynamically open setting with the spatially homogeneous boundary temperature is formally the same as the thermodynamically isolated setting with zero heat flux condition.

 However, we might notice that the tedious manipulation in~\eqref{eq:332} does hold even for the spatially inhomogeneous boundary data, which we might exploit later in the more general setting. We formulate our findings as a lemma.

 \begin{lemma}[Identity for the difference between the heat flux and the rescaled entropy flux]
   \label{lm:2}
   Let $\widehat{\temp}_{\steady}$ be a solution to the steady heat equation
   \begin{subequations}
     \label{eq:242}
     \begin{align}
       \label{eq:243}
       \divergence \left( \kappa \nabla \widehat{\temp}_{\steady} \right) &= 0, \\
       \label{eq:244}
      \left. \widehat{\temp}_{\steady} \right|_{\partial \Omega} &= \tempbdr,
    \end{align}
  \end{subequations}
  and let $\temp$ be an arbitrary temperature field in $\Omega$ such that it attains the same boundary data as $\widehat{\temp}_{\steady}$, that is 
  \begin{equation}
    \label{eq:245}
    \left. \temp \right|_{\partial \Omega} = \left. \widehat{\temp}_{\steady} \right|_{\partial \Omega}.  
  \end{equation}
  Then
  \begin{equation}
    \label{eq:246}
       \int_{\Omega}
   \left(
     \divergence \left( \kappa \nabla \temp \right)
     -
     \widehat{\temp}_{\steady}
     \frac{
       \divergence
       \left(
         \kappa
         \nabla \temp
       \right)
     }
     {
       \temp
     }
   \right)
   \, \cvolumee
   =
     -
   \int_{\Omega}
   \kappa
   \widehat{\temp}_{\steady}
   \vectordot{
     \nabla
     \left(
       \ln
       \frac{\temp}{\widehat{\temp}_{\steady}}
     \right)
   }
   {
     \nabla
     \left(
       \ln
       \frac{\temp}{\widehat{\temp}_{\steady}}
     \right)
   }
   \, \cvolumee
   .
  \end{equation}
 \end{lemma}

   \begin{summary}[Answer to a restricted Question~\ref{q:open}---stability of spatially homogeneous steady state in a thermodynamically open system]
    \label{summary:anwswer-to-stability-problem-in-open-system-restricted}
    Let us consider the thermodynamically open system introduced in Question~\ref{q:open}, where the function $\tempbdr$ in the temperature boundary condition
    \begin{equation*}
      \left. \tempsteady \right|_{\partial \Omega} = \tempbdr
    \end{equation*}
    is a \emph{constant}, which implies that the solution $\tempsteady$ of steady heat equation~\eqref{eq:32} is a constant spatially homogeneous temperature field.

    Let $\vec{W} = \left[\rho, \vec{v}, \temp \right]$ be a generic admissible state of the system, and let $\widehat{\vec{W}}_{\steady} = \left[\widehat{\rho}_{\steady}, \vec{0}, \widehat{\temp}_{\steady} \right] \equiv  \left[\rho_{\steady}, \vec{0}, \temp_{\steady} \right] $ denote the steady state~\eqref{eq:spatially-inhomogenoeus-steady-state-dirichlet}. The functional
    \begin{equation}
      \label{eq:235}
      \lyapunovfn{\nequilibrium}{\vec{W}}{\widehat{\vec{W}}_{\steady}}
      =
      \int_{\Omega}
      \frac{1}{2} \rho \absnorm{\vec{v}}^2
      \,
      \cvolumee
      +
      \int_{\Omega}
      \bregmandivergenceop{\volienergy}{\vec{A}}{\widehat{\vec{A}}_{\steady}}
      \,
      \cvolumee
      ,
    \end{equation}
    see~\eqref{eq:312}, where $\vec{A} = \left[ \volentropy, \rho \right]$ and $\widehat{\vec{A}}_{\steady} = \left[ \widehat{\volentropy}_{\steady}, \widehat{\rho}_{\steady} \right]$ denote the entropy--density variables corresponding to the states $\vec{W}$ and $\widehat{\vec{W}}_{\steady}$, has the following properties. The functional is nonnegative
    \begin{subequations}
      \label{eq:236}
      \begin{equation}
        \label{eq:237}
        \lyapunovfn{\nequilibrium}{\vec{W}}{\widehat{\vec{W}}_{\steady}} \geq 0,
      \end{equation}
      it vanishes if and only if the state $\vec{W}$ coincides with the steady state $\widehat{\vec{W}}_{\steady}$,
      \begin{equation}
        \label{eq:238}
        \lyapunovfn{\nequilibrium}{\vec{W}}{\widehat{\vec{W}}_{\steady}} = 0 \Leftrightarrow  \vec{W} = \widehat{\vec{W}}_{\steady},
      \end{equation}
      and the time derivative of the functional is nonpositive,
      \begin{equation}
        \label{eq:239}
        \dd{}{t} \lyapunovfn{\nequilibrium}{\vec{W}}{\widehat{\vec{W}}_{\steady}} =
        -
        \int_{\Omega}
        \kappa
        \widehat{\temp}_{\steady}
        \vectordot{
          \nabla
          \left(
            \ln
            \frac{\temp}{\widehat{\temp}_{\steady}}
          \right)
        }
        {
          \nabla
          \left(
            \ln
            \frac{\temp}{\widehat{\temp}_{\steady}}
          \right)
        }
        \, \cvolumee
        -
        \int_{\Omega}
        \widehat{\temp}_{\steady}
        \frac{
          \tilde{\lambda} \left(\divergence \vec{v}\right)^2
          +
          2 \nu \tensordot{\traceless{\gradsym}}{\traceless{\gradsym}}
        }
        {
          \temp
        }
        \, \cvolumee
        ,
      \end{equation}
      and the time derivative vanishes if and only if the velocity field vanishes and the temperature field coincides with the steady state temperature field,
      \begin{equation}
        \label{eq:240}
        \dd{}{t} \lyapunovfn{\nequilibrium}{\vec{W}}{\widehat{\vec{W}}_{\steady}} = 0
        \Leftrightarrow
        \vec{v} = \vec{0} \text{ and } \temp = \widehat{\temp}_{\steady}.
      \end{equation}
    \end{subequations}
    The functional $\lyapunovfn{\nequilibrium}{\vec{W}}{\widehat{\vec{W}}_{\steady}}$ is thus a Lyapunov type functional.
  \end{summary}

\subsubsection{Spatially inhomogeneous Dirichlet boundary condition}
\label{sec:spat-inhom-dirichl}
This is a more difficult case. For the \emph{incompressible} variant of this problem see~\cite{dostalk.m.prusa.v.ea:unconditional*1} and rigorous treatment in \cite{abbatiello.a.bulcek.m.ea:on*3,abbatiello.a.bulcek.m.ea:on}.

\section{Functional used in Feireisl's work on compressible Navier--Stokes--Fourier equations}
\label{sec:funct-used-feir}

The relative entropy/energy/ballistic free energy functional used in~\cite{feireisl.e.prazak.d:asymptotic} and previous/subsequent texts reads
\begin{equation}
  \label{eq:249}
  {\mathcal E}
  (\left. \rho, \vartheta, \vec{u} \right\| r, \Theta, \vec{U})
  =
  \int_{\Omega}
  \left(
    \frac{1}{2} \rho
    \absnorm{\vec{u} - \vec{U}}^2
    +
    H_\Theta(\rho, \vartheta)
    -
    \pd{H_\Theta(r, \Theta)}{\rho} \left(\rho - r\right)
    -
    H_\Theta(r, \Theta)
  \right)
  \,
  \cvolumee
  ,
\end{equation}
see~\cite[Equation 23]{feireisl.e:relative}, and the thorough discussion in~\cite[Section 1.3]{feireisl.e.novotny.a:mathematics} and~\cite[Section 3.1.1]{feireisl.e.novotny.a:mathematics}. (Note that this functional was originally used in weak-strong uniqueness analysis and in stability theory only with \emph{constant} $\Theta$.) The function
\begin{equation}
  \label{eq:250}
  H_\Theta(\rho, \vartheta) =_{\bydefinition} \rho \left( \ienergy(\vartheta, \rho) - \Theta \entropy(\vartheta, \rho) \right),
\end{equation}
where $\ienergy$ denotes the internal energy and $\entropy$ denotes the entropy, is referred to as the ballistic free energy. (We use $\entropy$ instead of~$s$ which is the original notation in~\cite{feireisl.e:relative}.) We can now work out the formula
$ 
    H_\Theta(\rho, \vartheta)
    -
    \pd{H_\Theta(r, \Theta)}{\rho} \left(\rho - r\right)
    -
    H_\Theta(r, \Theta)
$
according to the definition~\eqref{eq:250}. We get
\begin{multline}
  \label{eq:251}
  H_\Theta(\rho, \vartheta)
  -
  \pd{H_\Theta(r, \Theta)}{\rho} \left(\rho - r\right)
  -
  H_\Theta(r, \Theta)
  =
  \rho \left( \ienergy(\vartheta, \rho) - \Theta \entropy(\vartheta, \rho) \right)
  \\
  -
  \left[
    \left(\ienergy(\Theta, r) - \Theta \entropy(\Theta, r) \right)
    +
    \left.
      \rho \left( \pd{\ienergy(\vartheta, \rho)}{\rho} - \Theta \pd{\entropy(\vartheta, \rho)}{\rho} \right)
    \right|_{\left[\vartheta, \rho\right] = \left[\Theta, r\right]}
  \right]
  \left(\rho - r\right)
  \\
  -
  r \left( \ienergy(\Theta, r) - \Theta \entropy(\Theta, r) \right)
  .
\end{multline}
Using a thermodynamic identity we see that the middle term reduces as
\begin{equation}
  \label{eq:252}
  \left.
    \rho \left( \pd{\ienergy(\vartheta, \rho)}{\rho} - \Theta \pd{\entropy(\vartheta, \rho)}{\rho} \right)
  \right|_{\left[\vartheta, \rho\right] = \left[\Theta, r\right]}
  =
  \frac{
    \thpressure (\Theta, r)
  }
  {
    r
  },
\end{equation}
which yields
\begin{multline}
  \label{eq:253}
  H_\Theta(\rho, \vartheta)
  -
  \pd{H_\Theta(r, \Theta)}{\rho} \left(\rho - r\right)
  -
  H_\Theta(r, \Theta)
  \\
  =
  \rho \left( \ienergy(\vartheta, \rho) - \Theta \entropy(\vartheta, \rho) \right)
  -
  \left(\ienergy(\Theta, r) - \Theta \entropy(\Theta, r) \right)   \left(\rho - r\right)
  -
  \frac{
    \thpressure (\Theta, r)
  }
  {
    r
  }
  \left(\rho - r\right)
  -
  r \left( \ienergy(\Theta, r) - \Theta \entropy(\Theta, r) \right)
  \\
  =
  \rho
  \left[
    \left\{ \ienergy(\vartheta, \rho) - \Theta \entropy(\vartheta, \rho) \right\}
    -
    \left\{ \ienergy(\Theta, r) - \Theta \entropy(\Theta, r) \right\}
  \right]
  -
  \frac{
    \thpressure (\Theta, r)
  }
  {
    r
  }
  \left(\rho - r\right)
  .
\end{multline}
If we now identify
\begin{subequations}
  \label{eq:254}
  \begin{align}
    \label{eq:255}
    \rho &=_{\bydefinition} \widehat{\rho} + \widetilde{\rho}, \\
    \label{eq:256}
    r &=_{\bydefinition} \widehat{\rho}, \\
    \label{eq:257}
    \vartheta &=_{\bydefinition} \widehat{\temp} + \widetilde{\temp}, \\
    \label{eq:258}
    \Theta &=_{\bydefinition} \widehat{\temp}, \\
    \label{eq:259}
    \vec{u} &=_{\bydefinition} \widehat{\vec{v}} + \widetilde{\vec{v}}, \\
    \label{eq:260}
    \vec{U} &=_{\bydefinition} \widehat{\vec{v}},
  \end{align}
\end{subequations}
then we see that~\eqref{eq:249} reduces in virtue of~\eqref{eq:253} to
\begin{multline}
  \label{eq:261}
  {\mathcal E}
  (\left. \rho, \vartheta, \vec{u} \right\| \widehat{\rho}, \widehat{\temp}, \widehat{\vec{v}})
  =
  \int_{\Omega}
  \frac{1}{2} \left( \widehat{\rho} + \widetilde{\rho} \right)
  \absnorm{\widetilde{\vec{v}}}^2
  \,
  \cvolumee
  \\
    +
  \int_{\Omega}
  \left(\widehat{\rho} + \widetilde{\rho} \right)
  \left[
    \left\{
      \ienergy(\widehat{\temp} + \widetilde{\temp}, \widehat{\rho} + \widetilde{\rho})
      -
      \widehat{\temp}
      \entropy(\widehat{\temp} + \widetilde{\temp}, \widehat{\rho} + \widetilde{\rho})
    \right\}
    -
    \left\{
      \ienergy(\widehat{\temp}, \widehat{\rho})
      -
      \widehat{\temp}
      \entropy(\widehat{\temp}, \widehat{\rho})  
    \right\}
  \right]
  \,
  \cvolumee
  -
  \int_{\Omega}
  \frac{
    \thpressure(\widehat{\temp}, \widehat{\rho})
  }
  {
    \widehat{\rho}
  }
  \widetilde{\rho}
  \,
  \cvolumee
  .
\end{multline}
This is exactly the same functional as~\eqref{eq:325}---the integrands coincide pointwise.



\bibliographystyle{chicago}
\bibliography{vit-prusa}

\end{document}